\documentclass[11pt]{amsart}

\usepackage{amsmath,amssymb,amsfonts}

\usepackage{caption2}
\usepackage[]{hyperref}

\usepackage[pdftex]{color,graphicx}
\usepackage{multicol}
\usepackage{verbatim}

\numberwithin{equation}{section}

\setcaptionmargin{1cm}

\newcommand{\be}{\begin{equation}}
\newcommand{\ee}{\end{equation}}

\usepackage{color}
\definecolor{rosso}{cmyk}{0,1,1,0.4}
\definecolor{rossos}{cmyk}{0,1,1,0.55}
\definecolor{rossoc}{cmyk}{0,1,1,0.2}

\begin{document}

\begin{center}
{
 \bf

%
 

Investigation about a statement equivalent to Riemann
Hypothesis (RH) applied to Dirichlet primitive L functions

    }  
\vskip 0.2cm {\large

Giovanni  Lodone %
 \footnote{
giolodone3@gmail.com
}  
.} 
\\[0.2cm]

{\bf

ABSTRACT    } 
\\[0.2cm]

\vskip 0.2cm {\large

We  apply  a known equivalence,  for RH about Riemann $\zeta$   function,  to  Dirichlet L functions with primitive characters. 

\noindent  A key point is the study of positivity condition in phase derivative of primitive L functions. In order to write it in workable form  it is crucial to extend Euler product to critical strip. By exploiting the concept of $y-$friable integers this task is accomplished. Afterward, by use of the ratio limit in PNT,it is shown that this partial derivative expression cannot be negative along critical line. Fulfilling so the statement above. Then, if the mentioned equivalence is true, this implies that no off-critical line zeros are allowed for at least odd-primitive Dirichlet L functions.
  } 
\\[0.2cm]

\vskip 0.1cm

\end{center}

{

MSC-Class  :   11M06, 11M26, 11M99

  } 
\vspace{0.1cm}



     \tableofcontents

\section{ Introduction}

\noindent Using   Dirichlet characters $\chi(n)$, L functions, \cite [p.~249, and 262] {Apostol:1976}, and, specially,   (\cite [p.~3, and 37]  {Davenport1980} (whose conventions, i.e.   $ \chi_0$  as principal, is used throughout) 
   are defined as :
   
   \be \label {DirichletLFunc}
\L(s.\chi)=\sum _{n=1}^\infty  \   \chi(n)/n^s = \prod_{\forall p} \frac{1}{1-\frac{\chi(p)}{p^s}} \quad p \quad prime  \quad  ,  \quad  \quad \quad \Re(s)>1
\ee

\noindent The infinite product is called Euler Product, while $s : \ 0<\Re[s]<1$, and $s : \ \Re[s]=1/2$ are called respectively: critical strip, and, critical line.  The complex variable is $s=1/2+\epsilon+ it$ throughout. GRH statement is: 

\noindent { \it ``all the zeros  { \bf inside} critical strip of the  expressions  (\ref{DirichletLFunc}), or of their analytic continuations,  are 
  on critical-line''}. %

\noindent %
For reasons will be clear below we focus only on $\chi_{primitive}$, defined in  \cite [p.~168 -170]{Apostol:1976}.  
Companion function  like    $\xi(s)$,  (  \cite [p.~62] {Davenport1980} ) for $\zeta(s)$, is $\xi(s,\chi_{primitive} )$, (  \cite [p.~71] {Davenport1980} ) for $L(s,\chi_{primitive})$, below their definition :

     \be \label {XiForLfunc} 
\xi(s)= \Gamma \left(\frac{s+2}{2} \right)(s-1) \frac {\zeta(s)}{\pi^{s/2}} \ \ \ \ ;  \ \ \  \ \xi(s,\chi) = \left( \frac{q}{\pi}   \right)^{\frac{s+\alpha}{2}} \Gamma  \left( \frac{s+\alpha}{2} \right) L(s,\chi)  
 \ee

 \noindent  where $\Gamma(z)$ is defined in  \cite [p.~251]{Apostol:1976},  $ q= congruence \quad modulus $. While    $\alpha=0$ if $\chi_{primitive}(-1) =1 $, even character,  and,   $\alpha=1$ if $\chi_{primitive}(-1) =-1$, odd character. 
Both (\ref {XiForLfunc} ) have the same zeros of the $\zeta(s)$ or $L(s,\chi)$ functions  respectively in the critical strip  \cite{Davenport1980} \cite [p.~16]  {Edwards:1974cz}.
     \noindent       Following   %
     \cite [p.~71, and 79] {Davenport1980} 
      the functional equation for $\xi(s)$ and %
       $\xi(s,\chi)$  are:

        \be  \label {FuncEquForXiL}
       \quad \xi(1-s)=\xi(s)  \ \ \ ; \ \ \
        \xi(1-s,\bar {\chi})=\frac{i^\alpha q^{\frac{1}{2}}}{\tau(\chi)} \xi(s,\chi) \quad with \quad \chi_{primitive} \ \ %
         \quad
        \ee
          
           \noindent  Where $\bar {\chi}$ is the complex conjugate of $ \chi$ and
              \noindent $\tau(\chi)$ is the gauss sum $\tau(\chi) :=  \sum_{m=1}^q \chi(m) e^{2 \pi i m/q}$  
              see \cite [p.~65] {Davenport1980} ,           
        \noindent        or \cite [p.~165]  {Apostol:1976}. %

\noindent While zeros of $\zeta(s)$ have also a symmetry with respect to real axis, the $L(s,\chi)$ functions maintain this symmetry only for real characters $\chi$.

\noindent  There is  a link between, $\zeta(s)$, and, principal characters of congruence modulus  $q$, see  \cite [p.~232] {Apostol:1976}. It is  reported below for easy reading:

\be \label {Z2Princip}
L(s,\chi_{principal})=\zeta(s)\prod_{p \mid q}(1-p^{-s})
\ee

\noindent So  %
RH for $\zeta(s)$  (unproven until now) would prove also  that Dirichlet functions $L(s,\chi)$ with principal characters have same zeros  of $\zeta(s)$ inside critical strip. 
\noindent There is also  a   link between,    $\psi$, primitive  characters, and, character  $\chi$,  neither principal nor primitive  ( see \cite [p.~262]{Apostol:1976} ),  { \it ``inducing $\psi$''} : 

\be \label {Prim2NoPrim}
 L(s,\chi)=L(s,\psi)\prod_{p \mid q}(1- \psi(p) p^{-s}) \quad where \ \ \chi(n)=\chi_{principal} (n)\psi (n) 
\ee

\noindent %
 To prove completely GRH  for all  characters of Dirichlet L functions (\ref{DirichletLFunc}), it is enough  to prove   RH for $\zeta(s)$,(question unsolved till now),and,  to prove RH for  $L(s,\chi_{primitive})$ . Here we focus only on second point. In particular only on odd primitive characters.

\noindent The structure of the paper is:

\begin{itemize}

\item In section 2 we introduce functions useful for study phase behavior of $\xi(s,\chi_{primitive})$

\item  In section \ref{EPInCR} , thank to $y-$smooth integer concept (\cite {De Bruijn1966}),  we  show that $L(s,\chi_{ not \ principal})$ functions can be expressed by Euler Product in critical strip .

\item  In section \ref {ApplEuPr} we apply  previous section  to cmputation of partial derivative of   $L(s,\chi)$. And in section \ref {CONJ}, by the ratio in P.N.T. we prove  
  Theorem 1: $L(s,\chi_{ odd \ primitive})$   fulfills RH at least for  $| t | >T_{Asymp}(\alpha)$, where  $T_{Asymp}(\alpha)$  is defined in  (\ref  {Asympt} ).

\item in section \ref    {TowardGRH}     we  find that the value of $T_{Asymp}(\alpha)$, for primitive characters, is surprisingly about $\pi$. We extend Theorem 1  for  $| t |\le T_{Asymp}(\alpha) $ where it applies.

\item In appendix \ref {GammaDer} we give an unified treatment of derivatives along $t$ and $\epsilon$  in expressions where appears $\angle[\Gamma(z(s))]$ function excluding the $t=0$ neighborhood. The symbol $\angle[ a ]$ means phase of complex number  $a$.%
\item In appendix \ref{ApproxExpr} we compare an useful approximation to (\ref {DerLAngleSuTL}                                                                                                                                                                                                                                                                                                                                                                                                                                                                                                                                                                                                              ).

\item  In appendix \ref {EPinCS} is justified the use of $L(s,\chi_{primitive})$ in critical strip.

\item In appendix \ref {FEEDBACK} are discussed  some observation from readers.

\end{itemize}
%
\section{Angular Momenta and related Lemmas}

The derivative of the phase of $\xi(s)$ for constant $\epsilon$ with respect to $t$ is  :%
\begin{eqnarray}
\ \ \ \ \
\frac{\partial }{\partial t}  \arctan \left(  \frac {\Im \xi(s)}{ \Re \xi(s) } \right) &=&
 \frac{   1}{1+\left(  
 \frac {\Im\xi(s)}{\Re[\xi(s)]} 
  \right)^2}
   \frac{ \frac{\partial  \Im \xi(s) }{\partial t} \Re\xi(s)           -  \frac{\partial  \Re \xi(s) }{\partial t} \Im\xi(s) }{\left( \Re\xi(s) \right)^2}
\label {VarPhase} \\
&=& 
\frac{ \frac{\partial  \Im \xi(s)}{\partial t} \Re\xi(s)  -  \frac{\partial  \Re \xi(s) }{\partial t} \Im\xi(s)         }{ \left(\Re\xi(s)\right)^2 +\left( \Im\xi(s)\right)^2} = \frac{\partial \angle[\xi(s)]}{\partial t} \nonumber
\end{eqnarray}
\noindent The numerator, that determines the sign of (\ref {VarPhase}), can be seen as the angular momentum with respect to the origin of an unitary mass positioned in $(\Re\xi(s) , \Im\xi(s))$  at time $t$ for constant $\epsilon$. 
$$ DEFINITION \ 1 :  \ \ Angular \ Momentum \ for  \  \xi (s,\chi) $$

We can write  also for  $ \xi (s,\chi)$:

  \be \label {AngMomDef}
 \mathcal{L} [ \xi (s,\chi) ] :=
\det  \left(\begin{matrix} \Re[\xi(s,\chi)]  &  \Im[\xi(s,\chi)] \\ \frac{\partial}{\partial t}\Re[\xi(s,\chi)] & \frac{\partial}{\partial t}\Im[\xi(s,\chi)]  \end{matrix}\right)  = | \xi(s,\chi))|^2  \times \frac{\partial \angle[\xi(s,\chi)]}{\partial t}   \ \ : \  \ \chi_{primitive} %
  \ee

\

\subsection{LEMMA 1 }

%
%
%

\noindent%
 If $A(s)$ is a derivable  complex function and $F$ is a complex constant, then:

\be \label{lemma1}
\mathcal{L}[FA(s)] =
 |F|^2 \mathcal{L}[A(s)] 
\ee

\noindent PROOF: $\angle [F  \  A(s) ]=\angle [F ] +\angle [ A (s)]$ so $\frac {\partial \angle[F  \  A(s)]}{\partial t}=\frac {\partial \angle[  A(s)]}{\partial t}$. But $|FA(s)|^2=|F|^2|A(s)|^2$. From  %
(\ref{AngMomDef} ) 
 follows  (\ref{lemma1}).

\subsection{ LEMMA 2
   $ \ \ \ \ \ \ \ \ \ \mathcal{L} [ \xi (s,\chi_{primitive})]_{\Re(s)=\frac{1}{2}}=0$.}
 
  \noindent            Let us apply   Lemma 1 %
  to \ref  {FuncEquForXiL} ( remember $  \chi_{primitive}$).
   For \cite [p.~66] {Davenport1980}  $|\tau(\chi)  |^2= q$, then  $  \left|   \frac{i^\alpha q^{\frac{1}{2}}}{\tau(\chi)}  \right| =1  $. But,  %
    at $\epsilon= 0$ we have:  $1-s= \bar {s} $. Besides, as $\Gamma(\bar{s})=\int_0^\infty x^{\bar{s}-1} e^{-x}dx = \bar{\Gamma}(s)  $ (\cite [p.~8]  {Edwards:1974cz}), $ \left( \frac{q}{\pi}   \right)^{\frac{\bar{s}+\alpha}{2}} =\bar{ \left( \frac{q}{\pi}   \right)^{\frac{s+\alpha}{2}} }  $, and $L(\bar{s}, \bar{\chi})=\bar{L}(s,\chi)$ ( (\ref{DirichletLFunc})) , then for {\ref{XiForLfunc} ), and ( \ref{FuncEquForXiL}):  %

    \be \label {EtaDef}
    \xi(\bar{s},\bar{\chi})   e^{ - \frac{i}{2}\angle\left [  \frac{i^a q^{\frac{1}{2}}}{\tau(\chi)}  \right ]} =\xi(s,\chi)   e^{  \frac{i}{2}\angle\left [  \frac{i^a q^{\frac{1}{2}}}{\tau(\chi)}  \right ]}=\eta(t,\chi_{primitive}) 
    \ \  \in  \Re \ \  , \ \ \ \epsilon=0  \ee
    \noindent because if $\bar{z}= z$ then  $ z \in  \Re$. 
    So, posing $\eta'(t,\chi)=\frac{d\eta(t,\chi)}{d t}$: %
  $$\mathcal{L} [ \xi (t,\epsilon=0,\chi_{primitive})]= \left | e^{ -  \frac{i}{2}\angle\left [  \frac{i^a q^{\frac{1}{2}}}{\tau(\chi)}  \right ]} \right |^2 \mathcal{L}[ \eta(t )]=
  \left | e^{ -  \frac{i}{2}\angle\left [  \frac{i^a q^{\frac{1}{2}}}{\tau(\chi)}  \right ]} \right |^2
  \det  \left(\begin{matrix} \eta(t,\chi)  &  0  \\  \eta'(t,\chi) & 0 \end{matrix}\right)=
  0$$.

$$DEFINITION \ 2 $$

\noindent %
 Let us define ( see (\ref{FuncEquForXiL}) 
):

\be \label {EpsZero}
\eta(\frac{1}{2} +\epsilon +i t , \chi) =
\eta(s,\chi) =   \left( e^{  \frac{i}{2}\angle\left [  \frac{i^\alpha k^{\frac{1}{2}}}{\tau(\chi)}  \right ]} \xi(s,\chi) \right) ; \ for \  \epsilon=0 \quad it \quad is  \quad real \quad \forall \chi_{primitive}
\ee

\noindent Notice that for Lemma 1   $\forall s : \mathcal{L} [ \xi (s,\chi_{primitive})]=\mathcal{L} [ \eta (s,\chi_{primitive})]
$  as $\left | e^{ - \frac{i}{2}\angle\left [  \frac{i^\alpha q^{\frac{1}{2}}}{\tau(\chi)}  \right ]} \right |^2=1 \ \  \forall q$ and  $\forall \chi_{ primitive}$.

\subsection { LEMMA 3}

 \noindent 
  We have:

\be \label {EquCond}\left[  \frac{\partial   \mathcal{L}  [ \xi (s,\chi) ]  )}{\partial \epsilon} \right]_{\epsilon=0}
=
\left[  \frac{\partial   \mathcal{L}  [ \eta (s,\chi) ]  )}{\partial \epsilon} \right]_{\epsilon=0}
=
\eta(t,\chi ) \left[   -\frac{d^2 \eta(t,\chi )}{dt^2}  \right]+  \left[   \frac  {d \eta(t,\chi )}     {dt}        \right]^2  \ \ \; \ \ \ \chi = \chi_{primitive}
\ee
    
          \noindent  PROOF: 
          
           $$\left[  \frac{\partial   \mathcal{L}  [ \eta (s,\chi) ]  )}{\partial \epsilon} \right]=
           \frac{\partial}{\partial \epsilon} det \left(\begin{matrix}
  \Re[\eta (\frac{ 1}{2}+\epsilon +it)] & \Im[\eta (\frac{ 1}{2}+\epsilon +it )]   \\ 
  \frac{\partial \Re[\eta  (\frac{ 1}{2}+\epsilon +it  )]}{\partial t} &  \frac{\partial \Im[\eta(\frac{ 1}{2}+\epsilon +it )]}{\partial t}  \end{matrix}\right)=
           $$
            
           \be \label {rawCalc1}
            \frac{\partial}{\partial \epsilon}
            \left[  \Re[\eta (\frac{ 1}{2}+\epsilon +it)] \ \times \
              \frac{\partial \Im[\eta(\frac{ 1}{2}+\epsilon +it )]}{\partial t} - 
            \frac{\partial \Re[\eta  (\frac{ 1}{2}+\epsilon +it  )]}{\partial t} \ \times \
             \Im[\eta (\frac{ 1}{2}+\epsilon +it )]
             \right] 
           \ee  
           
    \noindent But for Cauchy-Riemann equations  \cite  [p.~19]  {KConrad:2018}, and, Lemma 2:

               \be \label {Cauchy-Riemann} 
  \frac{\partial \Im[ \eta  (\frac{ 1}{2}+\epsilon +it  ) ]}{\partial \epsilon} = -\frac{\partial \Re[ \eta  (\frac{ 1}{2}+\epsilon +it  ) ]}{\partial  t} \ \ ; \ \ \Im[\eta (\frac{ 1}{2} +it )] =0 \ \ ; \ \ \frac{\partial \Im[\eta(\frac{ 1}{2} +it )]}{\partial t}=0= \frac{\partial \Re[\eta(\frac{ 1}{2} +it )]}{\partial \epsilon}
  \ee
            \noindent so from (\ref{Cauchy-Riemann} ) we can equate , (for $\epsilon=0$),    (\ref  {rawCalc1}) to:  %
            
                      $$=  \Re[\eta (\frac{ 1}{2} +it)] \ \times \
              \frac{\partial^2 \Im[\eta(\frac{ 1}{2} +it )]}{ \partial \epsilon \partial t} 
              -
              \frac{\partial \Re[\eta  (\frac{ 1}{2} +it  )]}{\partial t} \ \times \ \frac { \partial
             \Im[\eta (\frac{ 1}{2} +it )]}{\partial \epsilon}=
                      $$   
            
             $$
             \eta(t,\chi ) \left[   -\frac{d^2 \eta(t,\chi )}{dt^2}  \right]+  \left[   \frac  {d \eta(t,\chi )}     {dt}        \right]^2  = [ \eta'(t,\chi)]^2- \eta''(t,\chi) \eta(t,\chi)
             $$

So Lemma 3 is proved.

  \section {Lemmas on   $L(s,\chi_{primitive}) $ phase variations along $t$ computed by Euler product } \label {StabAndAbsValForL}

\noindent If $\chi$ is not principal, as a primitive character, then $L(s,\chi)$ is an entire function (\cite [p.~255]{Apostol:1976} ), and, (\ref {DirichletLFunc})   converges absolutely for $ \Re(s)>1$, { \bf while converges conditionally  
(see \cite [p.~406] {Apostol:1967} 
 )
 for  $ \Re(s)>0$ ( \cite [p.~7] {Hindry:2010}, and, also unifomily \cite [p.~7] {KConrad:2018}). } So we take $\Re[s]>0$  throughout.
 
\noindent  In appendix   \ref  {EPinCS}  it is shown a proof for  that.

\noindent The idea to exploit Euler product also inside critical strip is not new. See \cite  { {LeClair 2015}}. 

\noindent Here we wonder if
, for primitive characters,  it makes sense  to think   also  of Euler product in $0<\Re(s)\le1$.  We are interested in  phase variations  $\Delta \angle[L(s,\chi_{primitive})]$ computed through Euler product.   %
 %

\subsection{ A preliminary analysis of  Euler product  }   \label {EPInCR}

  \noindent Let us consider the expressions whose we know that second one converges in critical strip ( appendix \ref{EPinCS}  ):                     
                        
                        $$EP(s,p_{max}\rightarrow \infty)=\prod_{\forall p} \frac{1}{1-\frac{\chi(p)}{p^s}} \quad p \quad prime  \quad \ \ \; \ \ \ 
                        \ L(s.\chi_{primitive},n_{max}=p_{max}\rightarrow \infty)=\sum _{n=1}^\infty  \   \chi(n)/n^s
                        $$

 \noindent     It is  known that    $\prod_{\forall p} \frac{1}{1-\frac{\chi(p)}{p^s}} \quad p \quad prime  \quad \ $ can be seen as the product of infinite geometrical serie, each one,  with common ratio $ \frac{\chi(p_j)}{p_j^s}$, ( with  $ \left|   \frac{\chi(p_j)}{p_j^s} \right| <1  $ ) i.e. $  \frac{1}{1-\frac{\chi(p)}{p^s}} = \sum_{n=0}^\infty  \left(\frac{\chi(p_j)}{p_j^s} \right)^n$. Let us consider only a finite number of primes  $j_{max}: \  p_1, p_2 . . .p_{max} \ , \ i.e.  j_{max}=\pi(p_{max}) $, and so also a finite numbers of geometrical series to be multiplied. Besides notice that geometrical series can be written by a finite number of terms:

\noindent                     We can multiply $ j=1,...j_{max}=\pi(p_{max}) $ geometric series with infinite terms, but  we can alternatively choose  to consider, in each geometrical series,   the sum beyond  a certain  exponent $\alpha_j =\left \lceil \frac {\ln(p_{max})}{\ln( p_{j })} \right \rceil$ as a whole. i.e.
           $$ \frac{1}{1-\frac{\chi(p_j)}{p_j^s}} = \sum_{n=0}^\infty  \left(\frac{\chi(p_j)}{p_j^s} \right)^n=
           \sum_{n=0}^{\alpha_j-1}  \left(\frac{\chi(p_j)}{p_j^s} \right)^n + \left(\frac{\chi(p_j)}{p_j^s} \right)^{\alpha_j} ( 1-\chi(p_j)p_j^{-s})^{-1}
           $$         
    So we can have a finite number of terms for each of the    $ j=1...j_{max}$ geometric series instead of infinite terms.

 \noindent For example as an $n \approx  p_{max}$ cannot have two divisors both $ >\sqrt{p_{max}}$   then all the  series with $p_j$ from $\approx \sqrt{p_{max}}$ to $p_{max}$ are  simplified as :

   $$1+ \frac{\chi( p_{j } )}    { p_{j}^s}   + \left(  \frac{\chi( p_{j } )}    { p_{j }^s} \right)^2 \frac{1}{1-\frac{\chi( p_{j } )}    { p_{j }^s}}
   \   \  ; \  \sqrt {p_{max}}   < p_j \le   p_{max}   $$
   
   \noindent At the end we get:
   $$EP(t,\epsilon,\chi,p_{max})= $$
       
                       \be \label {EulerTrickFiniteTerms}
                      \prod_{p=2}^{p_{max}}\frac{1}{1-\frac{\chi(p)}{p^s}}
                       =\sum_{n=1}^{p_{max}}\frac{\chi(n)}{n^s} 
                       +
                       \sum_{\rho>p_{max}} \frac{\chi(\rho)}{\rho^s} \prod_{\alpha_j >0} ( 1-\chi(p_j)p_j^{-s})^{-1} = L(s,\chi,p_{max} )+R(s,\chi,p_{max})
                       \ee
                      
                       Where  the integers 
                       \footnote { the greek letter $\rho$ in used for integers  $n> p_{max}$ }

                           \be \label {friableInt} \  \rho= p_1^{\alpha_1} p_2^{\alpha_2} . . .p_{max}^{\alpha_{max}} >p_{max} \ ;  \alpha_j = 0,1,2 . . .  \   \  and \ \  j_{max} = \pi(p_{max})\ \  \ \ \ \ \ \ \ \ \ \ \alpha_j \ ordered \  so \  \rho \ rises \ monothonically
                       \ee

                    \noindent   In (\ref{EulerTrickFiniteTerms}) we can have huge amount of terms :

                     \be \label {NOfTerms} N^{o } (terms)=  3 \prod_{j=1}^{j_{max}-1}   
                        \left(
                         \left  \lceil   \frac{\ln(p_{max})} {\ln(p_j)} \right  \rceil +1
                         \right)
                          %
                        \ \ ; \ \ 
                        last \  factor
                         \left  \lceil   \frac{\ln(p_{max})} {\ln(p_{max})} \right  \rceil +1   \ \ ; \ \  is  \ put \ to \ 3
                          \ee

          
           \noindent  For example for $p_{max}= 31$ we have in  (\ref{NOfTerms} )  that $ N^{o } (terms)= 787320$, whoose only 31 belong to $L(s,\chi,p_{max} )$. The others belong to 
           
           \be \label {Difference}
           R(s,\chi,p_{max})=R(t,\epsilon,\chi,p_{max})= 
            \prod_{p=2}^{p_{max}}\frac{1}{1-\frac{\chi(p)}{p^s}}-
            \sum_{n=1}^{p_{max}}\frac{\chi(n)}{n^s} 
            \ee
           
            \noindent  While $\rho_{max}(p_{max} )$, i.e. the greatest $\rho$  in (\ref{EulerTrickFiniteTerms}), is $>> Primorial(p_{max} )$, so an huge number too, though finite. Besides if we consider a successions of \footnote {$k$ is not a power but only an index.}
           \be \label {Succ} 
           p_{max}^k = \rho_{max}^{k-1}(p_{max}^{k-1})
           \ee
            it is apparent  that all the terms of $R(t,\epsilon,\chi,p_{max}^{k-1})$, are absent in $R(t,\epsilon,\chi,p_{max}^{k})$. So it is apparent that a build-up toward a certain value, from whatever $p'_{max}$ to  $p''_{max} \rightarrow \infty $,  it is to be excluded for  $R(t,\epsilon,\chi,p''_{max})$ when  $p''_{max} \ge \rho_{max}(p'_{max})$.  And we are speaking always of huges but  finite set of terms.

           \noindent So, for finite $p_{max}$,  we do not need to care, in applying Euler rearrangement for  $L(s,\chi,p_{max} )=\sum_{n=1}^{p_{max}}\frac{\chi(n)}{n^s}$, and $R(s,\chi,p_{max})$, about absolute convergence of the   $\prod_{\forall p} \frac{1}{1-\frac{\chi(p)}{p^s}} \quad p \quad prime
                        $, because we are dealing with  {\bf huges, but, finite amount of terms}.
                        
                        \noindent But there are others aspects to explore.


  \subsection {Exploitation of $y-$friable  (or $y-$smooth) number theory}                    
                      
\noindent Here we introduce the concept of integer numbers $x $  that can be factored with primes { \bf all less or equal  then } a given  prime $y$ like  $\rho $ and $p_{max}$ in (\ref  {friableInt}).

  \noindent The integer $n > p_{max}  $ so that $ \forall p$ prime, $>p_{max}  $ we have $ p  \nmid n$, are called $p_{max} -$smooth integers
   \footnote
  {or  $p_{max} -$friable integers. Last name seems more suggestive to convey the image of friable rocks that can be crumbled in smaller stones like $y-$friable integers can be decomposed in prime factors $\le y$. }
   They have been studied extensively, see \cite {Moree:2012}, \cite  {De Bruijn1966}. The standard terminology is that  the  number  of $y -$smooth positive integers till $x$, indicated by $\Psi(y,x)$ is given by $ x \times \rho\left(\frac{\ln(x)}{\ln(y)}\right)$ where   $\rho\left(\frac{\ln(x)}{\ln(y)}\right)$ is  the density 
   \footnote
    {Dicksons - de Bruijn curve }
    given ( \cite {Moree:2012}) by:
  \be \label {Dickman-deBruijnDensity}
  \rho\left(  \frac{\ln(x)}{\ln(y)}  \right) \approx  \left(\frac{\ln(x)}{\ln(y)}\right)^{-\left(\frac{\ln(x)}{\ln(y)}\right)} \ ; \  \ \ \   \frac{\ln(x)}{\ln(y)} \ge 3
  \ee

\noindent   Like for the  simplest case  $2 -$smooth integers (i.e. powers of 2)  they become very sparsely distributed at high $x$ values. If it  is convenient  to link in some way $x$ and $y$ (\cite  {Ravaswami1948}) we can consider $x=y^\delta$.  So ( \ref{Dickman-deBruijnDensity}) becomes $\rho\left(\frac{\ln(x)}{\ln(y)}\right) = \delta^{-\delta} $ i.e. , for example $ \delta=10$ gives  an extremely low density $10^{-10} $  for the $x=y^{10}$ case.  In this case the gap spacing  for $x^{\frac{1}{10}}-$smooth integers  at $x$  must be  not greater that $O(x)$ 
\footnote
 {In reality this bound is slightly less. See \cite [p.~302]  {granville:2008} citing  Tijdeman (1973; 1974).
 The average gap ,with $x>>y$, is $\approx \frac {x}{[\ln(x)]^c}$ where $c$ is a constant  depending from $y$. see better below}. 
 But  consider that 
$\rho\left(\frac{\ln(x)}{\ln(y)}\right) = \delta^{-\delta} =2^{-3}=0.037037... $ in  $x=y^3$ case.

\noindent An interesting fact  is reported in \cite  [p.~2]  {Heath-Brown : 2018} : {\it ``  . . .If $y$ is not too small compared to $x$ it turns out that questions about } [gaps from consecutives] { \it $y-$smooth numbers are no harder, and sometimes easier, that the corresponding questions about primes. For example recent work from  Matomaki and Radziwill  [\cite   [p.~2] {Matomaki 2016}]  shows that for any $\epsilon >0$ } [ here $ \epsilon= \delta^{-1}$ above] { \it  the gaps between consecutive  $x^{\epsilon}-$smooth numbers up to $x$ are at most $O_\epsilon(x^{1/2})$ . . . ``
}

\noindent  this fact  (see \cite   [p.~2] {Matomaki 2016}) happens in an interval  $x, x+C(\epsilon) x^{1/2}$ where $C(\epsilon)$  is a positive constant. 

\noindent But at $x>>y$ (i.e. $y$ fixed) the density of  $y-$smooth numbers evanish.
The  $\Psi(p_{max},x)$ distribution, for $x>> y$  in  (\ref {Dickman-deBruijnDensity}) , when $\rho\left(  \frac{\ln(x)}{\ln(y)}  \right)<<\frac{1}{x}$,  is so low-density, i.e. with gaps not far from the unreachable maximum  $O(x)$.

 \noindent More precisely, the average gap, when $x>>y$, is $\approx \frac {x}{[\ln(x)]^c}$ where $c$ is a constant  depending from $y$. In  (\cite  {Tijdeman:1974}) is shown that the average gap   is $ \asymp_{y} \frac{x} {[\ln(x)]^{c(y)}} $, 
 where $c(y)$ is a constant  dependent from $y$. 
 But $ \frac{x} {[\ln(x)]^{c(y)}}  \gg_{y} x^{1-\delta}$, so we can assume , at $x>>y$,  at least a gap of  $O(x^{1-\delta})$ where $\delta$ is small at will.
The above consideration is usefull to show that, at high $x$,and,  {\bf for constant $p_{max}$ ),} we have that (\ref{smoothSum}  )

\be \label {smoothSum}  
\left|\sum_{p_{max}  <n ,  \ \forall p>p_{max}  : \  p  \  \nmid n }^\infty \frac{\chi(n)} {n^{s}} \right|
\le
 \sum_{p_{max}  <n ,  \ \forall p>p_{max}  : \  p  \  \nmid n }^\infty \frac{|\chi(n)|} {n^{\Re(s)}}
 =\sum_{p_{max}  <n ,  \ \forall p>p_{max}  : \  p  \  \nmid n }^\infty \frac{1} {n^{\Re(s)}}
  \ \ \ ; \ \ \ p_{max}  <\infty
\ee 
is absolutely convergent . 
	$$PROOF$$
	\noindent Suppose  that tha  $y-$friable integer $n_k  $  is so big that  the gap  $ n_k-n_{k-1} \approx n_k^{1-\delta}$. we can write:

	 $$ \frac {1}{n_k^{1/2+\epsilon}}= \frac {n_k-n_{k-1}}{[n_k-n_{k-1}] n_k^{1/2+\epsilon}}
	 =
	 \frac {n_k^ {1-\delta}}
	 {[n_k^{1-\delta}] n_k^{1/2+\epsilon}}
	 <
	  \int_{n_{k-1}}^{n_{k}}\frac {dx}{x ^{3/2 +\epsilon -\delta}} 
	 $$

	\noindent Because for $n_{k-1}<x<n_{k}$, for $\Delta[n_k-n_{k-1}]=\Delta x$, we have 
	
	$$\frac{\Delta[n_k-n_{k-1}]}{[n_k-n_{k-1}] n_k^{1/2+\epsilon}} \ll \frac{\Delta[n_k-n_{k-1}]}{[n_k^{1-\delta}] n_k^{1/2+\epsilon}}  < \frac {\Delta x}{x ^{3/2 +\epsilon -\delta}}$$
	   
\noindent	   So  from a certain $n^*$ upward we have that  for (\ref {smoothSum}  ) we have:
	$$\sum_{p_{max}  <n ,  \ \forall p>p_{max}  : \  p  \  \nmid n }^\infty \frac{1} {n^{\Re(s)}} = 
	O(1)_{ p_{max} }+ O\left(\int_{x=n^*}^{x \rightarrow  \infty}    x^{-(3/2+\epsilon -\delta)}dx \right) \ \ \ , \ \ \ p_{max}<\infty \ \ and \ \ \delta> \approx 0 $$
	\noindent $O(1)_{ p_{max} }$ is constant for $p_{max}$ fixed, but grows as  $p_{max}$ grows.
	

		$$END \ of \ PROOF$$
			
\noindent		As  (\ref{smoothSum}  )  is absolutey convergent , for $\epsilon >-1/2$ (as $ \ \delta> \approx 0 $), then			
	 the order of the sum does not matter,( \cite [p.~411]  {Apostol:1967}) and, we can state that:

	\be \label {OrderNotMatters}  
	\sum_{p_{max}  <n ,  \ \forall p>p_{max}  : \  p  \  \nmid n }^\infty \frac{\chi(n)} {n^{s}} = 
	\sum_{\rho>p_{max}} \frac{\chi(\rho)}{\rho^s} \prod_{\alpha_j >0} \left[ 1-\chi(p_j)p_j^{-s}\right]^{-1}
                   = R(s,\chi,p_{max}) \ \ \ \, \ \ \ \ \  p_{max}  <\infty \ \  \ \ \epsilon> -1/2
	\ee
	
	\noindent The running index is $n=\rho$ , through  infinite values like(\ref{friableInt} ), on the left, and,  $\rho$, through huge but finite values,  on the right. The  absolute convergence of  (\ref {OrderNotMatters}) allows invariant result with a different computation order.  See  \cite [p.~411]  {Apostol:1967}.


 \subsection {Sequences uniformly distributed  modulo 1} \label {SeqUdMod1}

We define,  following  \cite[p.~1]  {Kuipers:1974}, a sequence  $(x_n)_{n=1}^\infty     $
to be equidistributed  mod.1 (modulo 1) if for  whatever sub interval  
$ I= (a,b)
 \subseteq [0,1]$, the proportion of points $(x_n) \in I$  fits asympotically  with $I$ amplitude, i.e.:

\be \label {EquiDistMod1}
\lim_{N\rightarrow  \infty  }    \frac{ \#\{ 1\le n \le N:x_n \in (a,b)\}}{N}=b-a
\ee

It is notewhorthy the following criteria (Weyl 1916) \cite[p.~7]  {Kuipers:1974}:

{\bf Weyl criterion. 

The sequence $(x_n)$, $n=1,2,3. . . .$ is uniformly distributed mod. 1 ,  if and only if,
$$\lim_{N\rightarrow \infty}\frac{1}{N}  \sum_{n=1}^N e^{2 \pi i h  x_n}=0  \ \ \ \ , \ \ \ \forall h \in Z  \setminus 0$$ }

i.e  for all integers $h \ne 0$.
                        \noindent 
                
    \noindent Through Weyl criterion  can be shown that (\cite[p.~8]  {Kuipers:1974}): \label {seqqByIrr}
    
                {  \bf  Let $\theta$ be an irrational number. Then the sequence  $\{ \theta n \}$, $n=1,2,3,....$ is uniformly distributed mod. 1  }            
                
                 \noindent $  \{ \theta n \}$ is called also Kronecker sequence.
          \subsection {Uniform distribution of occurrences of terms with same phase}  \label {UDOPh}

 \noindent The phase of each  term of  (\ref{smoothSum}  ) where $s=(1/2+\epsilon+it)$, and $n=\rho$ is 
:   
 $$\phi_n=-t\ln(n)= -t\ln[p_1^{\alpha_1} p_2^{\alpha_2} . . .p_{max}^{\alpha_{max}} ] = - t \sum_{j=1}^{j_{max}} \alpha_j \ln(p_j)
 $$
 \noindent Below  a Lemma that describes the distribution of the phases  of terms of  (\ref{smoothSum}  ) in $(0,2\pi]$ interval.

{ \bf Lemma  on Euler Product (EP) on critical strip :}  Consider the monothonically rising (in absolute value) sequence $(\phi_k)_{k=1}^\infty     $ with  $\phi_{k}= \phi_n=-t\sum_{j=1}^{j_{max}} \alpha_j \ln(p_j)$ from $n>p_{max}$ by wich we have $\sum_{j=1}^{j_{max}} \alpha_j \ln(p_j)  >\approx\ln(p_{max})$
 \footnote 
 {$>\approx$ means greater but closer as possible
 }.
  Then  $\phi_k=- t \sum_{j=1}^{j_{max}} \alpha_j \ln(p_j) $, 
with $t\ne 0$,   is uniformly distributed mod. $2\pi$ for $p_{max}<n= p_1^{\alpha_1} p_2^{\alpha_2} . . .p_{max}^{\alpha_{max}} \rightarrow \infty \ \ \ $, wich describes the  sequence   of $p_{max}-$smooth integers $> p_{max}$, where,  $$\ \ \  \ ;  \alpha_j = 0,1,2 . . .  \   \  and \ \  j_{max} = \pi(p_{max})\ \  \ \ \ \ \ \ \ \ \ \ \alpha_j \ ordered \  so \  \rho = n \ rises \ monothonically$$.

\noindent We indicate the  rising sequence   of $p_{max}-$smooth integers $> p_{max}$  with  $n_k$ where $k$ is the sequence of natural number. 

$$PROOF $$

 \noindent %
 We have that  $\ln(p_j)$ is irrational \footnote { for $k\ge 2$ we have $\ln(k)$ irrational because if $\ln(k)=\frac{p}{q} \rightarrow k^q=e^{p}$ , but $e$ cannot be a solution of $x^p-k^q=0$ because it is  trascendental.  so    for  
  $k\ge 2$ we have $\ln(k)$ irrational } .
 Suppose that $\ln(2) t$ is irrational too. Then, following \cite[p.~8]  {Kuipers:1974}, see section  \ref{seqqByIrr}
 \footnote{{ \bf  Let $\theta$ be an irrational number. Then the sequence  $\{ \theta n \}$, $n=1,2,3,....$ is uniformly distributed mod. 1  }}, also  $t \alpha_1 \ln(p_1)= t \alpha_1 \ln(2)$,  if $\ln(2)  t $ is irrational, is  uniformly distributed mod. $2 \pi$ as $\alpha_1= 1,2,3 . . .\infty$.
  \footnote { In  \cite[p.~8]  {Kuipers:1974} is considered  mod. 1 . Here we have changed mod. 1 with mod. $2 \pi$. This does not change  results about uniform distributions that with a suitable constant factor  can be translated to mod. 1 }. 
This means that,  \cite[p.~1]  {Kuipers:1974}, if we divide the interval $(0,2 \pi]$ in $M$ subinterval  $\Delta \phi_m = \left(
 \frac{m 2 \pi}{M} -\frac{ \pi}{M} \ , \  \frac{m 2 \pi}{M} + \frac{ \pi}{M}
 \right)$, 
with $m=0,1,2 . . M-1$  we tend o have same occurrences in above subintervals. 

 \noindent More precisely,  for the phase of $n_k^{-s}$, i.e.  $\phi_k$, the occurrence of result of   $\phi_k=- t  \alpha_1 \ln(p_1), \ \alpha_1=1,2,3 . . $        in every  subinterval $\Delta \phi_m$ is the same  as $\alpha_1=1,2, . . . \rightarrow \infty$ because we are dealing with a Kronecker sequence.

\be \label {EquiDistMod2pi}
\lim_{N\rightarrow  \infty  } 
\frac{C_m}{N}=
\lim_{N\rightarrow  \infty  }    \frac{ \#\{ n^*\le n \le N: \left[
\phi_k=\phi_{n-n^*} \in\left(
 \frac{m 2 \pi}{M} -\frac{ \pi}{M} \ , \  \frac{m 2 \pi}{M} + \frac{ \pi}{M}
 \right) \right]_{mod  \ 2\pi} \}}{N}=
 \frac{2 \pi}{M 2 \pi}= \frac{1}{M}
\ee
\noindent where $C_m$ is the counter for occurrences in $\Delta \phi_m$  phase-slot. 

 \noindent In other words if the counter $C_m$ counts    the occurrence of expression   $- t \ \alpha_1 \ln(p_1) $  mod. $2\pi$ in every   $  M$ subinterval of $0 . .2\pi$ interval , i.e. $\Delta \phi_m$, then,  it happens that for   $ -t \alpha_1 \ln(p_1) $ going to $-\infty$,  the $\frac{C_m}{N}$ ratio  tends to be the same  $\forall m \in \{ 0,1,  . . .M-1\}$.
 
 \noindent
Suppose now to use $p_j$ till $j_{max}=\pi(p_{max})$.  Suppose to choose  $\alpha_j$ such to get a monotically raising  $ \sum_{j=1}^{j_{max}} \alpha_j \ln(p_j)$ and $\ln(\rho) >\ln(p_{max})$, so  that    $ \phi_k =- t \sum_{j=1}^{j_{max}} \alpha_j \ln(p_j) $ tends monothonically to $-\infty$, while  the $\alpha_j$ increase orderly, describing monothonically
 \be \label {RhoRisingMonotonically}
 \rho= p_1^{\alpha_1} p_2^{\alpha_2} . . .p_{max}^{\alpha_{max}}  > p_{max} \ \  \ \ \ \ \ \ \ \ \ \ \alpha_j \ ordered \  so \  \rho \ rises \ monothonically
 \ee
 that are the all $p_{max} -$smooth integers. If $\alpha_1\ln(p_1)$ and $\alpha_2\ln(p_2)$  are  equidistributed for $\alpha_1,\alpha_2=1,2,3,...$, then also $\alpha_1\ln(p_1)+\alpha_2\ln(p_2)$ are equidistributed mod $2\pi$. The same is true for a finite sum like $\angle[n_k^{-\Im(s)}]= \phi_k =- t \sum_{j=1}^{j_{max}} \alpha_j \ln(p_j) $.
 If for some $t$ the  value $\ln(p_{j^*})  t $ is rational, then all others  $\ln(p_{j\ne j^*})  t $ are irrational
 \footnote 
 { If a single value of $t=t^*$, by multiplication,   would transform in rational number both  $\ln(p_{j_1})$ and $\ln(p_{j_2})$,  then   $\frac{\ln(p_{j_2})} {\ln(p_{j_2})}= \frac{m}{n}$ 
 with $p_{j_1} \ne p_{j_2}$, and $m,n \in Z \setminus 0$,  this would mean that $\frac{\ln(p_{j_1}) }{ \ln(p_{j_2})}$ is rational ,
 and, this implies that $p_{j_1} ^n=p_{j_2} ^m   $
 but   $p_{j_1} $, and $p_{j_2}$ are different  prime numbers. So it is impossible. 
 }
  and the  situation does not change, i.e. %
 we have an uniform distribution of occurrences of  $\angle[n_k^{-\Im(s)}]= \phi_k =- t \sum_{j=1}^{j_{max}} \alpha_j \ln(p_j) $  mod. $2\pi$ in everyone of the $ \   M$  phase intervals   $\Delta \phi_m$. 
  \noindent Indeed the consideration of characters phases $\angle[\chi(p_j)]$,  in (\ref {EulerTrickFiniteTerms}):
     \be \label {EulerTrickFiniteTerms2}
                     \sum_{p_{max}  <n ,  \ \forall p>p_{max}  : \  p  \  \nmid n }^\infty \frac{\chi(n)} {n^{s}} 
                       = 
                       \sum_{\rho>p_{max}} \frac{\chi(\rho)}{\rho^s} \prod_{\alpha_j >0} ( 1-\chi(p_j)p_j^{-s})^{-1} = R(s,\chi,p_{max})
                       \ee,  does not change the equidistribution of phases in section \ref{UDOPh}, because if  $- \sum_j t \alpha_j\ln(p_j)$ is equidistributed also $-t\ln(\rho)+ \angle[\chi(\rho)]=-\sum_j \ \alpha_j t \ln(p_j)+\angle[\chi(\rho)]$  is equidistributed too \footnote { The sum of not principal charecters is zero while the principal characters is always constant to 1. They give only a rotation to phase vectors that does not affect equidistribution}. 
 
 $$ END \ of \ PROOF$$

 \noindent Notice  that if we allow primes beyond $ p_{max} $ till $p'_{max} >  p_{max} $, then the  phase distribution of sequence   $ - \sum_{j=1}^{j'_{max}>j_{max}} 
  \alpha_j t  \ln(p_j) +\angle[\chi(\rho)]
  $ with $\ln(\rho) >\ln(p_{max})$, is no more  uniform  distributed mod $2\pi$. 

\noindent It becomes uniform mod.  $2\pi$  again after $p'_{max}$, the new $p_{max}$, i.e. the sequence   

$$  -\sum_{j=1}^{j'_{max}}  
 \alpha_j t  \ln(p_j) +\angle[\chi(\rho] 
 \ \  with \ \ \ln \left (\rho \right) =\ln\left( \prod p_j^{\alpha_j}\right) > \ln(p'_{max})$$

          \subsection {Uniform distribution of  the sum of amplitudes  of terms with same phase}   
          
       %
          

 \noindent If the counter $C_m = \#\{ - \sum_{j=1}^{j'_{max}} \left\{ \alpha_j t  \ln(p_j) +\angle[\chi(\rho] \right\} \in  \Delta \phi_m  \} $ counts    the occurrence of expression   $ - \sum_{j=1}^{j'_{max}}
 \alpha_j t  \ln(p_j) +\angle[\chi(\rho] 
 $  mod. $2\pi$ in every   $   M$ subinterval, $\Delta \phi_m $, of $0 . .2\pi$ interval, then  it happens that for   $  \sum_{j=1}^{j_{max}} 
  \alpha_j t  \ln(p_j)
 $ rising monothonically to $\infty$  ( a part the addition of $\angle[\chi(\rho]$ that does not affect equidistribution),  the $C_m$ tend to be the same  $\forall m \in \{ 0,1,  . . .M-1\}$, i.e.  in  each $\Delta \phi_m $ like in (\ref{EquiDistMod2pi})

\noindent The meaning of Weyl criterion  \cite[p.~7]  {Kuipers:1974} 
is that 
in every    $   M$ subinterval, $\Delta \phi_m$, of $0 . .2\pi$ interval,  
the phases of each terms of 
\be \label {smoothSum_0}  
\sum_{p_{max}  <n ,  \ \forall p>p_{max}  : \  p  \  \nmid n }^\infty \frac{\chi(n)} {n^{s}}
\ee 

\noindent, is evenly and homogenuosly distributed mod $2 \pi$,
 {\bf without  concentration in clusters or thinning}.
 
%
  
   \noindent 
   We want to prove that  the sum  of amplitudes ( $\angle[\chi(\rho]$ that does not affect equidistribution is omitted)
   \be \label {AmplSumInEachPhaseSlot-1}
   A_m=\sum_
   {n =\prod_j  p_j^{\alpha_j} \ > p_{max}  \ ,and, \  \  t \sum_{j=1}^{j_{max}} [  \alpha_j \ln(p_j) ] \in \Delta \phi_m}
   a_{n}   =\sum_{  n =\prod_j  p_j^{\alpha_j} \ > p_{max}  \ ,and, \ t \sum_{j=1}^{j_{max}} [  \alpha_j \ln(p_j) ] \in \Delta \phi_m}  
\left|   \frac{\chi(n)} {n^{s}}   \right| 
\ee 
   of occurred terms in each phase slot ($\Delta \phi_m $) %
    tend to be equalized $\forall m$.

\noindent with $M$ whatever big as $p_{max} <\infty$.

\noindent In other words 
: the sum of the amplitudes of addends with phase  $ \in  \Delta \phi_{m_1}$, and, the sum of the amplitudes of addends with phase  $ \in  \Delta \phi_{m_2}$ are such that $A_{m_1}  \approx A_{m_2}$ , and,  as $p_{max}\rightarrow \infty$ , with $M_{p_{max}} $ increasing accordingly  we have  $A_{m}\rightarrow K(p_{max},M_{p_{max}} ) \ , \forall m : 0< m< M-1$ \footnote { where $M_{p_{max}} $ is increased  in order  to contrast the growth of }. 
Let us see why
%
%
%

\begin {comment}

\noindent Let us prove it.
$$PROOF$$

 Weyl criterion affirms:

The sequence $(x_n)$, $n=1,2,3. . . .$ is uniformly distributed mod. 1 ,  iff,
$$\lim_{N\rightarrow \infty}\frac{1}{N}  \sum_{n=1}^N e^{2 \pi i h  x_n}=0  \ \ \ \ , \ \ \ \forall h \in Z  \setminus 0$$ 

i.e  for all integers $h \ne 0$. So, equidistribution happens necessarily not only in the $M$ slots of  interval $(0,1]$ but also in the corresponding  $M\times h^*$ slots of interval   $(0,h^*]$ where $h^*$ is an arbitrary integer.


\noindent Let us consider (in order to easy writing) the phase: $\{   \sum_{j=1}^{j_{max}} \left\{ \alpha_j t  \ln(p_j) +\angle[\chi(\rho] \right\}$.
If we neglect  $+\angle[\chi(\rho] $ the equidistribution is only  rotated  along the trigonometric circle it is not destroyed.


    \noindent For example let us suppose that in phase-slot $\phi_{m_1}$ occur more terms with phase $\approx \phi_{m_1}+ k \times 2\pi$ with odd  $k$  and less  terms with phase $\approx \phi_{m_1}+ k \times 2\pi$ with even  $k$  in such a way that they compensate when we apply Weil criterion with $h^*=1$. So   in the  M intervals of $2\pi$ we have same occurences.  But when we put $h^*=2$  Weil criterion fail in the $2\times 2\pi$ interval  in the corresponding $2 M$ phase slots.  At least  between phase-slots with central phase  $\approx \phi_{m_1}+ k \times 2\pi$  with $k$ odd or even we have for construction a failure of Weyl criterion..

\noindent In other words {\bf Weyl criterion}  requires equidistribution not only  on $[0,1]$ (in our case $2 \pi$) but also on $[0,1\times h^*]$ (in our case in  $2 \pi \times h^*$), with $h^* $ integer $\ne 0$.


\noindent Now let us come to  amplitude  $a_{n-n^*}=n^{-\Re(s)}$. If we take their natural logarithm we get : 

$$\ln(a_{n-n^*})=-\Re(s) \ln(n)=-\Re(s) \ln\left( \prod p_j^{\alpha_j}\right) , \ \ where \ \   \ln\left(\prod p_j^{\alpha_j} \right) >  \ln(p_{max})$$

\noindent So we get that the logarithm of amplitude differs from phase of same term i.e.:
\be \label {PHASE}
 \sum_{j=1}^{j'_{max}} \alpha_j t  \ln(p_j) +\angle[\chi(\rho] 
 \ee
by  $\angle[\chi(\rho]$,  that can be neglected for equidistribution issues. It is only a rigid rotation that does not affect equidistribution. Besides it differs
\noindent  
 by  factor  $\frac{\Re(s)}{t}$. So  we get , 
  the amplitude logarithm  of same term 
 $-\Re(s) \ln(\rho)=-\Re(s)[    \sum_{j=1}^{j_{max}}  \alpha_j   \ln(p_j)] =\ln(a_{n-n^*})$.

\noindent Then we must conclude that  also $\ln(a_{n-n^*})$ is uniformly distributed   mod.  $ \frac{\Re(s) 2 \pi}{t}$. 
So we can apply Weyl criterion  like in  
%
%
 phase-slots 
$\Delta \phi_m = \left(
 \frac{m 2 \pi}{M} -\frac{ \pi}{M} \ , \  \frac{m 2 \pi}{M} + \frac{ \pi}{M}
 \right)$
 to the corresponding log- amplitude slot :
$$\Delta a_\mu = \left(
 \frac{\mu 2 \pi \Re(s)}{M t} -\frac{\Re(s) \pi}{M t} \ , \   \frac{\mu 2 \pi \Re(s)}{M t}  + \frac{\Re(s) \pi}{M t}
 \right)$$ where $\mu= 0,1,2 . . .M-1$

Suppose that the sum of amplitudes of phases  in  $\Delta a_{\mu_1}$  to be different of amplitudes sum of   phases  in in $\Delta a_{\mu_2}$, so we are in a situation  similar to $\phi_{m_1}$ and $\phi_{m_2}$ above.
 
 \noindent This can happen only if in log-amplitude - slot  $\mu_1$, and, in log-amplitude - slot  $\mu_2$ ,  we have same number of  occurrencies  of amplitudes mod.  $ \frac{\Re(s) 2 \pi}{t}$ (in $m$ slots). But we  do not have  same occurrences 
 if we change with mod.  $ \frac{ h^* \times \Re(s) 2 \pi}{t}$  where $h^*\ne0$ is some integer, and  we consider  occurrences in all $h^* \times M$ slot. 
 
 \noindent This  {\bf  could } happens  if we have  {\bf concentration in clusters or thinning}
 
 \noindent But this is  in contradiction with  Weyl criterion that applies to phase-slots 
$\Delta \phi_m = \left(
 \frac{m 2 \pi}{M} -\frac{ \pi}{M} \ , \  \frac{m 2 \pi}{M} + \frac{ \pi}{M}
 \right)$ distribution with mod $2\pi h^*$ and {\bf must} apply also to 
 log-amplitude -slot 
 $$\Delta a_\mu = \left(
 \frac{\mu 2 \pi \Re(s)}{M t} -\frac{\Re(s) \pi}{M t} \ , \   \frac{\mu 2 \pi \Re(s)}{M t}  + \frac{\Re(s) \pi}{M t}
 \right)$$ distribution with mod . $ \frac{\Re(s) 2 \pi h^*}{t}$,  $ \forall h^*\ne 0$. Indeed they are the same distribution  apart a factor for the modulus.


\noindent So  the values of $A_m$  (\ref{AmplSumInEachPhaseSlot-1}) tends to be the same $\forall m$ because the assemblig  of aplitudes $n^{-\Re(s)}$     follows the same evenly and homogenuos distribution (without concentrations or thinning)  $\forall m$. As the phases (mod $2 \pi \times h^*$)  in $\in \Delta \phi_m $ sets, so also log-amplitude (mod.  $ \frac{ h^* \times \Re(s) 2 \pi}{t}$  where $h^*\ne0$ is some integer) . Otherwise Weyl criterion would fail for some $h^*$.

\noindent For $p_{max}<\infty$ the (\ref {smoothSum}  ) is absolutely convergent  to a limit $\lambda (p_{max},t,\epsilon)$. So the sum  in phase and amplitude of terms associated to   $ \Delta \phi_m $  tends to be :

  \be \label {vectorStar}
  \rightarrow
  \approx e^{ i \frac{2 \pi m}{M} } \frac{\lambda(p_{max},t,\epsilon)}{M} \ \ \ , \ \ \ m=0,1,2,3 . . .M-1
  \ee
 
 \noindent  i.e.   the sum of amplitudes tends to be the same, $\forall m$, under same  conditions of even and homogenuous distribution.

$$END \ OF \ PROOF$$

\end {comment}


\noindent Consider (\ref {smoothSum_0}  ). In the monotically increasing sequence of $p_{max}-$friable integers  $n_k=\prod_{j=1}^{j_{max}}p^{ \alpha_j } >  p_{max}$  lets choose a generic term $n_{k^*}$ with amplitude $ (n_{k^*}) ^{-\Re(s)*}$ and phase $-t \ln (n_{k^*}) +\angle[\chi((n_{k^*})] $.

\noindent We have that the neighboring terms  $n_{ k^* +\Delta k }$,  with  whatever $ \Delta k$,   have phase difference  with respect to $n_{k^*}$, i.e. 
$$\{ -t \ln (n_{ k^* +\Delta k }) +\angle[\chi((n_{ k^* +\Delta k })] \} -
\{ -t \ln (n_{k^*}) +\angle[\chi((n_{k^*})] \} =-t\sum_j \beta_j \ln(p_j) +\angle[\chi((n_{ k^* +\Delta k })]-
\angle[\chi((n_{ k^*  })] 
$$

\noindent Where $\beta_j=0,\pm1,\pm2... \ \ \forall j : 1 \le j \le j_{max}$. They distributes uniformly in $2 \pi$ interval as $\angle[\chi((n_{ k^* +\Delta k })]-
\angle[\chi((n_{ k^*  })]$ does not affect equidistribution. 

\noindent While, for the terms amplitude, with $|  \Delta k |$ relatively  small:

\be \label {LastApprox} 
 (n_{k^*}) ^{-\Re(s)} \approx
 (n_{k^* +\Delta k  }) ^{-\Re(s)}
 \ee

 \noindent Of course  (\ref{LastApprox} ) makes sense only  for $n_{k^*}$ and $n_{k^* +\Delta k  }$ not much bigger then $p_{max}$  where the density of $p_{max}-$friable integers is not far from 1.
 But for $p_{max} < \infty$,  but big,  we can say that  practically  all determinant contributions to  (\ref {smoothSum_0}  ) fall in this cathegory.
 Consider also that the relative difference  for unity of $n$, i.e. the derivative of amplitude of the   terms divided by the same term amplitute is : 
 
 \be \label  {DerLog} 
 \frac {  \frac{ d  \left(  
    n^{-\Re(s)}  \right)  } {d n} 
     } 
     { n^{-\Re(s)}  } 
 =\frac {1}{n} 
 \ee
 
 \noindent It goes  to zero as $n$ grows. So we can say that, as $p_{max}$ grows, the interval   $| \Delta k |$ for wich (\ref{LastApprox} ) makes sense  increases too.
 
 \noindent In other words around each $p_{max}-$friable integers, the amplitude of  corrisponding terms of  (\ref {smoothSum_0}  ), distributes almost uniformly in all the $M_{p_{max}}$ phase-slots.  And this statement is better verified as $p_{max}$ is bigger.

 \noindent Another point of view is to consider a particular phase-slot $\Delta \phi_{m^*}$, and to ask  how many terms  between $n_{k_1}$ and $n_{k_2 }$ with $k_2 >> k_1$  fall in $\Delta \phi_{m^*}$.
 As the phase  distributes uniformily  the answer is $ \approx \frac{k_2-k_1}{M} $. But each interval is also an interval  between amplitude  $n_{k_1}^{-\Re(s)}$ and $n_{k_2 }^{-\Re(s)}$. 
 
 \noindent As above we can say that if $n_{k_1}$ and $n_{k_2 }$ are not much greater then $p_{max}$, where the density of $p_{max}-$friable integers is not far from 1, we gather  from $(n_{k_1} \ , \ n_{k_2 })$ interval  in $\Delta \phi_{m^*} \, \ \forall m^*$, a sort of Monte-Carlo integral of the amplitude function  $x^{-1/2-\epsilon}$ apart a constant factor. This is the main contribution for $A_m$ in (\ref {AmplSumInEachPhaseSlot-1}), and implies that the  sum of amplitudes in each-phase slot tends to equalize
  \footnote {The reconstruction of the amplitude by $M^{-1}$ of the original samples  is reasonable because we are  notdealing with phase,   only  with amplitude  of terms whose relative variation tends to zero as $p_{max}\rightarrow \infty$. See  (\ref{DerLog} )}. 
   So (\ref {AmplSumInEachPhaseSlot-1}) tends to be the same $ \forall m^*$ as  $p_{max} \rightarrow \infty$.
  
  \noindent So as the phase is evenly and homogenuosly distributed mod $2 \pi$,  we reach same conclusion starting from a particular phase-slot $\Delta \phi_{m^*}$  as above  starting from a particular $n_{k^*}$ : for $p_{max} \rightarrow \infty$. Namely the sum of amplitudes occurred in each phase-slot where $M_{p_{max}}$ increases  conveniently together with $p_{max}$ must equalize $\forall m :  0 \le m \le M-1$ for $p_{max} \rightarrow \infty$.

  \noindent For $p_{max}<\infty$ the (\ref {smoothSum}  ) is absolutely convergent  to a limit $\lambda (p_{max},t,\epsilon)$ .
  \be \label {AbsSum}\lambda (p_{max},t,\epsilon)=\sum_{p_{max}  <n ,  \ \forall p>p_{max}  : \  p  \  \nmid n }^\infty \frac{| \chi(n)|} {|n^{s}|}
  \ee

\noindent  So the sum  in phase and amplitude of terms associated to  phase-slot   $ \Delta \phi_m $  tends to be :

  \be \label {vectorStar}
  \rightarrow
  \approx e^{ i \frac{2 \pi m}{M} } \frac{\lambda(p_{max},t,\epsilon)}{M_{p_{max}}} \ \ \ , \ \ \ m=0,1,2,3 . . .M_{p_{max}}-1
  \ee
 
 \noindent  because  the sum of amplitudes tends to be the same, $\forall m$ . 
  
   \noindent  (\ref{AbsSum}) grows without bound  for  $p_{max}\rightarrow \infty$, but we can increase $M_{p_{max}}$,  and, we get  more exact  (\ref{vectorStar}) evaluations.
 

\noindent So about (\ref  {AmplSumInEachPhaseSlot-1}) we can affirm also that:

\be \label {AmplSumInEachPhaseSlot}
A_m=
 \sum_{ n =\prod_j  p_j^{\alpha_j} \ > p_{max}  \ ,and, \   t \sum_{j=1}^{j_{max}} [  \alpha_j \ln(p_j)  \in \Delta \phi_m}  
\left|   \frac{\chi(n)} {n^{s}}   \right|  \rightarrow  \approx \frac{ 1     }{M}  \sum_{p_{max}  <n ,  \ \forall p>p_{max}  : \  p  \  \nmid n }^\infty \left|   \frac{\chi(n)} {n^{s}   } \right | 
\ee


\noindent with $p_{max} <\infty$ and $\Re(s)>0$ (i.e. $\epsilon>-1/2$). So with (\ref {AbsSum}) we can write:



   \be \label {Star}
   \sum_{p_{max}  <n ,  \ \forall p>p_{max}  : \  p  \  \nmid n }^\infty   \frac{\chi(n)} {n^{s}} 
   \approx 
   \sum_{m=0}^{m=M-1}  e^{ i \frac{2 \pi m}{M} } \frac{\lambda(p_{max},t,\epsilon)}{M} \ \ \ , \ \ \ 
   \ee
  It is apparent that so  (\ref{smoothSum}  )  can be written also as (\ref{Star}  ). In other words
  $R(s,\chi,p_{max})$ in (\ref {OrderNotMatters} )for $p_{max} \rightarrow \infty$, and, $M_{p_{max}}$ growing accordingly, becomes a star of equispaced and equi-modulus vectors whose sum is zero.
  
   \noindent So  the asymptotic behavior  of Euler product  minus $L(s,\chi_{not \ principal}) $ sum, both truncated, (\ref{Difference})   in critical strip  must tend to $0$ as $p_{max} \rightarrow \infty$ :

 \be \label {AsymBehav}\lim_{p_{max} \rightarrow \infty}       
    R(t,\epsilon,\chi,p_{max})= 
            \prod_{p=2}^{p_{max}}\frac{1}{1-\frac{\chi(p)}{p^s}}-
            \sum_{n=1}^{p_{max}}\frac{\chi(n)}{n^s}
           = 
           0
           \ee

\noindent So   section  \ref {StabAndAbsValForL} aim, i.e. to show that the Euler Product (EP) on critical strip ( at least with not principal characters)   makes sense, is accomplished
\noindent as appendix \ref {EPinCS},  allows  us to  write for $p_{max}\rightarrow \infty$:

\be \label {DirichletLFuncExpanded}
\L(s.\chi,p_{max})=\sum _{n=1}^{p_{max} \rightarrow \infty}  \   \chi(n)/n^s = \prod_{\forall p \le p_{max} \  \ \rightarrow \infty} \ \ \ \  \frac{1}{1-\frac{\chi(p)}{p^s}} \quad \ \ \ \ , \  p \quad prime  \quad  ,  \quad  \chi_{not \ principal}  \quad \quad \Re(s)>0
\ee

   \section{Application of Euler Product  to phase computation in critical strip}                                                 
          \label {ApplEuPr}

                             \noindent Inspired by (\ref{VarPhase}) ,and,  (\ref {WhereEPHolds}) we can consider, for $ s:  \  L(s, \chi) \ne 0$, that: 

\noindent $\angle\left[  \prod_{\forall p} \frac{1}{1-\frac{\chi(p)}{p^s}}  \right]= -\sum_{\forall p}\angle\left[
1-\chi(p) p^{-s}  \right]=   -\sum_{\forall p}\angle\left[
1-\chi(p) p^{-1/2-\epsilon-it}  \right] = $

\noindent $
=
-\sum_{\forall p}\angle\left[
\{ p^{1/2+\epsilon} -\chi(p)  p^{-it} \} p^{-1/2-\epsilon}   \right] =$

\noindent 
$
 -\sum_{\forall p}  \angle\left[    %
 p^{+1/2+\epsilon}  
 - \{   \cos( -\ln(p) t +\angle[\chi(p)] ) + i  \sin(-\ln(p) t +\angle[\chi(p)] )   \}       
\right]=$

$
 -\sum_{\forall p}  \angle\left[ 
\left (      p^{1/2+\epsilon} 
 -    \cos( \ln(p) t -\angle[\chi(p)] ) \right) + i  \left( \sin(\ln(p) t -\angle[\chi(p)] )   \right)         %
\right]
$

\noindent So the phase %
computed from Euler product in  (\ref{DirichletLFunc})  till $p_{max} \rightarrow \infty$ is:
%


\be \label {EulerPrPhaseA_grh}
-
 \sum_{   p<p_{max}    \ : \ gcd(p,q)=1  } \arctan \left( \frac{
 \sin(\ln(p) t    - \angle[\chi(p)]  )}{ p^{1/2+\epsilon}-  
\cos( \ln(p) t   - \angle[\chi(p)]  )  }\right)
\rightarrow
\angle [ %
EP(s,\chi_{primitive},p_{max})
]
 \
\ \ ; \ \ \Re(s) >0
\ee
where  $q$ is the congruence modulus. 

\noindent Far from $(t,\epsilon): \ L(t,\epsilon,\chi)=0$  we can apply (\ref{EulerPrPhaseA_grh}), and,  (\ref{WhereEPHolds}).

\noindent Expression  $\angle [L(s.\chi_{primitive}) ]= \angle \left[\sum _{n=1}^\infty  \   \chi(n)/n^{-1} \right]$,  where total phase is involved, converges for $\Re(s)>0$, but, in singular points.They are, see  (\ref   {VarPhase} ) ,         $s^*_k=1/2+\epsilon^*_k+i t_k^*$,where $\xi (s^*_k,\chi_{primitive})= L(s^*_k,\chi_{primitive})=0$ (see (\ref{XiForLfunc} ), and  (\ref{PhaseVarXiAndL})  ). We neglect, for the moment, double zeros. To be specific $\xi(s,\chi)$ is equal ( \cite [p.~82]  {Davenport1980}) to the infinite Hadamard product, that has single  factors  like \be \label {factor} 
  1-\frac{s}{s^*_k}= 1-\frac{(0.5+\epsilon^*_k)^2+t t^*_k + i (0.5+\epsilon^*_k)(t-t^*_k)}{(0.5+\epsilon^*_k)^2+( t^*_k)^2}=(t^*_k-t)
  \gamma \ ; \ \gamma \in C
  \ee
So, close to the single zero $\ (t_k^*,\epsilon_k^*) \ $, along $\epsilon$ with $-0.5<  \epsilon \le \epsilon^*_k \ge 0$ line, $\angle[\xi(s,\chi)]$  has, at $t= t^*_k$, an abrupt phase jump of $\pi$, due to the factor (\ref {factor}). And in (\ref {XiForLfunc} ), this abrupt  phase change,  can be attributed only to $\angle[L(s,\chi)]$. See  appendix \ref{GammaDer}. But singular points are  also on  half-line $-0.5<\epsilon<\epsilon^*_k \ , \ t=t_k^*$. Here there is a total phase discontinuity of, at least, $2\pi$. Because, as $\delta\rightarrow 0$,  from $(\epsilon<\epsilon_k^* \ , \ t=t^*_k-\delta)$ to $(\epsilon<\epsilon^*  \, \ t=t^*_k+\delta)$ there is a total phase jump of, at least,  $2\pi$, i.e.,  the phase change in circling  one or two zeros couterclockwise ( argument principle \cite  [p.~25]  {KConrad:2018} ). Notice this discontinuity is not present in  (\ref  {VarPhase}  ), where total phase is not used.

  \noindent  We will see in the following we are interested in the difference (or better the phase difference) in a finite interval, suitably far from $\L(s.\chi)=0$. Let us take:
\be \label {Estremi}
t_2=t+ \frac{\pi}{\ln(p^*)} \quad \quad  \mbox{ ; }  \quad \quad t_1=t - \frac{\pi}{\ln(p^*)}
 \mbox{ ; } \ \ \  \Delta t=t_2-t_1= \frac{2\pi}{\ln(p^*)}
\ee

\noindent Let us  see  if for some choice of $p_{max}$ and $p^*$   the difference $R(t_2,\epsilon,p_{max})- R(t_1,\epsilon,p_{max}) $ will be minimized so that

   \be  \label {WhereEPHolds}
       \left.  \Delta \angle [EP(s.\chi,p_{max}<\infty)]   \right|_{\Delta t} 
       \approx
        \left.  \Delta \angle [L(s.\chi,p_{max}<\infty)]   \right|_{\Delta t} \ \ ; 
       \ee

The cancellation will be almost  complete,  for a particular $\rho$, only if  

  \be \label {0Codition}
  \Delta t \times \ln(\rho) = \frac {2 \pi}{\ln(p^*)}  \ln(\rho)   = \mu \times 2 \pi  \ ; \ with \  \mu=  m \in N \ , \ i.e. \  \rho=(p^*)^m
  \ee

\noindent So in order to improve a little an extremely low convergence, pointing to bigger addends of  
, (i.e. at $\rho$ close to $p_{max}$ ),  we can choose:

\be \label {ConvImp}
\ \ p_{max}=(p^*)^m ; \ m \in N \ , \  m \rightarrow \infty
\ee

Thank to (\ref{DirichletLFuncExpanded}),  in section \ref{EPInCR},  we could write  (for big $p_{max}$):

\be \label {PhaseVarInDeltaT2}
\frac{  \left.  \Delta \angle [L(s.\chi),p_{max}]   \right|_{\Delta t} }{\Delta t}  \approx 
\frac{  \left.  \Delta \angle [EP(s.\chi,p_{max})]   \right|_{\Delta t} }{\Delta t}
\ \ \ \ ; \ \  \Delta t=\frac{2 \pi}{\ln(p^*)} \ ;
\ \ p_{max}=(p^*)^m ; \ m \in N \ , \  m \rightarrow \infty
\ee

\noindent See  fig. \ref{Q5-3Camp}, fig. \ref {Simmetries} or similar following  plots.

\noindent Now let us see how to evaluate phase variations of expression  $\angle [L(s.\chi_{primitive}) ]= \angle \left[\sum _{n=1}^\infty  \   \chi(n)/n^{-1} \right]$ using  (\ref    {WhereEPHolds} ).
 \noindent In order to trade resolution along $t$ with computational burden  it is convenient to address to mean phase variation in an interval $(t_1,t_2)$ (\ref {Estremi}) . Besides in computing the mean phase difference in the interval $\Delta t = \frac{2 \pi}{\ln(p^*)}$ (\ref {Estremi}) for $p_{max}\rightarrow \infty$ the oscillations of   $ R(t,\epsilon,p_{max}) $  
 are partially rejected also at $p_{max}=(p^*)^m <\infty$.

\noindent %
We consider,  %
$t$   increment  in order to compute        incremental ratio with respect to $t$, %
 and, define the value $p^*$ like in (\ref{Estremi}).

\noindent For $ \epsilon > - \frac{1}{2}$ (i.e. for $\Re(s)>0$) we can operate this way: { \bf we choose $p^*$ big at will, afterward we compute (\ref {DerLAngleSuTL}) till $p_{max}\rightarrow \infty$, or, better $m \rightarrow \infty$. We can iterate these operations  with a greater  $p^*$}.  %
In shortcut $[ \{ p_{max}=(p^* )^m \rightarrow \infty \}  p^*\rightarrow \infty]$, better  $[ \{p_{max}=(p^*)^m \ , \  m \rightarrow \infty\}  p^*\rightarrow \infty]$.
Notice that we do not interchange  differentiation with infinite sum. %
Instead we compute  the {\bf  mean incremental ratio} of the ``almost everywhere'' differentiable function  
 $\angle [L(s.\chi_{primitive}) ]$ using  (\ref  {WhereEPHolds})    between $t_1$ and $t_2$; i.e. $ \frac {\Delta \angle [ L_{EP}(s,\chi_{primitive})] }{\Delta t} $ where $\Delta t = \frac { 2 \pi}{\ln(p^*}$ using all primes. %
Afterward  every iteration  to  let  $p^*\rightarrow \infty$, requires new computation of mean incremental ratio.

 \be \label {DerLAngleSuTL}                                                                                                                                                                                                                                                                                                                                                                                                                                                                                                                                                                                                              \left[
  - \frac{\ln(p^*) }{2 \pi}
  \sum_{p=2}^{  p=p_{max} =(p^*)^m \ , \ m \rightarrow \infty  \ : \ gcd(p,q)=1  } \left\{    \arctan \left( \frac{
 \sin(\ln(p) t    - \angle[\chi(p)]  )}{ p^{1/2+\epsilon}-  
\cos( \ln(p) t - \angle[\chi(p)]  )  }\right)       \right \}_{t-\frac{\pi}{\ln(p^*)}   }^{t+\frac{\pi}{\ln(p^*)}  }
\right]_{p^* \rightarrow \infty}
 =
 \ee
 $$
 \left\{ \left[  \frac {\Delta \angle [ L_{EP}(s,\chi_{primitive})] }{\Delta t}  \right]_{ p_{max} =(p^*)^m \ , \ m \rightarrow \infty}  \right\}_{  p^*  \rightarrow \infty}
 \rightarrow
 \ \ \ 
 \frac {\partial \angle [ L(s,\chi_{primitive})] }{\partial t} 
 \ \ ; \ \ \{\Re(s) >0\}  \setminus \{sing. points\}
 $$

\noindent %
 It happens that , for $[ \{p_{max} \rightarrow \infty\}  p^*\rightarrow \infty]$, (\ref {DerLAngleSuTL}) grows without bounds not only on the zeros, but also on the %
 segment  $-0.5<\epsilon<\epsilon^*_k \ , \ t=t_k^*$, where a step discontinuity arises in total phase of (\ref {EulerPrPhaseA_grh}). So, when computing  (\ref{DerLAngleSuTL}),   we find a kind of ridge  stemming from the $k^{th}$ zero $s^*_k=1/2+\epsilon^*_k+i t_k^*  \  \ , \   0.5 >\epsilon^*_k \ge 0$ toward negative  reals. 
\noindent %
Notice that for $p$  big (  i.e.  $ \ \ p^{1/2+\epsilon} >> 1 \ge \left |
\cos( \ln(p) t - \angle[\chi(p)]  )\right | \ \ \ $ ), the content of  braces summed  in  (\ref {DerLAngleSuTL} tends to be:

%

\be \label {FdipparamChiH}
 \left\{ \arctan \left( \frac{
 \sin(\ln(p) t    - \angle[\chi(p)]  )}{ p^{1/2+\epsilon}-  
\cos( \ln(p) t - \angle[\chi(p)]  )  }\right) \right \}_{t-\frac{\pi}{\ln(p^*)}   }^{t+\frac{\pi}{\ln(p^*)}  } \rightarrow
\ee
$$\left\{  \left( \frac{
 \sin(\ln(p) t    - \angle[\chi(p)]  )}{ p^{1/2+\epsilon}-  
\cos( \ln(p) t - \angle[\chi(p)]  )  }\right) \right \}_{t-\frac{\pi}{\ln(p^*)}   }^{t+\frac{\pi}{\ln(p^*)}  } \rightarrow
\left\{  \left( \frac{
 \sin(\ln(p) t    - \angle[\chi(p)]  )}{ p^{1/2+\epsilon}   }\right) \right \}_{t-\frac{\pi}{\ln(p^*)}   }^{t+\frac{\pi}{\ln(p^*)}  } =$$

$$=  \frac{ \sin \left(  \ln(p) \left\{ t +\frac{\pi}{\ln(p^*)} \right\}-  \angle[\chi(p)]       \right) -
  \sin \left( \ln(p) \left\{ t - \frac{\pi}{\ln(p^*)} \right\}-  \angle[\chi(p)]   \right)}
 {p^{1/2+\epsilon}}  =
 $$

%

\be \label {FdaLgaritmicInt}
 \frac{2 \ \cos( \ln(p) t -  \angle[\chi(p)]  )   \sin\left( \frac{\pi \ln(p)}{\ln(p^*)} \right)  }  
  {p^{1/2+\epsilon}}  
  \ee

 \noindent %
  When there is a zero in $s=s^*=1/2+\epsilon^* +i t^*$, then (\ref{DerLAngleSuTL}) grows without bounds for $s\rightarrow (s^*)^+$ , with constant $t=t^*$. It must grow without bounds, for  $[ \{p_{max} \rightarrow \infty\}  p^*\rightarrow \infty]$, also for $-0.5< \epsilon <\epsilon^*$ with $t= t^*$.
 
  $$PROOF.$$
  We have that  (\ref{DerLAngleSuTL})  grows without bound for $[ \{p_{max} \rightarrow \infty\}  p^*\rightarrow \infty]$,  for  $s\rightarrow (s^*)^+$ also if we start the sum  (\ref{DerLAngleSuTL}) from $P_1>>2$ .

    \noindent       Moreover    for  $-0.5<\epsilon <  \epsilon^* \ , \ t= t^*$,   as $p^{1/2+\epsilon}-  \cos[ \ln(p) t - \angle[\chi(p)]  >0 \ , \ \forall p \ge 2$, $\forall t$ and $ \forall \epsilon >-1/2$,  (\ref{FdipparamChiH}) written as (\ref{FdaLgaritmicInt}) shows that   $ \forall p> P_1$, where $P_1$ is a prime big enough,   the content in braces  in (\ref{DerLAngleSuTL})  is bigger, because multiplied by $p^{\epsilon^*-\epsilon}$, and preserving sign (as $\epsilon>-1/2$).

          \noindent  We  can use  the below reported Abel summation formula  
          \cite [p.~407] {Apostol:1967}  for $p \ge P_1$. 
          \noindent Let $\{ a_k \}$ and $\{ b_k \}$ be two sequences of complex numbers, 
          \noindent and let $A_n=\sum_{k=1}^n a_k$, then $\sum_{k=1}^n a_k  b_k = A_n b_{n+1}+ \sum_{k=1}^n  A_k(b_k - b_{k+1})$.
          \noindent We want to start from a big prime $ P_1$ so we pose: 
          
          \noindent $A_n=\sum_{k=1}^n a_k= \sum_{p_k=P_1} ^ {p_n}\frac{2 \ \cos( \ln(p_k) t -  \angle[\chi(p)]  )   \sin\left( \frac{\pi \ln(p_k)}{\ln(p^*)} \right)  }  
  {p_k^{1/2+\epsilon^*}} \rightarrow  \infty$ for $p_n\rightarrow \infty \ , \  \forall P_1$.
   \noindent Besides $b_k=p_k^{\epsilon^*-\epsilon}$, with $p_1=P_1$.
   \noindent As  $ b_{n+1}= - [(0-b_1)+ (b_1-b_2)+ . .+(b_n-b_{n+1} ) ]$, we can write:
   \noindent $\sum_{k=1}^n a_k  b_k = A_n b_{n+1}+ \sum_{k=1}^n  A_k(b_k - b_{k+1})=
   $
   $ A_nb_1+\sum_{k=1}^n [ A_n (b_{k+1}-b_k)- A_k
    (b_{k+1}-b_k) ] =$
    $  A_n b_1+\sum_{k=1}^n [ ( A_n- A_k )   (b_{k+1}-b_k)  ]  $
    $=  \left [ \sum_{k=1}^n a_k \right ]b_1+
    \sum_{k=1}^n 
    [ \{ \sum_{k_1=k+1}^n  a_{k _1} \}  (b_{k+1}-b_k)  ]  $.
    \noindent  For hypothesis on $\{ a_k \}$ and $\{0< b_k <b_{k+1} \}$ we have : 
   \noindent $ \left [ \sum_{k=1}^n a_k \right ]b_1\rightarrow \infty$  for $n\rightarrow \infty$.
     \noindent  Also $ [ \{ \sum_{k_1=k+1}^n  a_{k _1} \}  (b_{k+1}-b_k)  ]  \rightarrow \infty$  for $n\rightarrow \infty \ ,  \ \forall k %
     $.
     \noindent So  (\ref{DerLAngleSuTL}) must grow without bounds also for $-0.5< \epsilon <\epsilon^*$ with $t= t^*$, and, $p_n\rightarrow \infty$.
            $$  END \  of \  PROOF$$

\noindent So we find, again, the ridge for (\ref{DerLAngleSuTL}) from  $s_k^*$, non-trivial zeros of $\L(s.\chi_{primitive})$,  toward  negative reals, but,  by a different approach.

\noindent The PNT in arithmetic progression, with $q$ as congruence modulus \cite[p.~148] {Apostol:1976}, states that  for each $h-$class of primes, ($ p \equiv h \pmod{q} $,  $ gcd(h,q)=1$, and $0<h<q$) :

\be \label {PNTForAritmSucc} \pi(x)_{h,q}   \sim  \frac {Li(x)}{\phi(q)} \quad ; \quad  \forall h : gcd(h,q)=1 \  ,and \  \ Li(x)=\int_2^x \frac{dt}{ \ln(t)}
\ee

\noindent  Where $\pi(x)_{h,q}$ is the count of $h-$class primes till $x$, and, aritmetic function $\phi(q)$  \cite [p.~25]{Apostol:1976}, is the Euler totient function. There are   $\phi(q)$ classes of primes $ p \equiv h \pmod{q} $ where $ gcd(h,q)=1$, and, $0<h<q$. %
For all primes in a particular  $h-$class,  $\angle[ \chi(p)_{p \equiv h \pmod{q}  \  , \  gcd(h,q)=1 \ , \ 0<h<q}] = \angle[\chi(h)]$ is constant.
\noindent %
   For non principal characters (in particular for primitives ones)   ( see \cite [p.~256]  {Apostol:1976}) we have:
 \be \label {XiNonPrincip}
  \ \sum_{h<q : gcd(h,q)=1}e^{i \ \angle[\chi(h)]}=0 \ \ for  \ \ \chi \ \ not  \ principal
 \ee

   \noindent    We want to point out at $h-$prime classes in (\ref {DerLAngleSuTL}). If we limit to a finite $p_{max}$ we can write  (\ref {DerLAngleSuTL}                                                                                                                                                                                                                                                                                                                                                                                                                                                                                                                                                                                                            ) as the sum of  the $\phi(q)$ sums involving the $h-$prime classes:

 \be \label {DerLAngleSuTL2}                                                                                                                                                                                                                                                                                                                                                                                                                                                                                                                                                                                                             
   - \frac{\ln(p^*) }{ 2\pi}
   \sum_{h<q \ : \ gcd(h,q)=1  } 
  \left\{
  \sum_{   p<p_{max}   \ : \  p \equiv h (mod \  q)
  } 
  \left[
  \arctan \left( \frac{
 \sin(\ln(p) t    - \angle[\chi(p)]  )}{ p^{1/2+\epsilon}-  
\cos( \ln(p) t - \angle[\chi(p)]  )  }\right)\right]_{t-\frac{\pi}{\ln(p^*)}   }^{t+\frac{\pi}{\ln(p^*)}  }  
  \right\}
\ee

\noindent In order to take $p_{max } =(p^*)^m \rightarrow \infty$ we have to fill  all the sums   in natural order so that partial sums are not affected.
{ \bf Grouping with braces must not alter the  order of the sum because there is not absolute convergence } for $0<\Re(s)<1$. In other words, it is true that rearranging terms  we can change the limit of a conditional converging series at will, \cite [p.~411]  {Apostol:1967}, but, if we rearrange within a finite number of terms (i.e. till $p_{max}$)  generated by natural order, then, the usual rules of addition must hold, like in (\ref {EulerTrickFiniteTerms}).
\noindent Notice that  in (\ref{DerLAngleSuTL2}), in the sum in braces $\angle[\chi(p)]=\angle[\chi(h)]$ because it deals only with one $h-$prime class.

       \subsection  {LEMMA 4: use of logarithmic integral }

    The sum of the $\phi(q)$  integrals  $h-$dependent  like (\ref    {HDepInt}) is null $\forall p_{max}$:

     \be \label {HDepInt} -
 \frac{\ln(p^*)} { \pi} 
  \sum_{h<q : gcd(h,q)=1 }
   \frac {1}{\phi(q)}  \int_{2}^{p_{max}}   \frac{\cos(\ln(y) \ t - \angle[\chi(h)])} {y^{1/2+ \epsilon}} \sin\left(  \frac{\pi \ln(y)}{\ln(p^*)} \right) d [ Li(y) ] = 0  \quad    \forall p_{max } %
   \ee 
\noindent %
\noindent PROOF.
\noindent  Each $h: 0<h < q $, defines an $h-$class   of primes  $mod \ q$

 $$ \sum_{h<q : gcd(h,q)=1 }
   \frac {1}{\phi(q)}  \int_{2}^{p_{max}}   \frac{\cos(\ln(y) \ t - \angle[\chi(h)])} {y^{1/2+ \epsilon}} \sin\left(  \frac{\pi \ln(y)}{\ln(p^*)} \right) d [ Li(y) ] =
 $$
 $$ \sum_{h<q : gcd(h,q)=1 }
   \cos( \angle[\chi(h)])
    \frac {1}{\phi(q)}  \int_{2}^{p_{max}}   \frac{\cos(\ln(y) \ t   )} {y^{1/2+ \epsilon}} \sin\left(  \frac{\pi \ln(y)}{\ln(p^*)} \right) d [ Li(y) ]  +
 $$
  $$ \sum_{h<q : gcd(h,q)=1 }
  \sin( \angle[\chi(h)])
    \frac {1}{\phi(q)}  \int_{2}^{p_{max}}   \frac{\sin(\ln(y) \ t   )} {y^{1/2+ \epsilon}} \sin\left(  \frac{\pi \ln(y)}{\ln(p^*)} \right) d [ Li(y) ] =
  $$
  Integrals are not $h-$dependent and so can be put into common factor obtaining:
  
  $$\frac{  \int_{2}^{p_{max}}   \frac{\cos(\ln(y) \ t   )} {y^{1/2+ \epsilon}} \sin\left(  \frac{\pi \ln(y)}{\ln(p^*)} \right) d [ Li(y) ] }{\phi(q)}
   \sum_{h<q : gcd(h,q)=1 }
   \cos( \angle[\chi(h)])
     +
 $$

  $$\frac{  \int_{2}^{p_{max}}   \frac{\sin(\ln(y) \ t   )} {y^{1/2+ \epsilon}} \sin\left(  \frac{\pi \ln(y)}{\ln(p^*)} \right) d [ Li(y) ] }{\phi(q)}
   \sum_{h<q : gcd(h,q)=1 }
   \sin( \angle[\chi(h)])
 $$
 
 \noindent But  for (\ref {XiNonPrincip})  $$ \sum_{h<q : gcd(h,q)=1 }
   \sin( \angle[\chi(h)]) =0 \ \ \ \ \ \  \ \ , and, \ \ \ \ \  \sum_{h<q : gcd(h,q)=1 }
   \cos( \angle[\chi(h)]) =0$$
   
  \noindent So $\forall p_{max} $ (\ref{HDepInt}) is null  independently of the values assumed by integrals at common factor.
       \noindent Notice , { \bf always complying with the natural  order of the sum in (\ref{DerLAngleSuTL2}) till $p_{max}$}, we could distribute (only for indexes less than running index) the $\phi(q)$ integrals  in (\ref{HDepInt}) into the $\phi(q)$ braces of (\ref{DerLAngleSuTL2}). If $\theta(z,h,t)=\ln(z)t - \angle[\chi(h)]$, we get  by  (\ref{Estremi}):

      \be \label {WithLogInt}
       - \frac{\ln(p^*) }{ 2\pi}
   \sum_{h<q \ : \ gcd(h,q)=1  } 
      \ee   
      $$ \left\{
  \sum_{   p<p_{max}   \ : \  p \equiv h (mod \  q)
  } 
  \left[ 
  \arctan \left( \frac{
 \sin( \theta(p,h,t) )}{ p^{1/2+\epsilon}-  
\cos( \theta(p,h,t) )  } \right)\right]_{t_1} ^{t_2} -
\frac{
\int_{2}^{p_{max}}   \frac{2 \cos(\theta(x,h,t))   } {y^{1/2+ \epsilon}} \sin\left(  \frac{\pi \ln(x)}{\ln(p^*)} \right) d [ Li(x)}{\phi(q)}
  \right\}
  $$

                     \noindent  as usual $[ \{p_{max} \rightarrow \infty\}  p^*\rightarrow \infty]$,  and, the  overall result does not change, provided {\bf all } operations for $p<p_{max}$ and $x<p_{max}$ are carried out at same time, while { \bf no}  computations with greater  $p_{max}$ is  just done.
       \noindent     This is true because   (\ref{DerLAngleSuTL}) converges $\forall s  \setminus \{  singular \ points\} $, and, for Lemma 4 the contribution of  (\ref{HDepInt}) till whatever  $ p_{max} $ is always  zero. Then is applicable  
        \cite [p.~385 Theorem 10.2] {Apostol:1967}:  if two series of complex terms converges, then converges also their sum to the sum of their limits.

  \subsection  {LEMMA 5: Euler Product on critical line  }  \label {Lemma5}

      \noindent We want to compute $ \left[\frac{\partial  \angle[\xi(s,\chi)] }{\partial t} \right]
$ around $\epsilon=0$  for $t \ne t^*_k$.

 \noindent Neglecting  minor contributions, from appendix \ref {GammaDer}, in particular (\ref{PhaseVarXiAndL}),  (\ref{Asympt}), and  (\ref {3Cases}), we have:
  
  
  $$ \left[\frac{\partial  \angle[\xi(s,\chi)] }{\partial t} \right]_{\epsilon \approx 0 \ , \ t \ne t^*_k}=
$$


\be \label {PhaseVarXiAndL_}
   \left[\frac {\partial \angle[L(s,\chi_{primitive})]}{\partial t}\right]_{\epsilon=0}+ \epsilon  \times \left[\frac {\partial^2 \angle[L(t,\epsilon=0,\chi_{primitive})]}{ \partial \epsilon \partial t}\right]_{\epsilon=0} +\ln\left[ \sqrt{\frac{t q }{2 \pi}} \right]
    \quad    \pm \frac {\epsilon}{4t^2} \ \ ; \   \ \forall t \ne t^*_k  
    \ee

  \noindent%
\noindent 
Where $\pm$ refers rispectively to odd and even characters  (\ref {3Cases}).
 We want to use (\ref {DerLAngleSuTL}) to compute previous equation for:
  $\epsilon=0$. In order to avoid  singular points  we choose for  $\epsilon \le \epsilon^*_k$
  $ \quad  
   |t -t^*_k| > \frac{2 \pi}{\ln(p^*)} \ 
 $,
   while this condition can be dropped for  $ \  \epsilon> \epsilon^*_k $.
 So we can write ($t>T_{Asymp}(\alpha)$ (\ref{Asympt}) ):

  \be \label {Uguale0Tris_} 
  \ln\left[ \sqrt{\frac{t  \times q}{2 \pi}} \right]+
   \left\{ \left[  \frac {\Delta \angle [ L_{EP}(s,\chi_{primitive})] }{\Delta t}  \right]_{ p_{max}  \rightarrow \infty}  \right\}_{  p^*  \rightarrow \infty}
  \rightarrow  \left[ \frac {\partial \angle[\xi(t,\epsilon,\chi )]  } {\partial t}  \right]_{\epsilon=0} = 0
   \ee

  \noindent  PROOF. 
  From Lemma 2, appendix \ref {GammaDer} and  (\ref  {DerLAngleSuTL}                                                                                                                                                                                                                                                                                                                                                                                                                                                                                                                                                                                                              ) follows (\ref{Uguale0Tris_} ). End of proof.
  
 \noindent In half-plane $\Re[s]>0$ (see also figg. \ref{Q3-Q4-m3}, \ref   {Q5-3Camp}, and , \ref  {Simmetries})) , we can remark that: 
 \begin{itemize}
\item (1) $\mathcal{L} [ \xi (s,\chi) ]$ is defined $\forall s$, see (\ref{AngMomDef}), 
\item (2) $\frac {\partial \angle[\xi(t,\epsilon,\chi )]  } {\partial t}$  is defined everywhere but points $s^*_k=1/2+\epsilon^*_k + i t^*_k$  where $L(s^*_k,\chi)=0$. See (\ref{VarPhase}), where, only phase difference between  vector of real and imaginary components with respect to vector of  their $t-$derivative matters.
\item (3)  $\frac {\Delta \angle [ L_{EP}(t,\epsilon,\chi) ] }{\Delta t}   $ instead is defined everywhere but points $ ( \epsilon \le \epsilon^*_k \ , t=t^*_k )$ , see (\ref  {DerLAngleSuTL}                                                                                                                                                                                                                                                                                                                                                                                                                                                                                                                                                                                                              ). Here is used total phase.
\end{itemize}



     \subsection {   LEMMA 6:  Sign estimate  of $\frac{\partial ^2\angle[L(s,\chi_{primitive})]}{\partial \epsilon \partial t}$   }  \label {Lemma6}

Suppose $ \epsilon >0$ but arbitrarily small, $ |t -t^*_k| > \frac{2 \pi}{\ln(p^*)}$, %
 with (\ref {Estremi}), and  
 $p^*$ big (fixed), and $p_{max} = (p^*)^m 
 < \infty$.
   \noindent Then,for (\ref {Uguale0Tris_}) we have that (\ref {DerLAngleSuTL})
 evaluated in $\epsilon=0$ is tending to $ - \ln\left[ \sqrt{\frac{t  q}{2 \pi}} \right] <0$ as $[ (p_{max} \rightarrow \infty) p^* \rightarrow \infty]$', 
  in order to verify Lemma 2. We affirm that the same (\ref {DerLAngleSuTL}),  
 evaluated for small $\epsilon >0$  is always negative but with lower absolute value, so that (\ref {Uguale0Tris_}  ) has a positive increment for  small $ \epsilon >0$. i,e
 \footnote { 
 at least for odd primitive characters, see (\ref {PhaseVarXiAndL_}),
  because we have $
    \epsilon  \times \left[\frac {\partial^2 \angle[L(t,\epsilon=0,\chi_{primitive})]}{ \partial \epsilon \partial t}\right]_{\epsilon=0} +\
     \frac {\epsilon}{4t^2} >0 $. We are aiming to the result of following Lemma 7:
     $$ \left[\frac{\partial^2  \angle[\xi(s,\chi)] }{\partial \epsilon \partial t} \right]_{\epsilon \approx 0 \ , \ t \ne t^*_k} >0
$$  
     }
     :
  
   \be \label {LiDeltaXEps_2_GRH}
\left|
 \frac {\Delta \angle [ L_{EP}(t,\epsilon=0,\chi,p_{max}) ] }{\Delta t}    \right| >%
\left|
 \frac {\Delta \angle [ L_{EP}(t, \epsilon>0,\chi ,p_{max}) ] }{\Delta t}    \right|%
 \ee
  
  \noindent with $ \left\{ \left[  \frac {\Delta \angle [ L_{EP}(t,\epsilon=0,\chi_{primitive}, p_{max})] }{\Delta t}  \right]_{ p_{max}  \rightarrow \infty}  \right\}_{  p^*  \rightarrow \infty} \rightarrow -  \ln\left[ \sqrt{\frac{t  \times q}{2 \pi}} \right]$

$$ PROOF$$

   %
   \noindent  Our aim is difference between (\ref {DerLAngleSuTL}) and (\ref {HDepInt}), i.e. (\ref {WithLogInt}), computed in natural order. For  Lemma 4 the result for  $[ \{p_{max} \rightarrow \infty\}  p^*\rightarrow \infty]$  does not change from  (\ref {DerLAngleSuTL}) that converges 
   apart singular points. See (\ref{EulerPrPhaseA_grh}).

   \noindent      Let us analyze the peculiar oscillations in  (\ref{WithLogInt}) integrands  and, to simplify, suppose, at first, that $x< p^*$, i.e $ \sin\left(  \frac{\pi \ln(x)}{\ln(p^*)} \right) > 0$. Considering that
 the zero crossing of $\cos[  \ln(x)t   - \angle[\chi(h)]  ] $ are at $   \ln(x)t   - \angle[\chi(h)] =2 k \pi  \mp \frac{\pi}{2}\rightarrow  \ln(x)=\frac{2 k \pi \mp\frac{\pi}{2} +\angle[\chi(h)] }{t}$, let us see zero transitions (subscript $0t$ means `zero transition')
 \footnote{Notice the $\pm\pi/2$ at exponent}.

    \begin{itemize}

\item Zero transition at increasing $\cos[\ln(x)t   - \angle[\chi(h)]]$ :  

   \be \label {TransizMenoPiMezzi__} 
     \quad x_{0t}(k,h)= e^{ (2 \pi k- \pi/2    + \angle[\chi(h)]) / t}   \quad increasing \quad cosine
   \ee
\item Zero transition  at decreasing   $\cos[\ln(x)t   - \angle[\chi(h)]]$ : %
    
    \be \label {TransizPiuPiMezzi__} : 
     \quad x'_{0t}(k,h) = e^{ (2 \pi k+\pi/2   + \angle[\chi(h)]) / t}   \quad decreasing \quad cosine
     \ee
        \end{itemize}

        \noindent      Corresponding points to $ x_{0t}(k+1,h), x'_{0t}(k,h)$   defined on (\ref {HDepInt})   can be defined on   (\ref{DerLAngleSuTL}) : $ y_{0t}(k+1,h), y'_{0t}(k,h)$,  thinking 
         (\ref {FdipparamChiH}) as continuous function of $y=p$.
  \noindent In more detail: let us choose $\phi(q)$ functions of $y$, continuous variable, in place of discrete $p$:  
    \be \label {ThoughtContinuous} 
    F(y,\chi(h),\epsilon)= \left\{ \arctan \left( \frac{
 \sin(\ln(y) t    - \angle[\chi(h)]  )}{ p^{1/2+\epsilon}-  
\cos( \ln(y) t - \angle[\chi(h)]  )  }\right) \right \}_{t-\frac{\pi}{\ln(p^*)}   }^{t+\frac{\pi}{\ln(p^*)}  }
\ee   
      
 \noindent   $ y_{0t}(k+1,h), y'_{0t}(k,h)$ are defined as positive slope zero transitions   and negative slope zero transition respectively  (similarly to    $ x_{0t}(k+1,h), x'_{0t}(k,h)$) 
     
          \noindent  Besides  for $k$ big (i.e. $p$ or $y$  big) $ x_{0t}(k+1,h) \rightarrow  y_{0t}(k+1,h)$ and $ x'_{0t}(k,h) \rightarrow  y'_{0t}(k,h)$
   \noindent  as we want. To  build    (\ref {DerLAngleSuTL}) minus (\ref {HDepInt}), i.e (\ref {WithLogInt}),  then, negative contributions for (\ref {HDepInt}) become positive contributions for final result.

\noindent  Negative  contribution for  (\ref {HDepInt}) 
( i.e. positive contribution for final result)
 is given by   the  (\ref {HDepInt}) integral in the intervals  $\Delta x_{k,h}$ :

 \be \label {HalfRotation_GRH}
 x'_{0t}(k,h) -  x_{0t}(k,h)=    x_{0t}(k,h)  [e^{ ( \pi  / t)}-1 ]= \left( \frac {\pi}{t}+  \frac{\pi^2}{2t^2} +
  . .\right) e^{ (2 \pi k -\pi/2 + \angle[\chi(h)] ) / t} =\Delta x_{k,h} \approx \frac {\pi}{t} \hat{x}_{0t}(k,h)    
 \ee
 
 \noindent   where $\hat {x}_{0t}(k,h)$ is a suitable value inside $\Delta x_{k,h} $

 \be \label {SingleOillLi1}  
   O^+_{k,h}(\epsilon)_{\int dLi}  = \left|  %
  \frac {\ln(p^*)} {\pi \phi(q)}  
           \int_{ x >x_{0t} (k,h)  \ : \  p \equiv h (mod \  q)}^{ x <x'_{0t} (k,h)} 
            \frac{ 
              \cos [\ln(x)t- \angle[\chi(h)] ]\sin [\pi \ln(x) /\ln(p^*)] }{\sqrt{x} x^\epsilon}   dLi(x)
              \right|
 \ee

 \noindent Negative contribution of  (\ref {DerLAngleSuTL})  (and of final result) is given by the 
 $h-$class primes (weighted by (\ref  {ThoughtContinuous} ) function)  inside  the analogous intervals %
 $\Delta y_{k,h}$ ( that for $k$ big tends to the previous one $\Delta x_{k,h}$):
 
 \be \label {SingleOillSum1} 
  O^-_{k,h} (\epsilon)_{\sum}  = \left|
  \frac {\ln(p^*)} {2 \pi}  
           \sum_{ p >y_{0t} (k,h) \ : \  p \equiv h (mod \  q)}^{ p <y'_{0t} (k,h)} 
            F(p,\chi(p),\epsilon)
              \right|
 \ee


                                    \noindent  Positive  contribution for (\ref {HDepInt})  (i.e. negative contribution for final result)  is given by   the  (\ref {HDepInt}) integral in the intervals  $\Delta x_{k,h}'$ (see (\ref {HalfRotation_GRH}) ) : %

      \be \label {HalfRotation2_GRH}
   x_{0t}(k+1,h)-x'_{0t}(k,h)=   x'_{0t}(k,h) [e^{ ( \pi  / t)}-1]
     =
     \left( \frac {\pi}{t}+  \frac{\pi^2}{2t^2} +
  . .\right)    e^{ (2 \pi k +\pi/2 + \angle[\chi(h)] ) / t}
      =\Delta x_{k,h}' \approx \frac {\pi}{t} \hat{x}'_{0t}(k,h)  
   \ee

      \be \label {SingleOillLi2}
    O^-_{k,h}(\epsilon)_{\int dLi} = \left|   %
  \frac {\ln(p^*)} { \pi \phi(q)}  
           \int_{ x >x'_{0t} (k,h) \ : \  p \equiv h (mod \  q)}^{ x <x_{0t} (k+1,h)} 
            \frac{ 
              \cos [\ln(x)t- \angle[\chi(p)] ]\sin [\pi \ln(x) /\ln(p^*)] }{\sqrt{x} x^\epsilon}  dLi(x)
              \right|
 \ee

 \noindent Positive contribution of  (\ref {DerLAngleSuTL})  (and of final result) is given by the 
 $h-$class primes (weighted by (\ref  {ThoughtContinuous} ) function)  inside  the analogous intervals %
 $\Delta y_{k,h}'$ ( that for $k$ big tends to the previous one $\Delta x_{k,h}'$)

   \be \label {SingleOillSum2}
    O^+_{k,h}(\epsilon)_{\sum} = \left|
  \frac {\ln(p^*)} {2 \pi}  
           \sum_{ p >y'_{0t} (k,h) \ : \  p \equiv h (mod \  q)}^{ p <y_{0t} (k+1,h)} 
           F(p,\chi(p),\epsilon)
              \right| 
 \ee       
                

                  \noindent Until now  we have dealt with $p_{max} < p^*$. If we let $p^*$ fixed and $p_{max}$ increases without bound  we have that in the intervals :
 
 \be \label {Campate}
 (p^*)^m <p_{max}<(p^*)^{m+1}  \ee

 $$-$$

 \noindent 
if $m$ is even, the factor $ \sin\left(  \frac{\pi \ln(x)}{\ln(p^*)} \right)$ in (\ref {HDepInt}) is positive, so,    the contributions in corresponding  $x-$intervals are as above. 
In other words, the sign contributions of  final result, i.e. (\ref{DerLAngleSuTL}                                                                                                                                                                                                                                                                                                                                                                                                                                                                                                                                                                                                       ) minus  (\ref {HDepInt}), or , (\ref {WithLogInt}),
 \noindent  in intervals $\Delta x ,\Delta x' ,\Delta y ,\Delta y'    $,    is as that   stated   in ( \ref {SingleOillSum1}),  (\ref {SingleOillSum2}) , and,  (\ref {SingleOillLi1}),     (\ref{SingleOillLi2}). 
 \noindent  Instead where $m$, in (\ref {Campate}),  is odd, then,  in every  interval above  
 the sign changes concurrently \footnote{  intermediate situations are avoided by the choice $p_{max}=(p^*)^m$, but, in any case 
  they  
 can be neglected.}. 
      %
                  \noindent So we could compute    (\ref {DerLAngleSuTL})  minus       (\ref {HDepInt}), till $p_{max}$, ( interchanging $O^+  \leftrightarrow O^-$ when $m$       is odd in (\ref {Campate}) ) as in   :

                  \be \label {FinRes}  
                   \sum_k  \sum_h \left[ O^+_{k,h }(\epsilon)_{\int dLi}  + O^+_{k,h }(\epsilon)_{\sum} 
                   \right]
                   -
     \sum_k  \sum_h \left[ O^-_{k,h } (\epsilon)_{\sum}  + O^-_{k,h }(\epsilon)_{\int dLi} \right]
                  \ee

     \noindent  
     {\bf We stress again, (\ref{FinRes} ) is computed without changing the order  given by $p$ or $x=y$, like in (\ref{DerLAngleSuTL2}                                                                                                                                                                                                                                                                                                                                                                                                                                                                                                                                                                                                           ) , i.e. partial sums are unaffected}. 
     \noindent Limit difference  in  (\ref{FinRes} ) is in   form:    $\infty- \infty$, and, we know that if $|t-t_k^*|>\frac {2 \pi}{\ln(p^*)}$ limit result is finite 
     . 
    %
     %
     %
     %
     %
     %
    \noindent  In other words: positive quantities:  $ O^+_{k h} ( \epsilon)_{\int dLi} 
 \ \ , \ \  O^-_{k,h }( \epsilon)_{\int dLi} $,and , $  O^+_{k,h} (\epsilon)_{\sum}  \ \ , \ \  O^-_{k ,h}(\epsilon)_{\sum} $  like   in ( \ref {SingleOillSum1}), and  (\ref {SingleOillSum2}) , and,  (\ref {SingleOillLi1}),     (\ref{SingleOillLi2}) are  only { \bf ``place-holders''} filled according to progress of $x=y<p_{max}$, and with chosen  sign following the $m$ parity in (\ref {Campate}).
                  %
               %
                 %
          %

      \noindent  Consider now the ratios :

            \be \label {Ratios} 
            \frac{ O^-_{k,h }(\epsilon')_{\sum} }{ O^+_{k,h }(\epsilon')_{\int dLi}} \ \ \ \  \ \ ; \ \ \ \ \ \  \frac{ O^+_{k,h }(\epsilon')_{\sum} }{ O^-_{k,h }(\epsilon')_{\int dLi}} 
            \ee
              for each $h$ and for $k$ big,       defined respectively 
              in  intervals $\Delta y_{k,h} \rightarrow  \Delta x_{k,h} $ ( see (\ref {SingleOillLi1}) and (\ref   {SingleOillSum1} )   ) , 
               and 
                in  intervals           $\Delta y_{k,h}' \rightarrow  \Delta x_{k,h}'$ ( see (\ref {SingleOillLi2})  and (\ref       {SingleOillSum2}) )
         . $$LEMMA \ 7$$
         
         \noindent {\it LEMMA 7}:  $\forall h,k$  as $k \rightarrow \infty$ both (\ref {Ratios} ) tend to unity.
       
       %
       $$ PROOF$$
      \noindent  PNT ratio   in small intervals \cite  [p.~23]  {Heath-Brown : 
1988}  assures that  
%
   \footnote {If Riemann hypothesis is true , then  (\ref  {PNTInSmallInt}) will be true with $\theta=1/2$ instead of the Huxley unconditional reached result $\theta=7/12=0.58333. . .$(\cite {Heath-Brown : 1988}).Recenly a breackthrough has improved the unconditional result to  $\theta=17/30=0.5666....$ (\cite  {Guth : 2026}) value. But if $\theta=1/2 $ this does not imply RH see Terence Tao Blog (\cite {TerenceTao:2026})}:

    \be \label {PNTInSmallInt}
    \pi(x+\Phi(x) )-\pi(x)  \sim \frac{\Phi(x)}{\ln(x)}   \ \ with \ \ x \rightarrow \infty \ \ , if \ \Phi(x)\ge x^{7/12-\epsilon(x)}
    \ee
    
    \noindent where $\epsilon(x)$ in (\ref {PNTInSmallInt}) tend to zero as $x\rightarrow \infty$.  Imposing this on equi-sign  intervals, for
   (\ref {HalfRotation2_GRH}) and  (\ref {HalfRotation_GRH})
    \footnote {In \cite [p.~78]  {Edwards:1974cz} :  $ |\pi(x)-Li(x)| < Li(x)e^{-\sqrt{c \ln(x)}}$ , for $c$ positive constant and sufficiently large x.
      In  \cite  {Trudgian: 2016} : 
      $ |\pi(x)-Li(x)| \le  \frac{0.2795 x} {(\ln(x))^{3/4} } e^{-\sqrt {\frac{\ln(x)} {6.455} } } \ $  from $\ \; x>229 $ 
      }:

 \be \label {ToFulfillHB} \Delta x_{k,h} =\frac{\pi}{t}x >x^{7/12} 
 \rightarrow \frac{\pi}{t}x^{5/12} >1 \rightarrow 
   x> \left( \frac {t}  {\pi}  \right)^{12/5} 
 \ee

\noindent i.e.  for $x$ big enough the distribution of primes in equal-sign intevals $\Delta x_{k,h} \rightarrow\Delta y_{k,h}$, and, for (\ref  {PNTForAritmSucc}) also $h-$class primes, can be close to theorethical one as we want.
          \noindent   In 
\cite  [p.~65] {RosserSchoenfeld1962} it is dealt about the sum of prime sampled functions: { \it
``in the neighborhood of the number $x$ the average density of the primes
is $1/\ln( x)$. On this basis, if one should wish an estimate for the sum of $F(p)$
over all primes $p < x$, the natural approximation would be''}

  
  \be \label {shoenfeld}  \sum_{p<x} F(p)  \sim
   \int_2^x \frac{F(t) dt}{\ln(t)}  = \int_2^x F(t) d[Li(t)]  \rightarrow  \frac{\sum_{p<x} F(p)}{\int_2^x F(t) d[Li(t)] }  \rightarrow 1  \ \ for \ \ x \rightarrow \infty
  \ee

  \noindent   PNT  theorem  is same statement of  (\ref{shoenfeld}) with  $F(p) =1$.

\noindent We use then (\ref{shoenfeld}) in intervals  (\ref {HalfRotation2_GRH}) and  (\ref {HalfRotation_GRH})  for $k\rightarrow \infty$, with  $F(p)$ function as in (\ref{FdaLgaritmicInt}) or as  in (\ref{ThoughtContinuous}).   For $p,x,y ,k \rightarrow \infty$
(\ref{FdaLgaritmicInt})  and (\ref{ThoughtContinuous}) become the same function in intervals $\Delta x  \rightarrow \Delta y \ \  ,  \ \ \Delta x'  \rightarrow    \Delta y'    $ that grow without bouds, fulfillig so (\ref {ToFulfillHB}). So  thesis is proved. 
 \noindent Besides   if $\Delta x_{k,h}  \approx \Delta y_{k,h}$, and,  (\ref {ToFulfillHB}) is fulfilled,  the $h-$class primes are $\approx \frac{ \Delta x}{\phi(q) \ln(\hat{x})}$, so, far from $p_{max}\approx m  p^*$,  neglecting $|\cos()|$ and  $| \sin()|$ factors from (\ref {FdaLgaritmicInt}), or equivalent shaping from (\ref{ThoughtContinuous}), and, with $\hat{x}$ as a suitable value inside $\Delta x_{k,h}$, $O^\pm_{k,h}$  can be evaluated  as:
 \be \label {DiMassima}
  O^\pm_{k,h}   \approx \frac {\pi}{t}\frac{ 2 \hat{x}^{0.5   } }{\phi(q)\ln(\hat{x})}   \ee

\noindent So, %
numerators and denominators of  (\ref{Ratios} ) grow without bound as $x=y \rightarrow \infty$ as (\ref{DiMassima}), i.e. as $ \approx K \times \frac{ \sqrt{x}}{\ln(x)}$.
  
  $$ END \ OF \ PROOF$$

    \noindent Now let us build the following  ratios   between all positive,$\sum_{+}^{p_{max}}(\epsilon'')$, and all negative contributions,$ \sum_{-}^{p_{max}}(\epsilon')$, to final result (\ref{FinRes} ):

   \be \label {RappSommeGRH}
    \frac{ \sum_{+}^{p_{max}}    (\epsilon'')}{ \sum_{-}^{p_{max}}    (\epsilon')}
      =
    \frac { \sum_k \sum_h O^+_{k,h }(\epsilon'')_{\int dLi}  + O^+_{k,h }(\epsilon'')_{\sum}  }{\sum_k\sum_h  O^-_{k,h }(\epsilon')_{\sum} +O^-_{k,h }(\epsilon')_{\int dLi}   }=\rho(\epsilon''=\epsilon' =0,p_{max } \rightarrow \infty) \ \ \  \rightarrow  \ \ 1 \ \ \,   
   \ee

  \noindent We can expect a  faster convergence   then PNT ratio, because, alternating   between  $ \sum$ and $ \int$ at nominator and denominator, results in a   reshufling that  eases the convergence of the ratio to 1 \footnote{ {\bf This   is true in a continuous way  also considering 
  that $p_{max}=[p^*]^m $ and $m\rightarrow \infty$. Besides equal-sign intervals at higher $p$, or, $k$,  have greater weight then equal-sign intervals at lower $p \approx x,$ or $y$ }.  See (\ref {DiMassima}). So the rare zero transitions of  factor $ \sin\left(  \frac{\pi \ln(x)}{\ln(p^*)} \right)$ in $\Delta x_{h,k}$ are without influence in the ratio (\ref{RappSommeGRH}).
  \noindent  The differences between $\Delta x$ and $\Delta y$ fades away with big $x=y \approx p$, so  they do not influence the ratio  (\ref {RappSommeGRH}). Besides    when $p_{max}$  satisfy (\ref{Campate}) with odd $m$ then in ratio  (\ref {RappSommeGRH}) we add $O^-_{k,h }(\epsilon'')_{\sum} +O^-_{k,h }(\epsilon'')_{\int dLi} $ at  numerator and $O^+_{k,h }(\epsilon')_{\int dLi}  + O^+_{k,h }(\epsilon')_{\sum} $ at denominator, and, we have always the ratio beween all positive versus all negative contributions }
 

  
    \noindent  Let us call all negative contributions till $p_{max}$: 
     \be \label {AllNegContr}   \sum_{-}^{p_{max}}    (\epsilon)     =\sum_k  \sum_h O^-_{k,h } (\epsilon)_{\sum}  + O^-_{k,h }(\epsilon)_{\int dLi}
     \ee

   \noindent and all positive contributions till $p_{max}$ ($t\ne t*_k$ is implied) : 
     \be \label {AllPosContr}  \sum_{+}^{p_{max}}    (\epsilon)     =\sum_k  \sum_h O^+_{k,h } (\epsilon)_{\sum}  + O^+_{k,h }(\epsilon)_{\int dLi}
     \ee   
     \noindent so  we can write  (\ref{FinRes} ) as                  $=  \sum_{+}^{p_{max}}    (\epsilon)  -\sum_{-}^{p_{max}}    (\epsilon)  $, and, for  $\epsilon=0$, for Lemma 2,   
    it is tending to $ - \ln\left[ \sqrt{\frac{t  q}{2 \pi}} \right] <0$ as $[ (p_{max} \rightarrow \infty) p^* \rightarrow \infty]$

   \noindent In (\ref{RappSommeGRH}) we can change $\epsilon$. in particular we use $0$ and $\epsilon>\approx0$. With definition (\ref{AllPosContr}) and (\ref{AllNegContr}) we can write ($t\ne t_k^*$):
   
   $$
   \frac{\sum_+^{p_{max}} (\epsilon)}{\sum_-^{p_{max}} (\epsilon)} \rightarrow_{p_{max}\rightarrow \infty} 1 \ \ ; \ \   \frac{\sum_+^{p_{max}} (0)}{\sum_-^{p_{max}} (0)} \rightarrow_{p_{max}\rightarrow \infty} 1
   $$

  \noindent So we have also

     \be \label {RatioUnequalEps}
     \rho_+(\epsilon,p_{max})=
     \frac{\sum_+^{p_{max}} (\epsilon)}{\sum_+^{p_{max}} (0)} \rightarrow_{p_{max} \ big}   \ \ \rightarrow  \ \   \frac{\sum_-^{p_{max}} (\epsilon)}{\sum_-^{p_{max}} (0)}=
      \rho_-(\epsilon,p_{max})=   \rho(\epsilon,p_{max})
   \ee
   
   \noindent Besides  for $p_{max}$  big $ \rho_+(t,\epsilon,p_{max})\rightarrow \rho_-(t,\epsilon,p_{max})$ as we want.

\noindent We need to be cautius in 
(\ref {RatioUnequalEps}). If we fix $\epsilon$ and let $p_{max} \rightarrow \infty$ we get $ \rho(\epsilon,p_{max}) \rightarrow 0$. This  because,  neglecting modulation with absolute value cosine, characters phases, and , neglecting continuous versus prime-powers sampled interval, we deal with an object something like:

 \be \label {RatioPrototype}
 \frac{\int_0^ { p_{max}  }  \frac{dx}{x^{1/2+\epsilon}}}{  \int_0^  { p_{max}  }  \frac{dx}   {x^{1/2}} } =
\frac{
\frac{2  { p_{max}  } ^{1/2-\epsilon}}{1-2\epsilon}
}{\frac{2  { p_{max}  } ^{1/2}}{1}} =
\frac{ { p_{max}  } ^{-\epsilon}}   {1-2\epsilon}
\ee
 
\noindent  But for each $p_{max}$ big, but,  {\bf fixed}, if we choose $0<\epsilon< O( \ln^{-1}(p_{max}) )$ we can , using (\ref{RatioUnequalEps}),   { \bf approach 1 from below } as we want.

\noindent So generalising (\ref{FinRes}) we can write for each $p_{max}$ big  $ \exists \ \ \delta  , \gamma >\approx 0$:

\be \label {ApproachFromBelow} 
\left |  \sum_+^{p_{max}} \left(  0<\epsilon<\frac{\delta}{\ln(p_{max})}    \right)-\sum_-^{p_{max}}  \left(0<\epsilon<\frac{\delta}{\ln(p_{max})}\right) \right|=
 X \left( 0<\epsilon<\frac{\delta}{\ln(p_{max})} ,p_{max}\right) =
\ln\left[ \sqrt{\frac{t }{2 \pi}} \right]   \ \  - \ \gamma
\ee

\noindent Where (dipendence on $t\ne t^*_k$ is implied) .  In short we can write:

\be \label {Xdef}
\sum_+^{p_{max}} (\epsilon)-\sum_-^{p_{max}} (\epsilon) \rightarrow  -  X(\epsilon,p_{max}) > - \ln\left[ \sqrt{\frac{t }{2 \pi}} \right] 
\ee


 \noindent In the sense of (\ref {ApproachFromBelow}), we can refer to 
 the ratio of all positive  contributions to final result (\ref{FinRes} ) at different $\epsilon$ value. In detail :  $\forall h < q \ , gcd(h,q)=1$,   $\epsilon'' >0$, and $\epsilon'=0$:

      \be \label {RappSommePosGRH} 
      \frac{ \sum_{+}^{p_{max}}    (\epsilon)}{ \sum_{+}^{p_{max}}    (0)}
      =
        \frac { \sum_k  \sum_h O^+_{k,h }(\epsilon)_{\int dLi}  + O^+_{k,h }(\epsilon)_{\sum}  }{\sum_k  \sum_h O^+_{k,h }(0)_{\int dLi}  +O^+_{k,h }(0)_{\sum}  }=\rho^+(\epsilon,p_{max})  \  <  \\ 1  ;   \ \  \epsilon >\approx 0
   \ee

    \noindent and the same for negative contributions ( we do not mention the 
     interchange we do  when $m$  is odd in (\ref {Campate})). Besides in intermediate cases; i.e. when $p_{max}  \approx (p^*)^m$,  we could have a split  in intervals   $\Delta x ,\Delta y$ or $\Delta x'  , \Delta y'    $,  that, in any case, do not alter the limit of  (\ref{RappSommeGRH})  or (\ref {RappSommePosGRH}) ): %

    \be \label {RappSommeNegGRH}
    \frac{ \sum_{-}^{p_{max}}    (\epsilon) } { \sum_{-}^{p_{max}}    (0)}
      =
     \frac { \sum_k  \sum_h O^-_{k,h } (\epsilon)_{\sum}  + O^-_{k,h }(\epsilon)_{\int dLi}  }{\sum_k  \sum_h O^-_{k,h }( 0)_{\sum} +O^-_{k,h } (0)_{\int dLi}   }=
     \rho^-(\epsilon, p_{max}) \ \ \  <  \ \ 1  ;   \ \  \epsilon >\approx 0
   \ee

 
 \noindent   (\ref {RappSommeNegGRH} ) and   (\ref {RappSommePosGRH} ), if $ \epsilon>0 $ is fixed and $ p_{max}\rightarrow \infty$ goes to zero, because if $\epsilon>0$ and $k$ i.e. $p_{max} \rightarrow \infty$, then $p_{max}^{- \epsilon}\rightarrow 0$ in (\ref{RatioPrototype}).

\noindent In other words, in oder to clarify. 

\noindent Starting from following identity (we have $\forall p_{max} <\infty$):

  \be \label {Id} -
  X (t,\epsilon>\approx0 , p_{max}<\infty ) =\left \{
  \frac{    \sum_{+}^{p_{max} < \infty}    (t,\epsilon >\approx0)  }
  {  \sum_{+}^{p_{max} < \infty}    (t,\epsilon=0)  }  \right \}
    \left[  \sum_{+}^{p_{max} < \infty}    (t,\epsilon=0)  \right]
  - \ee

      $$\left\{
  \frac{  \sum_{+}^{p_{max} <\infty}    (t,\epsilon >\approx0)  
  +X(\epsilon >\approx0 ,p_{max})  }
  {    \sum_{+}^{p_{max} <\infty}    (t,\epsilon=0)  
  +X (\epsilon=0 , p_{max} )   }  
  \right\} 
   \left[      \sum_{+}^{p_{max} <\infty}    (t,\epsilon=0)  
  +X (\epsilon=0 , p_{max} )    \right]=
  $$
  $$ =\rho_+(t,\epsilon,p_{max})\left[  \sum_{+}^{p_{max} < \infty}    (t,\epsilon=0)  \right]-
  \rho_-(t,\epsilon,p_{max})  \left[      \sum_{+}^{p_{max} <\infty}    (t,\epsilon=0)  
  +X (\epsilon=0 , p_{max} )    \right]
  $$
  
\noindent And considering that,  from (\ref{RatioUnequalEps}) , and, (\ref{Xdef}) we have:

  \be \label {R+R-}
    \rho_-(t,\epsilon,p_{max})=  \rho_+(t,\epsilon,p_{max})
  \left( 1+\frac {X(t,\epsilon >\approx0 ,p_{max}) }{ \sum_{+}^{p_{max} <\infty}    (\epsilon >\approx0)}   \right)
   \left( 1-\frac {X(t,\epsilon =0 ,p_{max}) }{ \sum_{+}^{p_{max} <\infty}    (t,\epsilon =0)}   \right)=
  \ee
  \be \label {R+R-2}
  \rho_+(t,\epsilon,p_{max})
  \left( 1+
  \frac {X(t,\epsilon >\approx 0 ,p_{max}) }    { \sum_{+}^{p_{max} <\infty}    (t,\epsilon >\approx0)} -
  \frac {X(t,\epsilon =0 ,p_{max}) }{ \sum_{+}^{p_{max} <\infty}    (t,\epsilon =0)}  - 
   \frac {X(t,\epsilon >\approx 0 ,p_{max}) }    { \sum_{+}^{p_{max} <\infty}    (t,\epsilon >\approx0)} \times
  \frac {X(t,\epsilon =0 ,p_{max}) }{ \sum_{+}^{p_{max} <\infty}    (t,\epsilon =0)}  
  \right)
  \ee

  \noindent Now, what is critical is the cancellation of  the big  contribution of $ \left[  \sum_{+}^{p_{max} < \infty}    (t,\epsilon=0)  \right]$  in second term of identity  (\ref{Id}).  When eliminated $\rho_+(t,\epsilon,p_{max})\left[  \sum_{+}^{p_{max} < \infty}    (t,\epsilon=0)  \right]$, apart lesser terms, what is left is:
  $$-
  X (t,\epsilon>\approx0 , p_{max}<\infty ) =-
   \rho^+(t,\epsilon,p_{max})
X (\epsilon=0 , p_{max} )  -
$$
$$
 \left(  \frac {X(t,\epsilon >\approx 0 ,p_{max}) }    { \sum_{+}^{p_{max} <\infty}    (t,\epsilon >\approx0)} 
-
  \frac {X(t,\epsilon =0 ,p_{max}) }{ \sum_{+}^{p_{max} <\infty}    (t,\epsilon =0)}   \right) \left[      \sum_{+}^{p_{max} <\infty}    (t,\epsilon=0)  
  +X (\epsilon=0 , p_{max} )    \right]
  $$
  \noindent But 
  $$\left(  \frac {X(t,\epsilon >\approx 0 ,p_{max}) }    { \sum_{+}^{p_{max} <\infty}    (t,\epsilon >\approx0)} 
-
  \frac {X(t,\epsilon =0 ,p_{max}) }{ \sum_{+}^{p_{max} <\infty}    (t,\epsilon =0)}   \right) 
  \times \sum_{+}^{p_{max} <\infty}    (t,\epsilon=0)  =
  $$
  $$\left(  \frac {X(t,\epsilon >\approx 0 ,p_{max}) }    {   \rho^+(t,\epsilon,p_{max})} -
  \frac {X(t,\epsilon =0 ,p_{max}) }{ 1}   \right)=0$$

  \noindent 
  Then  for  (\ref {R+R-2})  for $p_{max}$  big enough   and $\epsilon$   small enough, but positive, we  can consider  $  \rho_-(t,\epsilon,p_{max}= \rho_+(t,\epsilon,p_{max}) =\rho(t,\epsilon,p_{max}$. So the identity (\ref{Id} ) for conditions stated above  becomes: 
  $$ -X (t,\epsilon>\approx0 , p_{max}<\infty ) =-
  \rho(t,\epsilon >\approx0,p_{max})
X (\epsilon=0 , p_{max} )    
  $$

  \noindent  Where $\rho $ is for $\rho^\pm$. Lemma 6  thesis through (\ref  {ApproachFromBelow} )  is so  proved, i.e.   $ \rho(t,\epsilon >\approx0,p_{max})<1$ approaching 1 from below.

  $$  \ Corollary \  of \  Lemma \ 6$$

 \noindent  
 In Lemma 6 we considered $p_{max}=(p^*)^m$ big but finite.
  From (\ref {WithLogInt} ) 
 we can write for $\epsilon\rightarrow 0^+$:
  
    \be \label {L6ALContinuo}
    \frac{1
 }
 { \epsilon}
\left [
\frac {\Delta \angle [ L_{EP}(t,\epsilon>0,\chi) ] }{\Delta t} -\frac {\Delta \angle [ L_{EP}(t,\epsilon=0,\chi) ] }{\Delta t} 
\right]_{ \ p_{max} \rightarrow \infty,  \ p^* \ big  } 
   \rightarrow 
   \left[ \frac {\partial^2 \angle[L(s,\chi_{
    primitive}) )]  } { \partial \epsilon\partial t} \right]_{ \epsilon =0}  \ge 0
   \ee
$$PROOF.$$

\noindent 
The difference in square brackets is between mean values in $t \pm \frac{\Delta t}{2}$ (\ref {Estremi})
of same function with different $\epsilon$, that in  $|t-t^*_k|>\frac{2 \pi}{\ln(p^*)} \ $ , for $  \ \epsilon \rightarrow 0$ tends to a smooth  function  in  $t$ :i.e. $- \ln \left(       \sqrt{    \frac { q t } { 2 \pi}      }    \right)  $
 \footnote { i.e. the value (\ref {WithLogInt} )  must assume for Lemma 2 for $\epsilon=0$)}. 
  So, for $p^*$ big the difference of the mean in the intervals  (\ref{Estremi}) is close to the function on the right in (\ref{L6ALContinuo}). Besides (\ref{LiDeltaXEps_2_GRH}) holds for every big $p^*$,  otherwise we could find a $p_{max}=(p^*)^m$  that makes Lemma 6  not true. But this is impossible.
  
  \noindent  So (\ref{L6ALContinuo}) must be true.

  |$$ END \ of \ PROOF$$
  
  $$ Lemma \ 7$$

   \noindent For odd primitive characters,for $t\ne t^*_k$,
   we have :
   $
    \epsilon  \times \left[\frac {\partial^2 \angle[L(t,\epsilon=0,\chi_{primitive})]}{ \partial \epsilon \partial t}\right]_{\epsilon=0} +\
     \frac {\epsilon}{4t^2} >0 $. So we conclude that:
     
     \be \label {whereWePointed}
      \left[\frac{\partial^2  \angle[\xi(s,\chi)] }{\partial \epsilon \partial t} \right]_{\epsilon \approx 0 \ , \ t \ne t^*_k} >0
\ee  
  
  \noindent PROOF:  From   (\ref {PhaseVarXiAndL_}), and , corollary of Lemma 6 it is self evident.


 \subsection {LEMMA 8 :    
 $ \left[  \frac{\partial   \mathcal{L}  [ \xi (t,\epsilon,\chi_{odd \ primitive}) ]  )}{\partial \epsilon} \right]_{\epsilon=0}  > 0  \  ,  \forall t \ne t^*_k \  ,  \forall | t | > T_{Asymp}(\alpha)   $ } \label {TheoOne}

\ 

 \noindent For $T_{Asymp}(\alpha) $ see  (\ref{Asympt} ). Referring to  (\ref{EpsZero}) (where $\eta(s,\chi)$ has been defined to be real on critical line) we have:

 \be \label {THEO1}
 \forall t \ \ with \ \   \eta(t,\chi) \ne 0 \rightarrow \left(  \frac{d \eta(t,\chi)}{dt} \right)^2- \eta(t,\chi) \times \frac{d^2 \eta(t,\chi)}{dt^2} >0 \ \ where  \ \ \eta(t,\chi)= \eta(1/2+it,\chi ), \ and \ \chi_{odd \ primitive}
 \ee 

$$PROOF$$ 
 \noindent From definition (\ref {AngMomDef}), and, (\ref{EpsZero}):
  $$\mathcal{L} [ \eta (s,\chi) ] =
 \mathcal{L} [ \xi (s,\chi) ]  = | \xi(s,\chi))|^2  \times \frac{\partial \angle[\xi(s,\chi))]}{\partial t} =
   | \eta(s,\chi))|^2  \times \frac{\partial \angle[\eta(s,\chi))]}{\partial t} 
    \ \ : \  \ \chi_{primitive} %
 $$

 \noindent
  From  (\ref {Uguale0Tris_} ) , as  for Lemma 2:
 
 $$\left[\frac{\partial \angle[L(s,\chi_{primitive})]}{\partial t}\right]_{\epsilon=0}= - \ln \left( \sqrt{\frac{t q }{2 \pi}} \right)
 \ ; \ \forall t \ne t^*_k $$
 \noindent we have that, from (\ref{XiForLfunc} ), (\ref{DerAngFattXiChiPrim}),  (\ref{3Cases}) and  (\ref{PhaseVarXiAndL_}) : 

$ \left[\frac{\partial^2  \angle[\xi(s,\chi)] }{ \partial \epsilon\partial t} \right]_{\epsilon = 0} \  =
   \left[\frac {\partial^2 \angle[L(s,\chi_{primitive})]}{ \partial \epsilon \partial t}\right]_{\epsilon=0}
    +\left[
  \frac {\partial^2 \angle\left[ \left( \frac{\pi}{q}   \right)^{-\frac{s+\alpha_1}{2}} \Gamma  \left( \frac{s+\alpha}{2} \right)    \right]}{ \partial \epsilon \partial t} 
  \right]_{\epsilon=0}
  =$
    
    \be \label {Uno}\left[\frac {\partial^2 \angle[L(s,\chi_{primitive})]}{ \partial \epsilon \partial t}\right]_{\epsilon=0} 
    \quad    \pm \frac {1}{4t^2} \ \ ; \   \ \forall t \ne t^*_k \ \ ; + \ odd \ \ ; - \ \ even  \ \ \chi  
    \ee

\noindent On the other hand:


 $$
\left [\frac{ \partial \mathcal{L} [ \xi (s,\chi) ] }{\partial \epsilon}  \right ]_{\epsilon=0}
= \left\{ \frac {\partial | \xi(s,\chi) |^2}{\partial \epsilon} \times \frac{\partial \angle[\xi(s,\chi))]}{\partial t}  \right\} _{\epsilon=0} +
\left[
 | \xi(s,\chi))|^2  \times \frac{\partial^2 \angle[\xi(s,\chi))]}{\partial \epsilon \partial t}  
 \right]_{\epsilon=0} \ \ : \  \ \chi_{primitive} %
 $$
 \noindent  but, the expression in braces $ \left\{ \frac {\partial | \xi(s,\chi) |^2}{\partial \epsilon} \times \frac{\partial \angle[\xi(s,\chi))]}{\partial t}  \right\} _{\epsilon=0} =0  $ because for $t\ne t^*_k$:  $ \left\{  \frac{\partial \angle[\xi(s,\chi))]}{\partial t}  \right\} _{\epsilon=0} =0  $, and  $ \left\{ \frac {\partial | \xi(s,\chi) |^2}{\partial \epsilon}  \right\} _{\epsilon=0}   $ is finite  ( as $  \xi(s,\chi)$ is holomorphic). The same is true if we replace $\xi$ by $\eta$.

 
 \noindent   %
 So, referring only  to odd  primitive characters ( see (\ref{3Cases}) ), and,  using
  corollary to Lemmas 6 (\ref{L6ALContinuo})  that states $  \left[ \frac {\partial^2 \angle[L(s,\chi_{odd \ primitive}) )]  } { \partial \epsilon\partial t} \right]_{ \epsilon =0} \ge 0$, then, (\ref{Uno}), 
  implies that:  $\left[\frac{\partial^2 \angle[\eta(s,\chi))]}{\partial \epsilon \partial t}  
 \right]_{\epsilon=0} =
  \left[\frac{\partial^2 \angle[\xi(s,\chi))]}{\partial \epsilon \partial t}  
 \right]_{\epsilon=0}  > 0$  %
 (for $t>T_{Asymp}(\alpha)$, and $t \ne t^*_k$). 
 So also (\ref{EpsZero}):

$$ \left\{  \frac {\partial}{\partial \epsilon}
 \mathcal{L} [ \xi(s,\chi_{primitive}) ]\right\}_{\epsilon=0} =
 \left\{  \frac {\partial}{\partial \epsilon}
 \mathcal{L} [ \eta(s,\chi_{primitive}) ]\right\}_{\epsilon=0}
 = $$

\be \label {ProBombieriP6ForL}
  \left\{ 
 |\eta(s,\chi_{primitive})|^2  \frac {\partial}{\partial \epsilon}\frac {\partial}{\partial t}  \angle [\eta(s,\chi_{primitive})]  \right\}_{\epsilon=0} > \delta >0 \ \ ; \ \ \forall t: t \ne t^*_k
 \ee

\noindent Where $\delta= |\eta(t,\epsilon=0,\chi_{primitive})|^2 \times   \left( 
  \frac{1}{4t^2}  \right)>0$. 
  
  $$END \ of  \ PROOF$$
   
  {\it LEMMA 9 : $ \left[  \frac{\partial   \mathcal{L}  [ \xi (t,\epsilon,\chi_{odd \ primitive}) ]  )}{\partial \epsilon} \right]_{\epsilon=0}  > 0 $  also   $  \forall t = t^*_k , \  if \ \epsilon^*_k >0\  ,  \forall | t | > T_{Asymp}(\alpha)   $ 
   }.

\noindent   Suppose      $\eta(t^*_k,\epsilon^*_k>0,\chi) = 0$ , so  $\eta(t^*_k,\epsilon=0,\chi) \ne 0$.
 \noindent If   we  shrinks $|t-t^*_k| <\frac{ 2 \pi}{\ln(p^*)} $ interval with$[ \{p_{max} \rightarrow \infty\}  p^*\rightarrow \infty]$,  then (\ref  {ProBombieriP6ForL}) is valid everywhere around $t^*_k$ for $\epsilon=0$. But  also in $t^*_k$, because with $p^*$ big enough the validity of (\ref  {Uguale0Tris_}  ) can be reached close to the phase discontinuity line ($(t=t^*  \ , \  \epsilon <\epsilon^*)$)  as we want, and,  see section  \ref {Lemma5}, %
  $[\eta'(t,\chi)]^2 -\eta(t,\chi)  \eta''(t,\chi) =\left\{  \frac {\partial}{\partial \epsilon}
 \mathcal{L} [ \xi (s,\chi) ]\right\}_{\epsilon=0} $ is a continuous function.
 %

%


  \noindent
 \subsection {THEOREM 1  : GRH is true for the odd primitive L-functions  except for at most finitely many exceptions whose imaginary parts live in  a finite  interval $|t|<T_{Asymp}(\alpha) $ 
 }
      \label {CONJ}

      
      \noindent In \cite[p.~6] {Bombieri:2000cz} is reported a statement as equivalent to R.H. for $\zeta(s)$: { \it ''The  Riemann hypothesis is equivalent to the statement that all local maxima
of  { \bf $ \xi(t,\epsilon=0)$ } are positive and all local minima are negative.''}.  Let us see why.  
 Suppose we are close to extremal points. Then $\Re  [\xi(s )]$, is an harmonic function, obeying to Laplace equation $ \frac { \partial ^2\Re[ \xi  (\frac { 1}{2}+\epsilon +it  ) ]}{\partial \epsilon^2}
   +\frac{\partial^2 \Re[ \xi  (\frac { 1}{2}+\epsilon +it  ) ]}{\partial  t^2} =0 $,   with null total curvature. See (\ref {Cauchy-Riemann}).  
  As   $ \frac { \partial \Im[ \xi  (\frac { 1}{2}+\epsilon +it  ) ]}{\partial \epsilon}
   = -\frac{\partial \Re[ \xi  (\frac { 1}{2}+\epsilon +it  ) ]}{\partial  t}  $, 
  and  $ \Im[ \xi  (\frac { 1}{2} +it  ) ]=0$, then  extremal points on critical line (that are saddle points)
  are first candidate to look for off-critical line zeros.  
 So, at extremal points (relative maxima end minima), if the curvature is toward   $[\epsilon ,t]$ plane along $t$, it will be in the opposite direction along $\epsilon$. Simmetrically if the curvature is toward  $[\epsilon ,t]$ plane along $\epsilon$,  it will be in the opposite direction along $t$. 
Only in second case off-critical line zeros are possible because we have for example $\Re[ \xi(t_0,\epsilon_0=0)]>0$ and $\Re[ \xi(t,\epsilon_0+ \Delta \epsilon>0)] <\Re[\xi(t_0,\epsilon_0=0)]$, i.e. we are tending toward the plane $[\epsilon ,t] $  along $\epsilon$. Instead in first case we have for example $\Re[ \xi(t_0,\epsilon_0=0)]>0$ and $\Re[ \xi(t,\epsilon_0+ \Delta \epsilon>0)] >\Re[\xi(t_0,\epsilon_0=0)]$, i.e. we are increasing distance from plane $[\epsilon ,t] $  along $\epsilon$. In first case the curvature along $t$ must be negative, i.e. $\frac{d^2 \xi(t,\epsilon=0)}{d t^2} <0$ so $\xi(t,\epsilon=0) \frac{d^2 \xi(t,\epsilon=0)}{d t^2} <0$.  In second case we have instead  $\frac{d^2 \xi(t,\epsilon=0)}{d t^2} >0$ so $\xi(t,\epsilon=0) \frac{d^2 \xi(t,\epsilon=0)}{d t^2} >0$. Notice the sign of the product   $\xi(t,\epsilon=0) \frac{d^2 \xi(t,\epsilon=0)}{d t^2}<0 $ is invariant also if we choose $\Re[ \xi(t_0,\epsilon_0=0)]<0$  if all local maxima
of  { \bf $ \xi(t,\epsilon=0)$ } are positive and all local minima are negative.
This means that if  { \bf all} extremal points of  $\Re[ \xi(t_0,\epsilon_0=0)]$ comply with this rule, 
 %
no off-critical line zero is possible because 
 $ \frac{\partial \Im[ \xi  (\frac{ 1}{2}+\epsilon +it  ) ]}{\partial \epsilon} = -\frac{\partial \Re[ \xi  (\frac{ 1}{2}+\epsilon +it  ) ]}{\partial  t} \ \ $. So $|\Im[ \xi  (\frac{ 1}{2}+\epsilon +it  ) ] |$ grows with $\epsilon$ far from $\Re[ \xi  (\frac{ 1}{2} +it  ) ]$ extremal points.  We do not say that the set of points of $(t,\epsilon)$ plane where  $\Im[\xi(s)]=0$ are perpendicular to critical line at $t_0$ where $ \Re[ \xi  (\frac{ 1}{2} +it_0  )]$ is extremal, but, the truth is not far from it \cite [p.~3] {Giovanni Lodone 2021} , and \cite [p.~9,fig 6] {Giovanni Lodone 2024}.

\noindent  General facts of holomorphic behavior are present also for  $\eta(t,\epsilon,\chi_{primitive})$, besides it  is real on critical line  as $\xi(t,\epsilon)$ so \cite[p.~6] {Bombieri:2000cz} equivalence must be valid also for it.  Besides for Lemma 8 and Lemma 9 we have:
$[\eta'(t,\chi)]^2 -\eta(t,\chi)  \eta''(t,\chi) =\left\{  \frac {\partial}{\partial \epsilon}
 \mathcal{L} [ \xi (s,\chi) ]\right\}_{\epsilon=0} >0  \ , \forall t : \   \eta(t,\chi) \ne 0$. But this means that at extremal points, i.e.  $ \frac  {d \eta(t',\chi )}     {dt}   =0$,   %
local maxima are positive and local minima negative because:

  \be \label {maxminrel}
  \eta(t',\chi ) \left[   \frac{d^2 \eta(t',\chi )}{dt^2}  \right] < 0
\ee

\

\noindent So for $ \eta(t,\epsilon=0,\chi) = 0$   the zero ( or zeros if double) are on  critical line. But also for 
 $ \eta(t,\epsilon=0,\chi) \ne0 $   expression (\ref {maxminrel})  implies that relative maxima and minima must comply with no zeros off-critical line. So,  for Lemma 8 and Lemma 9 ,  RH is granted at least  in $|t|>T_{Asymp}(\alpha) $,  see  (\ref{Asympt} ). So for a given $\L(s, \chi_{odd \ primitive})$ only a finite numbers of off-critical line zeros  exist in $|t|<T_{Asymp}(\alpha) $, because the function is analytic in $ |t| < T_{Asymp}(\alpha)$, if,  there were infinitely many zeros there, it would be an accumulation point and the function would be identically zero, a contradiction.

\section {  Toward the proof  $\forall t$ and $\forall$ modulus $q$ of primitive (odd)  characters}
 \label {TowardGRH}

\noindent With argument in appendix \ref  {GammaDer} and Lemmas of previous section, %
we proved RH for odd primitive characters in $|t| >T_{Asymp}(\alpha) %
$. In other words if  $|t| >T_{Asymp}(\alpha)$ we are allowed to write (\ref {DerAngFattXiChiPrim})  i.e. $
   \frac {\partial \angle\left[ \left( \frac{\pi}{q}   \right)^{-\frac{s+\alpha}{2}} \Gamma  \left( \frac{s+\alpha}{2} \right)    \right]}{\partial t} = \ln\left( \sqrt{\frac{t q}{2 \pi}}\right) + O(t^{-2})$, and  $O(t^{-2}) \rightarrow +(-1)^{\alpha+1} \frac{\epsilon}{4 t^2}$ for $\epsilon \approx 0$ ($\alpha=0$ even, and $ \alpha=1$ odd characters).%

\noindent Here we find  %
 same results of appendix \ref  {GammaDer} without  Stirling formula like in  expression (\ref {Pi&Gamma3}). The purpose is  to evaluate  (\ref {DerAngFattXiChiPrim})  in  the $t=0$ neighborhood.
\noindent The starting point is the   formula (\cite[p.~8] {Edwards:1974cz}) )

\be \label {GaussWeistrassFor} 
\Gamma(z)=\lim_{N \to \infty} \frac{N!  \ \ (N+1)^z}{z(z+1) . . .(z+N)} %
\ee

 \noindent  where %
  $z= \frac{s + \alpha}{2}= \frac{1/2+\epsilon+ i t + \alpha}{2}$ (see (\ref {XiForLfunc}), and (\ref{FattorePerXIAChi}´) ) .

\noindent%
For $N\rightarrow \infty$ we sum and subtract: $\sum_{n=1}^N   \frac{1}{n} \rightarrow \ln(N) + \gamma \approx \ln(N+1) +\gamma  $,  where $\gamma= 0.577216 . . \ \ $ is the Euler - Mascheroni constant.

\be \label {Fase} 
\Im[\ln(\Gamma(z))]=
\angle[ \Gamma(z)] = \angle[ z\ln(N+1)] -\angle[z]-\angle[z+1] . . .-\angle[z+N] =
\ee

 $$
 \left( \ln(N+1)  -  \sum_1^N   \frac{ 1}{ n}\right)  \frac{ t}{2} -\arctan \left(      \frac {\frac{t}{2} }{\frac {1/2+ \epsilon+ \alpha}{2}}  \right)+
\sum_1^N \left \{ \frac{ t}{2 n}  -   \arctan \left(      \frac {\frac{t}{2n} }{ 1+  \frac {1/2+ \epsilon+ \alpha}{2n}}  \right)  \right \}=
$$

 $$
  -\gamma \frac{ t}{2} -\arctan \left(      \frac {\frac{t}{2} }{\frac {1/2+ \epsilon+ \alpha}{2}}  \right)+
\sum_1^N \left \{ \frac{ t}{2 n}  -   \arctan \left(      \frac {\frac{t}{2n} }{ 1+  \frac {1/2+ \epsilon+ \alpha}{2n}}  \right)  \right \}
$$

Let us fix $N$ in (\ref {GaussWeistrassFor} ), and,  in (\ref{Fase}), then if we take derivative  with respect to $t$:

$$\frac {\partial \angle \left[\Gamma\left( \frac {s+\alpha}{2} \right) \right]}{\partial t}  = 
   \frac{ -\gamma}{2} - \left(      \frac {1}{   1/2+ \epsilon+ \alpha   }  \right) \left(  \frac{1}{1+\left(\frac{t}{ 1/2+ \epsilon+ \alpha}  \right)^2} \right)+
$$

\be \label {DerGammaSuTcircaZero}+
\sum_1^{N\rightarrow \infty} \left \{ 
\frac{ 1}{  2 n}  -    \left(      \frac {1}  { 2n+  1/2+ \epsilon+ \alpha}   \right) \left(  
  \frac{1}{1+   \left(   \frac  {t}{  2n+  1/2+ \epsilon+ \alpha  }\right)^2} \right)
 \right \}
\ee

 \noindent Without Stirling formula as in  (\ref {Asympt} )  , the equivalent relation   of $\ln\left[ \sqrt{\frac{t q }{2 \pi}} \right] $ in (\ref   {Uguale0Tris_} ) or  (\ref {PhaseVarXiAndL_}), is (\ref {DerAngleDeT}). See fig.  (\ref{Q3}), and  fig.(\ref{Q5}) for respectively odd and even characters and fig \ref{ForZfunc} for $\zeta(s)(s-1)$ .

\be \label {DerAngleDeT}
   \frac {\partial \angle\left[ \left( \frac{\pi}{q}   \right)^{-\frac{s+\alpha_1}{2}} \Gamma  \left( \frac{s+\alpha}{2} \right)    \right]}{\partial t} =
  \frac {\partial \Im \left\{  \ln  \left[ \left( \frac{\pi}{q}   \right)^{-\frac{s+\alpha_1}{2}} \Gamma  \left( \frac{s+\alpha}{2} \right)    \right]  \right\}}{\partial t} 
   =
   \frac{1}{2} \ln\left(\frac{q}{\pi}\right) +\frac {\partial \angle \left[\Gamma\left( \frac {s+\alpha}{2} \right) \right]}{\partial t}  
   \ee

\begin{figure}[]
\begin{center}
\includegraphics[width=1.0\textwidth]{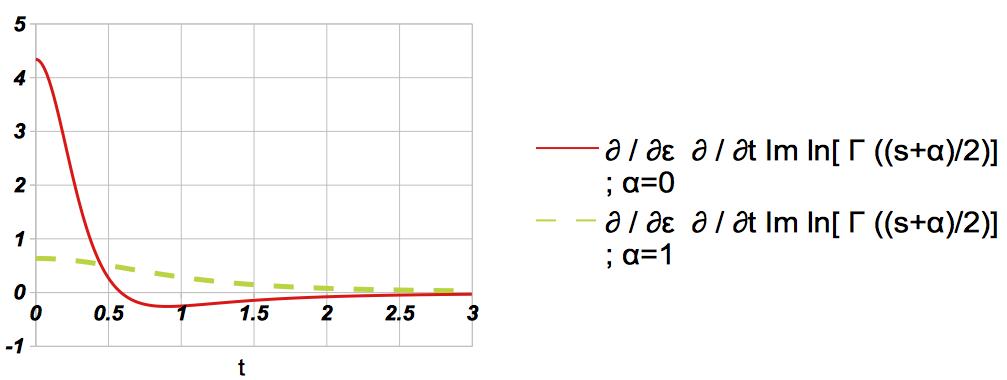} 
\caption{\small {\it  
Derivative $\frac{\partial}{\partial \epsilon} $ of (\ref  {DerAngleDeT}). With $\alpha=0$ even,  $\alpha=1$  odd characters.
 The point of crossing with horizontal axis in case $\alpha=0$ is: $0.585 <t_{cross}<0.588$. The curves are independent from  congruence modulus $q$.
}}
 \label {DerEpsDerTDiGam}
\end{center}
\end{figure}

\subsection {Odd characters ($\alpha=1$), $\forall t$, and $\forall q$}
\label {GRHOdd}

In fig \ref {Q3} is plotted  (\ref{DerAngleDeT}) For $q  = 3 ;  \alpha =1$. For $q'>q$  we must shift the curve upward. For example, with $q'=11$,  adding $\frac{\ln(q'/\pi)-\ln(3/\pi)}{2}  >0.64$.
 So   for $q \ge 11$,   (\ref {DerAngleDeT}) is positive $\forall t$. %
Then we are in situation of Theorem 1,  
\noindent because   (\ref{DerAngleDeT}) is always $>0$, and, 

  \noindent $\frac{\partial}{\partial \epsilon} \frac {\partial \Im \left\{  \ln  \left[ \left( \frac{\pi}{q}   \right)^{-\frac{s+\alpha_1}{2}} \Gamma  \left( \frac{s+\alpha}{2} \right)    \right]  \right\}}{\partial t}= \left\{ \frac{\partial}{\partial \epsilon}
\frac{\partial}{\partial t} \angle \left[ \Gamma \left(\frac{s + \alpha}{2 }  \right) \right] \right\}_{\alpha=1, } >0$ too. See fig.  \ref{DerEpsDerTDiGam} ).
 So, as in previous
section  we conclude that for   $q \ge 11 ;  \alpha =1$:  RH is true $\forall t$.   The sufficient condition (\ref {SuffOdd}) 

  \be \label {SuffOdd}
t :\left[ \frac{1}{2} \ln\left(\frac{q}{\pi}\right) +\frac {\partial \angle \left[\Gamma\left( \frac {s+\alpha}{2} \right) \right]}{\partial t}\right]_{\alpha=1} >0 \  \  ; 
\ \  \left\{ \frac{\partial}{\partial \epsilon}
\frac{\partial}{\partial t} \angle \left[ \Gamma \left(\frac{s + \alpha}{2 }  \right) \right] \right\}_{\alpha=1, }>0 \
\ee
    
  \noindent      can be used also locally if the zeros %
 happen  at $t$ values in which  same situation   (\ref {SuffOdd}) is locally met.
\noindent The zeros  of  $L(s,\chi_{primitive})$  distribute %
like the peaks of the expression (\ref{DerLAngleSuTL}                                                                                                                                                                                                                                                                                                                                                                                                                                                                                                                                                                                                             ) that we know  occurs  exactly at the zero of $L(s,\chi_{primitive})$. Irrespective if RH is verified or not.

\begin{table}[bht]
\begin{center}
\begin{tabular}{c|c|c|c|c|c|}
Characters phase     & $ n=0 $    & $n=1$     &      $n=2$     &     $n=3$     & $n=4$        \\
\hline
$\angle[ \chi_0(n)]$ & undefined &0& 0 & 0 & 0 \\
\hline
$\angle[ \chi_1(n)]$ & undefined &0 & $\pi/2$ &  $(3/2) \pi $ & $\pi  $\\
\hline
$\angle[\chi_2(n)]$ & undefined & 0& $\pi$ & $\pi$ & 0\\
\hline
$\angle[\chi_3(n)]$ & undefined & 0& $(3/2) \pi$& $\pi/2$ & $\pi$ \\
\hline
\end{tabular}
\end{center}
\caption{\small Phase of characters components for $L(s,\chi(n) )$ .Arithmetic congruence  modulus $q= 5$. In column $n=0$ the  modulus  ( we mean : $\sqrt{\Re^2[\chi] +  \Im^2[\chi]  }$ ) is zero. Elsewhere it is 1. When $q$ is prime only  principal  and  primitive characters are present. See \cite [p.~168]  {Apostol:1976}.Notice in each, but the first, row   $\sum_{h<q : gcd(h,q)=1}e^{ i \ \angle[\chi_h]}=0$ for $\chi$ not principal see \cite [p.~256]  {Apostol:1976}. Table  loaded from \cite {CharactersTablesTillQ62} and adapted for (\ref{DerLAngleSuTL}                                                                                                                                                                                                                                                                                                                                                                                                                                                                                                                                                                                                              ) .
}
\label {tab1}
\end{table}

\begin{table}[bht]
\begin{center}
\begin{tabular}{c|c|c|c|c|c|}
 q     & Odd chrs Suff. Cond. (\ref{SuffOdd})    &  closest odd chrs correlation  peak distance  from $t=0$       \\
\hline
$3$ & $| t | > 2$ &$| t | > 8$ \\
\hline
$4$ & $| t | > 1.5$    & $| t | >6$ \\
\hline
$5$ & $ | t |>1.25$ & $| t | >4$\\
\hline
$7$ & $| t | >0.75$& $| t | >2$ \\
\hline
$8$ & $| t | >0.5$& $| t | >2.5$ \\
\hline
$9$ & $| t | >0.25$& $| t | >2.4$ \\
\hline
\end{tabular}
\end{center}
\caption{\small  For odd characters,  see fig \ref  {Q3},  sufficient condition to apply the argument of  Theorem 1 pag. 15, %
  is $\left[\frac {\partial \angle \left[\Gamma\left( \frac {s+\alpha}{2} \right) \right]}{ \partial t}\right]_{\epsilon=0} +\frac{\ln(q/ \pi)}{2} > 0$ with $\alpha=1$. Here are reported (second column)  his crossing with horizontal axis in fig  \ref  {Q3}  at various $q'$, while in third column the $t$ interval without peaks (zeros).  For $q=6$ and $q=10$  no primitive characters are present \cite  {Mathar}.
}
\label {tabODD}
\end{table}

\begin{figure}[]
\begin{center}
\includegraphics[width=1.0\textwidth]{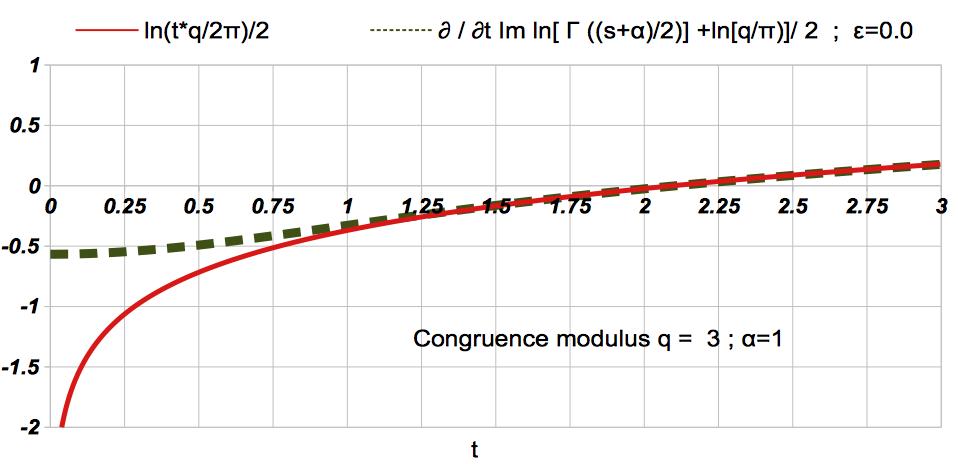} 
\caption{\small {\it  
%
Plot of (\ref{DerAngleDeT}). For $q'>3$  we must shift the curve by the quantity $\frac{\ln(q'/\pi)-\ln(3/\pi)}{2}$. For $q \ge 11$ the curve is positive $\forall t$.Notice that the  asymptotic result (\ref{Asympt} )  is well verified at least  from $|t|>3$. Here  $T_{Asymp}(\alpha) \approx 3$.
}}
\label {Q3}
\end{center}
\end{figure}

\noindent For $q<11$  see for example  fig.  %
\ref     {Simmetries} .%
Here by direct verification,  we see that for $q=5$  all zeros are far from negative values of (\ref {DerAngleDeT}).  We see that  the minimum distance of the odd characters  ($\chi_1$ and $\chi_3$) from $t=0$ happen to be $\approx 4 >1.25 $ where (\ref{DerAngleDeT}) with $q=5$ crosses horizontal axis. See also tab. \ref {tabODD}. So RH is verified for $q=5$, $\chi_1$ and   $\chi_3$,  because at each  odd character zero we can locally apply argument of Theorem 1 
as sufficient condition  (\ref{SuffOdd}) is locally met. For even primitive character $\chi_2$  this is not applicable.

\noindent %
Considering data of  tab. \ref {tabODD}, built as for case $q=5$, RH is verified  for all primitive characters   with $3 \le q < 10$.

\begin{figure}[]
\begin{center}
\includegraphics[width=1.0\textwidth]{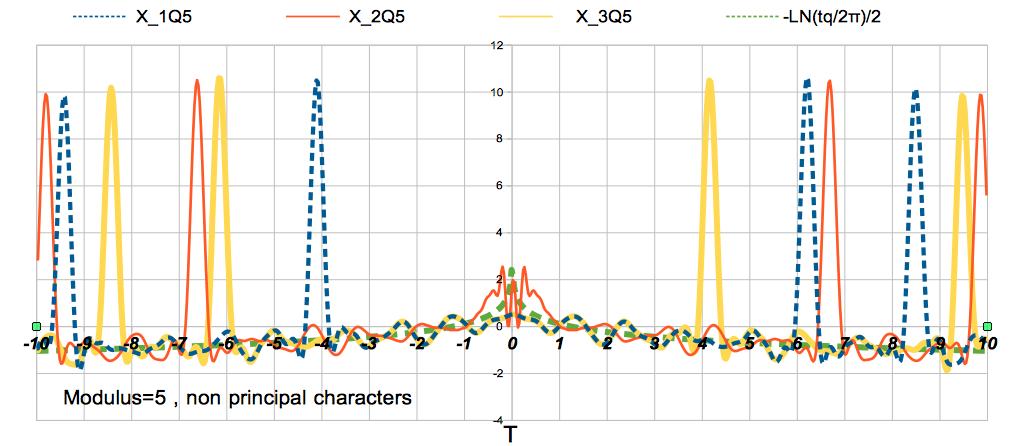} 
\caption{\small {\it   The peaks  symmetries of (\ref  {DerLAngleSuTL}                                                                                                                                                                                                                                                                                                                                                                                                                                                                                                                                                                                                              ) with  $\epsilon=0$, respect to real axis,  in not principal character case  $q =5$ ( see tab. \ref {tab1}) ,  with $p_{max}=p^*=158 \times 10^6$ ( here are used the first 8868881 primes in (\ref{DerLAngleSuTL} ) ).%
In wolfram demostration project - Dirichlet L-Functions and their zeros we can compute,  by zeroing real and imaginary parts,  the zeros  closer to $t=0$ of what  here are called X1 and X3,  the result is: $\pm9.443$, $\pm8.457   $, $\pm6.184,  \   \ \pm 4.133 \  $. Very comparable with results in figure.
}}
\label {Simmetries}
\end{center}
\end{figure}


\begin{figure}[]
\begin{center}
\includegraphics[width=1.0\textwidth]{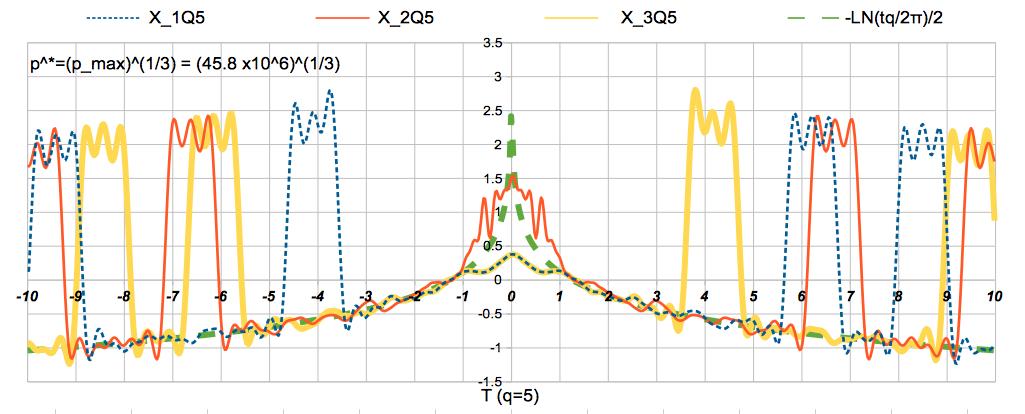} 
\caption{\small {\it   The peaks  symmetries of (\ref  {DerLAngleSuTL}                                                                                                                                                                                                                                                                                                                                                                                                                                                                                                                                                                                                              ) with  $\epsilon=0$, respect to real axis,  in  character case  $q =5$,  with $p_{max}=(p^*)^{3}=45.8 \times 10^6$ ( here are used the first 2763823 primes in (\ref{DerLAngleSuTL} )  ).%
The peaks are lower and fatter with respect to fig. \ref  {Simmetries}, but, the inter-peaks behavior is much more close to $\ln\left( \sqrt{\frac{q  t }{2 \pi}}  \right)$, or, better,  to (\ref {DerAngleDeT}), that follows different behavior close to $t=0$ if $\alpha=0$ (even chrs) or $\alpha=1$ (odd chrs).
}}
\label {Q5-3Camp}
\end{center}
\end{figure}

\begin{figure}[]
\begin{center}
\includegraphics[width=1.0\textwidth]{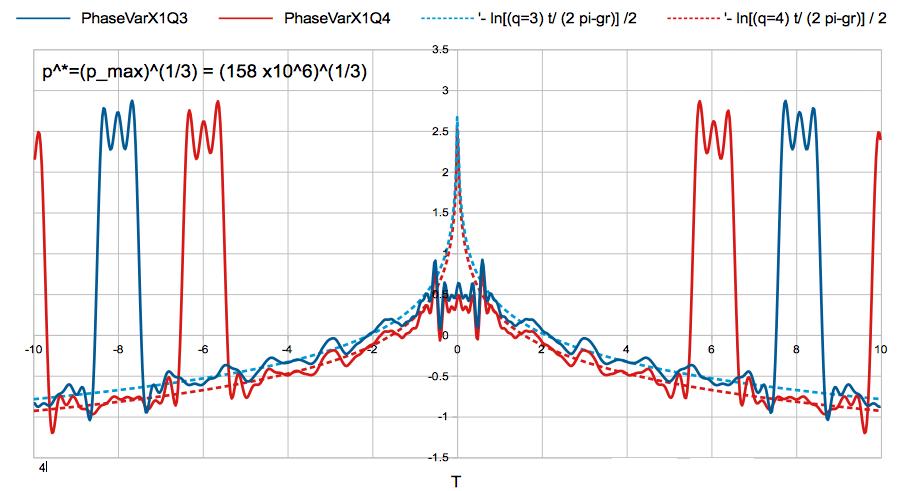} 
\caption{\small {\it   The peaks  symmetries of (\ref  {DerLAngleSuTL}                                                                                                                                                                                                                                                                                                                                                                                                                                                                                                                                                                                                              ) with  $\epsilon=0$, respect to real axis,  in  character case  $q =3$ and   $q =4$,  with $p_{max}=(p^*)^{3}=158 \times 10^6$ ( here are used the first 8868881 primes).%
For  $q=3$, and,  $\chi_1(n)=(0,1,-1)$  or  $\angle[\chi_1(n)]=(undefined, 0 , \pi)$ to comply with format of tab. \ref {tab1}.  For $q=4$,  and $\chi_1(n)=(0,1,0,-1)$, or,  $\angle[\chi_1(n)]=(undefined, 0 , undefined, \pi)$ to comply with format of tab. \ref {tab1}.
Notice close to $t=0$ the characteristic behavior of odd primitive L-functions. Compare with fig.  \ref {Q3}. Until now the only example of even primitive L-function is  phase variation of X2Q5  in fig. \ref  {Simmetries}, and fig. \ref {Q5-3Camp}. In this case compare the inter-peaks level with fig. \ref{Q5} and asymptotic formula (\ref {Asympt} ). 
}}
\label {Q3-Q4-m3}
\end{center}
\end{figure}

\noindent Putting all pieces together we can state { \bf that   $\xi( s,\chi_{odd \ primitive})$,  is RH compliant $\forall t$ and $\forall q$ }.

\begin{table}[bht]
\begin{center}
\begin{tabular}{c|c|c|c|c|c|c|c|c|c|c|c|c|c|}
Characters / n=      & $ 0 $    & $1$     &      $2$     &     $3$     & $4$    & $5$      & $6$  & $7$  & $8$  & $9$  & $10$  & $11$  & $12$  \\
\hline
$ \chi_0(n)$ & 0 &1 & 1 & 1 & 1 &1& 1 & 1 & 1 & 1 & 1 & 1 & 1\\
\hline
   $ \chi_1(n)$ & 0 & $1$ &$w$ & $w^4$ &  $w^2 $ & $-i$  & $w^5$ &$-w^5$ &$i$  &$-w^2$   &$-w^4$ &$-w$ &$-1$\\
   \hline
    $\chi_2(n)$ & 0 & $1$  &$w^2$ & $-w^2$ &  $w^4 $ & $-1$  & $-w^4$ &$w^2$ &$-1$  &$w^4$   &$-w^2$ &$w^2$ &$1$\\
    \hline
    $\chi_3(n)$ & 0& $1$  &$i$ & $1$ &  $-1 $ & $i$  & $i$ &$-i$ &$-i$  &$1$   &$-1$ &$-i$ &$-1$\\
    \hline
    $\chi_4(n)$ & 0& $1$  &$w^4$ & $w^4$ &  $-w^2 $ & $1$  & $-w^2$ &$-w^2$ &$1$  &$-w^2$   &$w^4$ &$w^4$ &$1$\\
    \hline
    $\chi_5(n)$ & 0 & $1$ &$w^5$ & $-w^2$ &  $w^4 $ & $-i$  & $w$ &$-w$ &$i$     &$w^4$ &$w2$  &$-w^5$   &$-1$\\
    \hline
    $\chi_6(n)$ & 0 & $1$  &$-1$ & $1$ &  $1$ & $-1$  & $-1$ &$-1$ &$-1$  &$1$   &$1$ &$-1$ &$1$\\
    \hline
    $\chi_7(n)$ & 0& $1$  &$-w$ & $w^4$ &  $w^2 $ & $i$  & $-w^5$ &$w^5$ &$-i$  &$-w^2$   &$-w^4$ &$w$ &$-1$\\
    \hline
    $\chi_8(n)$ & 0& $1$  &$-w^2$ & $-w^2$ &  $w^4 $ & $1$  & $w^4$ &$w^4$ &$1$  &$w^4$   &$-w^2$ &$-w^2$ &$1$\\
    \hline
     $\chi_{9}(n)$ & 0& $1$  &$-i$ & $1$ &  $-1$ & $-i$  & $-i$ &$i$ &$i$  &$1$   &$-1$ &$i$ &$-1$\\
    \hline
    $\chi_{10}(n)$ & 0& $1$  &$-w^4$ & $w^4$ &  $-w^2 $ & $-1$  & $w^2$ &$w^2$ &$-1$  &$-w^2$   &$w^4$ &$-w^4$ &$1$\\
    \hline
     $\chi_{11}(n)$ & 0 & $1$ &$-w^5$ & $-w^2$ &  $-w^4 $ & $i$  & $-w$ &$w$ &$-i$  &$w^4$   &$w^2$ &$w^5$ &$-1$\\
    \hline
\end{tabular}
\end{center}
\caption{\small Characters components for $L(s,\chi(n) )$ .Arithmetic congruence  modulus $q= 13$, $w=e^{i \pi/6}$. I Table  loaded from \cite {CharactersTablesTillQ62} .
}
\label {tab2}
\end{table}

\begin{table}[bht]
\begin{center}
\begin{tabular}{c|c|c|c|c|c|c|c|c|c|c|c|c|c|}
Characters / n=      & $ 0 $    & $1$     &      $2$     &     $3$     & $4$    & $5$      & $6$  & $7$  & $8$  &     \\
\hline
$ \chi_0(n)$ & 0 &1 & 1 & 0 & 1 &1& 0 & 1 & 1 &     \\
\hline
   $ \chi_1(n)$ & 0 & $1$ &$w$ & $0$ &  $w^2 $ & $-w^2$  & $0$ &$-w$ &$-1 $  & \\
   \hline
    $\chi_2(n)$ & 0 & $1$  &$w^2$ & $0$ &  $-w $ & $-w$  & $0$ &$w^2$ &$1$  &  \\
    \hline
    $\chi_3(n)$ & 0& $1$  &$-1$ & $0$ &  $1 $ & $-1$  & $0$ &$1$ &$-1$  &not primitive   \\
    \hline
    $\chi_4(n)$ & 0& $1$  &$-w$ & $0$ &  $w^2 $ & $w^2$  & $0$ &$-w$ &$1$  &   \\
    \hline
    $\chi_5(n)$ & 0 & $1$ &$-w^2$ & $0$ &  $-w $ & $-w$  & $0$ &$w^2$ &$-1$     &   \\
    \hline
   
\end{tabular}
\end{center}
\caption{\small Characters components for $L(s,\chi(n) )$ .Arithmetic congruence  modulus $q= 9$, $w=e^{i \pi/3}$. The table  is loaded from \cite {CharactersTablesTillQ62} .
}
\label {tab3}
\end{table}



\begin{figure}[]
\begin{center}
\includegraphics[width=1.0\textwidth]{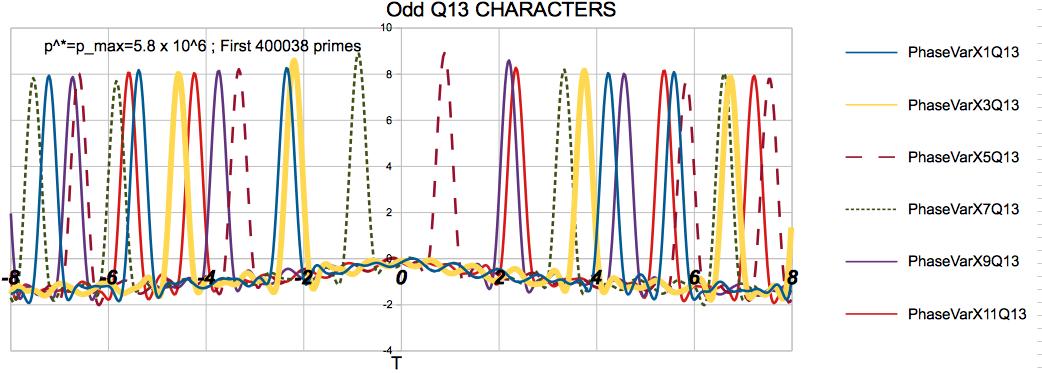} 
\caption{\small {\it   The peaks  symmetries of (\ref  {DerLAngleSuTL}                                                                                                                                                                                                                                                                                                                                                                                                                                                                                                                                                                                                              ) with  $\epsilon=0$, respect to real axis,  in  character case  $q =13$ (see tab. \ref{tab2}), and   ,  with $p_{max}=(p^*)=5.8 \times 10^6$ ( here are used the first 400038 primes  ).%
}}
\label {OddQ13}
\end{center}
\end{figure}

\begin{figure}[]
\begin{center}
\includegraphics[width=1.0\textwidth]{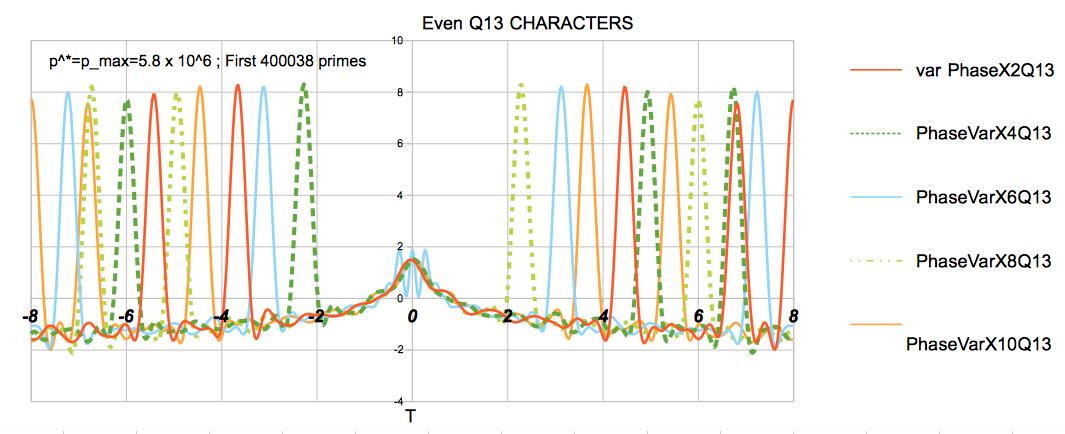} 
\caption{\small {\it   The peaks  symmetries of (\ref  {DerLAngleSuTL}                                                                                                                                                                                                                                                                                                                                                                                                                                                                                                                                                                                                              ) with  $\epsilon=0$, respect to real axis,  in  character case  $q =13$ (see tab. \ref{tab2}), and   ,  with $p_{max}=(p^*)=5.8 \times 10^6$ ( here are used the first 400038 primes).%
}}
\label {EvenQ13}
\end{center}
\end{figure}

\begin{figure}[]
\begin{center}
\includegraphics[width=1.0\textwidth]{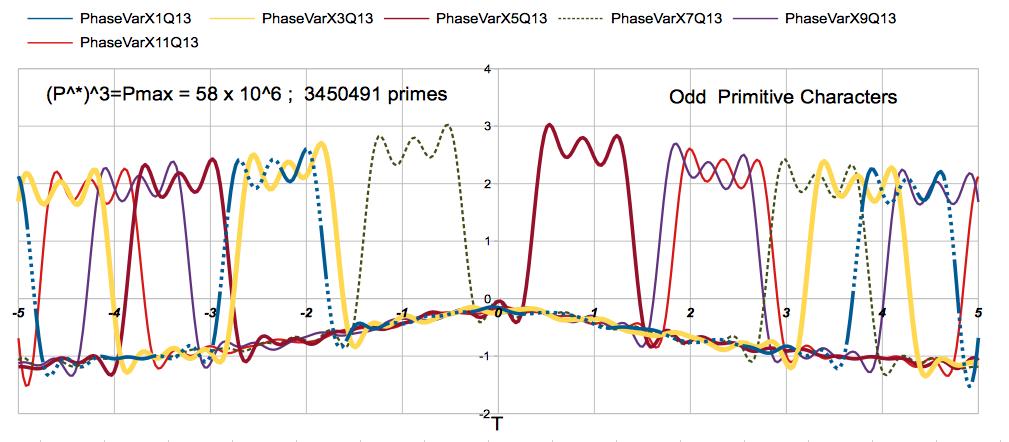} 
\caption{\small {\it   The peaks  symmetries of (\ref  {DerLAngleSuTL}                                                                                                                                                                                                                                                                                                                                                                                                                                                                                                                                                                                                              ) with  $\epsilon=0$, respect to real axis,  in the even  character $\chi_{10}$ case for   $q =13$ and   ,  with $p_{max}=(p^*)^{3}=58 \times 10^6$ ( here are used the first 3450491 primes).%
Notice that inter-peaks level is  the same of  fig, (\ref{Q3}) flipped vertically. So Lemma 2 is verified.%
}}
\label {3COddQ13}
\end{center}
\end{figure}


\begin{figure}[]
\begin{center}
\includegraphics[width=1.0\textwidth]{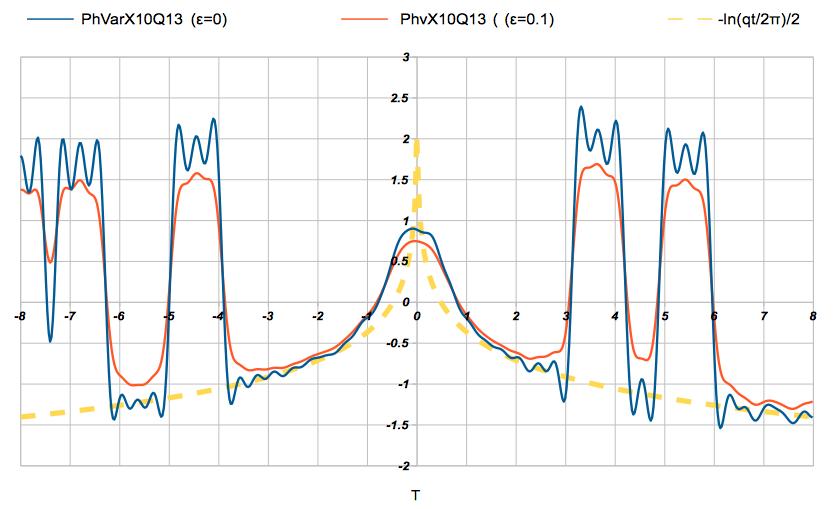} 
\caption{\small {\it   The expression (\ref  {DerLAngleSuTL}                                                                                                                                                                                                                                                                                                                                                                                                                                                                                                                                                                                                              )   for  even character  $\chi_{10}$, in   case  $q =13$ ( tab. \ref {tab2}). It  has been computed  with  $\epsilon=0$ and  $\epsilon=0.1$. Notice that  with greater $\epsilon$  the curve is close to $t$ axis as predicted in Lemma 6 . See (\ref  {LiDeltaXEps_2_GRH} ).  The choices are: $p_{max}=(p^*)^{3}=58 \times 10^6$ ( here are used the first 3450491 primes). 
}}
\label {Q13X10Eps}
\end{center}
\end{figure}

\begin{figure}[]
\begin{center}
\includegraphics[width=1.0\textwidth]{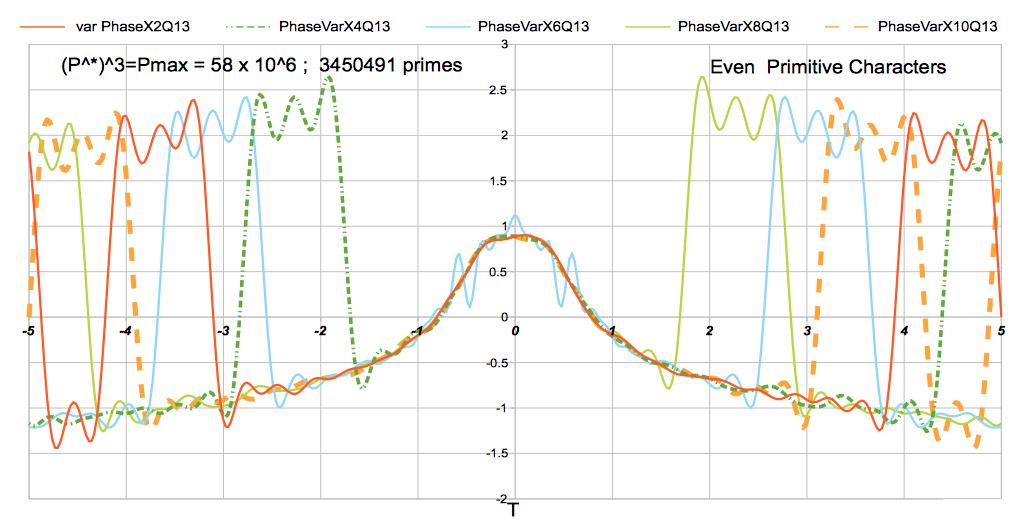} 
\caption{\small {\it   The peaks  symmetries of (\ref  {DerLAngleSuTL}                                                                                                                                                                                                                                                                                                                                                                                                                                                                                                                                                                                                              ) with  $\epsilon=0$, respect to real axis,  in  character case  $q =13$ and   ,  with $p_{max}=(p^*)^{3}=58 \times 10^6$ ( here are used the first 3450491 primes). Notice that inter-peaks level is  the same of  fig, (\ref{Q5}) flipped vertically. So Lemma 2 is verified.%
}}
\label {3CEvenQ13}
\end{center}
\end{figure}


\begin{figure}[]
\begin{center}
\includegraphics[width=1.0\textwidth]{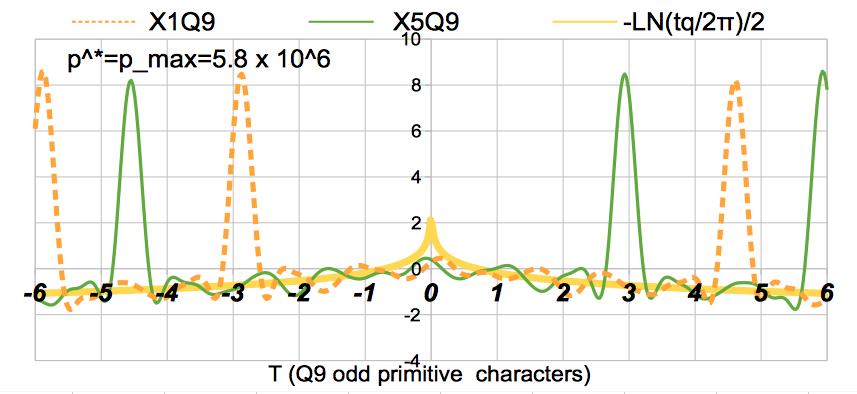} 
\caption{\small {\it   The peaks  symmetries of (\ref  {DerLAngleSuTL}                                                                                                                                                                                                                                                                                                                                                                                                                                                                                                                                                                                                              ) with  $\epsilon=0$, respect to real axis,  in  character case  $q =9$ and   ,  with $p_{max}=p^*=5.8 \times 10^6$ ( here are used the first 400038 primes).%
Notice that inter-peaks level tend to be  the same of  fig, (\ref{Q3}) with a shift  by the quantity $\frac{\ln(9/\pi)-\ln(3/\pi)}{2}$  and afterward      flipped vertically. So Lemma 2 is verified.%
}}
\label {Q9Odd5M8}
\end{center}
\end{figure}

\begin{figure}[]
\begin{center}
\includegraphics[width=1.0\textwidth]{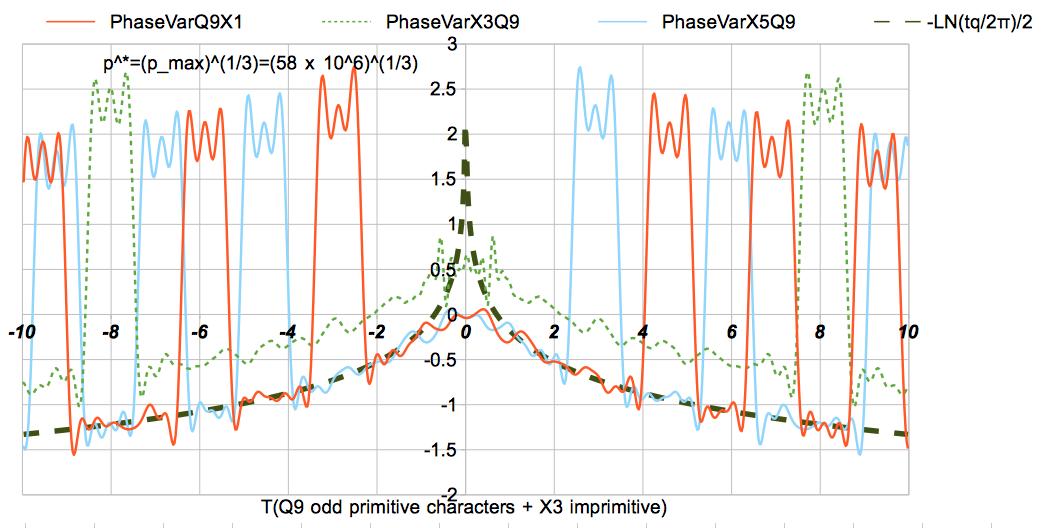} 
\caption{\small {\it   The peaks  symmetries of (\ref  {DerLAngleSuTL}                                                                                                                                                                                                                                                                                                                                                                                                                                                                                                                                                                                                              ) with  $\epsilon=0$, respect to real axis,  in  character case  $q =9$ and   ,  with $p_{max}=(p^*)^3=58 \times 10^6$ ( here are used the first 3450491 primes).%
Notice that inter-peaks level is  the same of  fig, (\ref{Q3}) with a shift  by the quantity $\frac{\ln(9/\pi)-\ln(3/\pi)}{2}$  and afterward      flipped vertically. So Lemma 2 is verified.%
}}
\label {Q9Odd58M}
\end{center}
\end{figure}

\begin{figure}[]
\begin{center}
\includegraphics[width=1.0\textwidth]{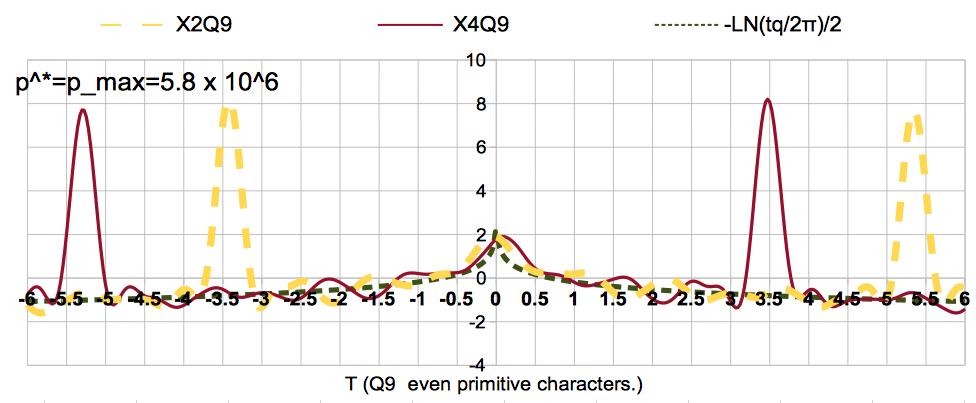} 
\caption{\small {\it   The peaks  symmetries of (\ref  {DerLAngleSuTL}                                                                                                                                                                                                                                                                                                                                                                                                                                                                                                                                                                                                              ) with  $\epsilon=0$, respect to real axis,  in  character case  $q =9$ and   ,  with $p_{max}=p^*=5.8 \times 10^6$ ( here are used the first 400038 primes).%
Notice that inter-peaks level tend to be  the same of  fig, (\ref{Q5}) with a shift  by the quantity $\frac{\ln(9/\pi)-\ln(5/\pi)}{2}$  and afterward      flipped vertically. So Lemma 2 is verified.%
}}
\label {Q9Even5M8}
\end{center}
\end{figure}

\begin{figure}[]
\begin{center}
\includegraphics[width=1.0\textwidth]{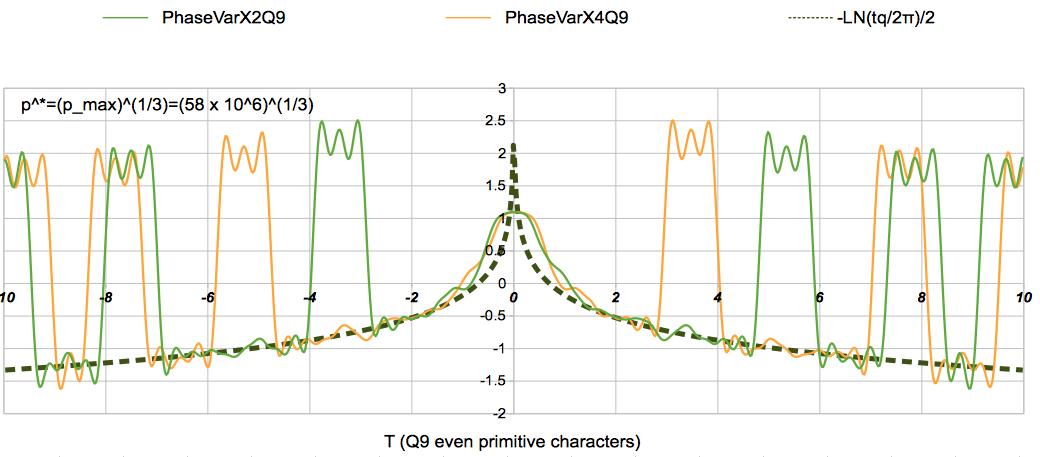} 
\caption{\small {\it   The peaks  symmetries of (\ref  {DerLAngleSuTL}                                                                                                                                                                                                                                                                                                                                                                                                                                                                                                                                                                                                              ) with  $\epsilon=0$, respect to real axis,  in  character case  $q =9$ and   ,  with $p_{max}=(p^*)^3=58 \times 10^6$ ( here are used the first 3450491 primes).%
Notice that inter-peaks level is  the same of  fig, (\ref{Q5}) with a shift  by the quantity $\frac{\ln(9/\pi)-\ln(5/\pi)}{2}$  and afterward      flipped vertically. So Lemma 2 is verified.%
}}
\label {Q9Even58M}
\end{center}
\end{figure}

\subsection {Even characters ($\alpha=0$) %
} \label {GRHEven}

For $q \ge 220 ;  \alpha =0$ we have that (\ref{DerAngleDeT}) is positive $\forall t$, see fig \ref {Q5} with  $q'=220$ , and, following ordinate increase of $\frac{\ln(q'/\pi)-\ln(5/\pi)}{2} >1.89$. So  for  $q \ge 220 ;  \alpha =0$  the argument of  Theorem 1 %
 is applicable but only for $|t|<t_{cross} \approx 0.58$, see fig.  \ref{DerEpsDerTDiGam},  because $\ \  \left\{ \frac{\partial}{\partial \epsilon}
\frac{\partial}{\partial t} \angle \left[ \Gamma \left(\frac{s + \alpha}{2 }  \right) \right] \right\}_{\alpha=0, }>0$ only for  $|t|<t_{cross} \approx 0.58$.

\noindent %
In other words, considering also tab. \ref  {tabEVEN},
 if  $L(s,\chi_{  primitive})$  does not comply with RH,  in $| t | < t_{cross} \approx 0.585$, then,  this must occurs necessarily within following constraints.

 \be \label {Xarea}
 |t|  \le 0.58 \ \ , \ \ and  \ \ , \ \
 21  \le
  q  \le 220 \ \  ; \ \  \chi_{even \ primitive}
  \ee
 
 \noindent For $q<21$  no zeros are present in $|t|  \le 0.58 $. For  $ q  > 220$ Theorm 1 applies in $|t|  \le 0.58 $.

\begin{table}[bht]
\begin{center}
\begin{tabular}{c|c|c|c|c|c|}
 q     &    &  closest even chrs correlation peak distance from $t=0$       \\
\hline
 $5$&   & $|t|>6$ \\
\hline
$7$ &   & $|t|>4$ \\
\hline
$8$ & & $|t|>4$\\
\hline
$9$ & & $|t|>3$ \\
\hline
$11$ & & $|t|>2$ \\
\hline
$12$ & & $|t|>3$ \\
\hline
$13$ & & $|t|>3$ \\
\hline
$15$ & & $|t|>2$ \\
\hline
$16$ & & $|t|>2$ \\
\hline
$17$ & & $|t|>1.5$ \\
\hline
$19$ & & $|t|>1.4$ \\
\hline
$20$ & & $|t|>2$ \\
\hline
\end{tabular}
\end{center}
\caption{\small  For even characters,  see fig \ref  {Q5}, no correlation peaks, i.e, zeros are present in   $|t|<t_{cross} \approx 0.58$ at least till $q=20$.
For $q=14$ and $q=18$  no primitive characters are present.
See also \cite {Mathar}.
}
\label {tabEVEN}
\end{table}

\begin{figure}[]
\begin{center}
\includegraphics[width=1.0\textwidth]{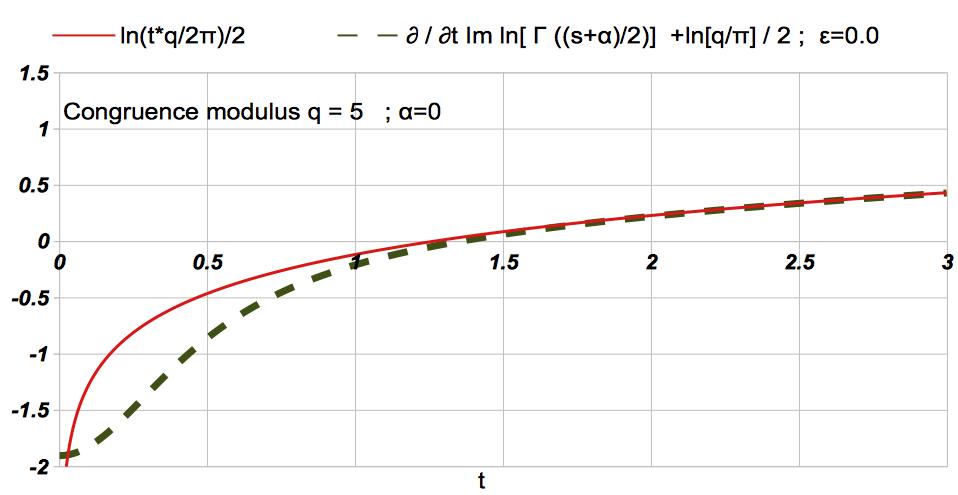} 
\caption{\small {\it   
 Plot of (\ref{DerAngleDeT}). The conditions  of Theorem 1: the curve (\ref {DerAngleDeT})  positive, like his $\frac{\partial}{\partial \epsilon}$ derivative, are  verified only for  %
$q \ge 220$  and $| t | < t_{cross }$,  see fig. \ref {DerEpsDerTDiGam} . %
Here  $T_{Asymp}(\alpha) \approx 3$.
}}
\label {Q5}
\end{center}
\end{figure}

\section { Conclusions}

\noindent For Theorem 1 in section \ref {CONJ}}, and with extentions in section \ref {TowardGRH}  we can affirm that RH  
for $\xi( (t,\epsilon,\chi_{primitive})$  
 is  true 
at least 
for all $L(s,\chi_{odd \ primitive \  } ) $ 
$ \ \ \forall q \ , \ \forall |t| %
$. %
We can say the same of $L(s,\chi_{even \ primitive \  } ) $ but only within the constraints $$|t|  \le 0.58 \ \ , \ \ and  \ \ , \ \
  q  >  220 \ \  ; \ \  \chi_{even \ primitive}$$
Application of main idea to $\zeta(s)$ is proposed in a  pre-print  paper \cite {Giovanni Lodone Dec2024}.

%

%

$$Acknowledgments
$$
I thank Paolo Lodone for useful discussions, and,  contributions.


\vspace{0.3cm}

 %

%

\appendix

\section {Phase variations of $\xi(s)$,  $\xi(s,\chi)$, 
$\zeta(s)$ and $L(s,\chi)$ } \label {GammaDer}

Phase variations along t  for $\xi(s)$ and $\zeta(s)$ (see (\ref {XiForLfunc} ) )are connected by:

\be \label {PhaseVarXiAndZeta}
\left[ \frac {\partial \angle[\xi(1/2+\epsilon+i t  )]  } {\partial t} \right]
=
    \frac {\partial \angle[\zeta(s)(s-1)]}{\partial t} +\ln\left[ \sqrt{\frac{t }{2 \pi}} \right]
    \quad + O(t^{-2}) %
      \ee
Phase variations along t  for $\xi(s,\chi_{primitive})$ and $L(s,\chi_{primitive})$ (see (\ref {XiForLfunc} )) are connected by:

\be \label {PhaseVarXiAndL}
\left[ \frac {\partial \angle[\xi(s,\chi_{primitive}) )]  } {\partial t} \right]%
=
    \frac {\partial \angle[L(s,\chi_{primitive})]}{\partial t} +\ln\left[ \sqrt{\frac{t q }{2 \pi}} \right]
    \quad   + O(t^{-2})%
    \ee
     ( $q$ is the congruence modulus). %
To prove eq. (\ref {PhaseVarXiAndZeta}),  and (\ref  {PhaseVarXiAndL}) notice that:
\be \label {FattorePerXIA}
\angle \left[   \Gamma\left(\frac{s}{2}  \right)\frac {s(s-1)} {2}\pi^{-\frac{s}{2}} \zeta(s) \right]=
  \angle \left[  \zeta(s)(s-1)  \   \Gamma\left(\frac{s}{2} +1 \right)   \pi^{-\frac{s}{2}} \right] =  \angle \left[  \zeta(s)(s-1) \right ]+\angle \left [   \Gamma\left(\frac{s}{2} +1 \right)   \pi^{-\frac{s}{2}} \right] 
\ee

While for eq. (\ref {PhaseVarXiAndL}), and,  (\ref{XiForLfunc} ) we have:
\be \label {FattorePerXIAChi}
 \theta_{primitive}(t) =   \angle \left[  \left( \frac{\pi}{q}   \right)^{-\frac{s+\alpha_1}{2}} \Gamma  \left( \frac{s+\alpha}{2} \right)  \ \  L(s,\chi) \right]= 
  \angle \left[  \left( \frac{\pi}{q}   \right)^{-\frac{s+\alpha_1}{2}} \Gamma  \left( \frac{s+\alpha}{2} \right)  \right] + \angle \left [  L(s,\chi) \right]
\ee
With $\alpha_1=\alpha=0  $, even $\chi_{primitive}$, or  $\alpha_1=\alpha=1  $, odd $\chi_{primitive}$. See (\ref  {XiForLfunc} ).
\noindent If we consider the factor $(s-1)$ attached to $\zeta(s)$ like in (\ref {PhaseVarXiAndZeta}),
 \noindent we can use  (\ref {FattorePerXIAChi}) again putting :%
  $\alpha_1=0 \ , \ \alpha=2 \ , \ q =1$. As $z+1=  \frac{s+\alpha}{2} \rightarrow z=\frac{2\epsilon+2\alpha-3}{4}+i \ \frac{t}{2}$,  Stirling formula ( \cite[p.~109,112]{Edwards:1974cz}) allows us to write  (\ref {FattorePerXIAChi}) for  each of the 3 choices of $(\alpha,\alpha_1)$:

\be\label {Pi&Gamma3}
\ln( \Gamma(z+1) )  +\ln\left[\left( \frac{\pi}{q} \right)^{ - \frac{s+\alpha_1}{2}} \right]=
\ln( \Gamma(z+1) )  -\frac{s+\alpha_1}{2} \ \ln\left(\frac{\pi}{q}\right) =
\ee
$$\ln\left( e^{-z}z^{   z+\frac{1}{2}   }(2\pi)^{\frac{1}{2}}  \right )  +  \left( \sum_{k=1}^{K-1} \frac{B_{2k}   }{  2 k (2k-1) z^{2k-1}} \right) +R_{2K}(z)  -\ln\left(\frac{\pi}{q}\right)\frac{s+\alpha_1}{2} =\Re_1+\Re_2+\Re_3 + i(\Im_1+\Im_2+\Im_3)
$$
Although the expression (\ref {Pi&Gamma3})  is a not-convergent asymptotic expansion, it  can be used  with a finite $K$  keeping the modulus of error term    $|R_{2K(z)}|$ to a suitable level.  Let us call    $\Im_1$    the imaginary  of  asymptotic  (for $t $ big) part of   $\ln\left( e^{-z}z^{   z+\frac{1}{2}   }(2\pi)^{\frac{1}{2}}  \right ) $. And $\Re_1$ is the real part.  Let us call    $\Im_2$   the left imaginary part  of  $\ln\left( e^{-z}z^{   z+\frac{1}{2}   }(2\pi)^{\frac{1}{2}}  \right ) $ wich goes to zero  for $t\rightarrow \infty$. And $\Re_2$ is the real part.And finally let us call $\Re_3+i \Im_3 =    \left( \sum_{k=1}^{K-1} \frac{B_{2k}   }{  2 k (2k-1) z^{2k-1}} \right) +R_{2K}(z)   $.

\noindent The $B_j$  are the Bernoulli numbers  ( \cite[p.~11] {Edwards:1974cz} ), that vanish for odd $j$ while: 
$$
B_{2} =\frac{1}{6}
\quad ; \quad
 B_{4} =-\frac{1}{30}
\quad ; \quad
B_{6} =\frac{1}{42}
\quad ; \quad
B_{8} =-\frac{1}{30}
\quad ; \quad
B_{10} =\frac{5}{66} 
\quad ; \quad
 B_{12} =-\frac{691}{2730}
\quad ...
$$

The modulus of  the error term $|R_{2K}(z)|$ is bounded by: 
\be \label  {RestoInApp}
 |R_{2K}(z)| < \left(    \frac{B_{2K}   }{  2 K (2k-1) z^{2K-1}}          \right)
\frac{1}{       \left[     \cos\left ( \frac{\arg(z)}{2}\right)         \right]^{2K}        }
\ee
where $\arg(z)$ is taken in the interval: $-\pi <\arg(z) <\pi$ (see  \cite[p.~112]{Edwards:1974cz}.

 \noindent  In (\ref{Pi&Gamma3}) we may put : $z=\frac{2\epsilon+2\alpha-3}{4}+i \ \frac{t}{2}$ , so,   leaving apart only $\Im_3$ for now, we have:

 $$
\ln\left[   \left( \frac{\pi}{q}   \right)^{-\frac{s+\alpha_1}{2}} \Gamma  \left( \frac{s+\alpha}{2} \right) \right] 
\approx  \Re_1 + \Re_2+ i ( \Im_1 + \Im_2) = %
 $$
 $$
 -\left( \frac{s+\alpha}{2} -1 \right)+\left(  \frac{s+\alpha}{2} -\frac{1}{2} \right)\ln\left( \frac{s+\alpha}{2} -1  \right)+\ln(\sqrt{2 \pi} )  -\ln\left(\frac{\pi}{q}\right)\frac{s+\alpha_1}{2}  =
 $$
 
 $$
 \frac{3}{4} - \frac {\epsilon+ \alpha}{2}-\frac{i t}{2} +\left (\frac {\epsilon+ \alpha}{2}- \frac{1}{4} + \frac{i t}{2}\right) \left\{     \ln  \sqrt{\left( \frac {\epsilon+ \alpha}{2}-\frac{3}{4} \right)^2+\frac{t^2}{4}}  +i \arctan \left[ \Im=\frac{ t}{2} ; \Re = \frac {\epsilon+ \alpha}{2}-\frac{3}{4} \right]  \right\} +
 $$
 \be \label {NoApproxSalvo Bernoulli}
 \frac{\ln(2 \pi)}{2}-\ln\left(\frac{\pi}{q} \right)\left(    \frac {\epsilon+ \alpha_1}{2} +\frac{1}{4}   \right)  +i\frac{t}{2}\ln\left(  \frac{q}{\pi}   \right)
 \ee

\noindent   Taking into account that:
 $
  \arctan \left[ \Im=\frac{ t}{2} ; \Re = \frac {\epsilon+ \alpha}{2}-\frac{3}{4} \right]=\frac{\pi}{2}+\arctan \left[  \frac{3-2(\epsilon+\alpha)}{2t}  \right] 
 $, we get:
  \be \label {NoApproxSalvo Bernoulli2}
 \frac{3}{4} - \frac {\epsilon+ \alpha}{2}-\frac{i t}{2} +\left (\frac {\epsilon+ \alpha}{2}- \frac{1}{4} + \frac{i t}{2}\right) \left\{     \ln  \sqrt{\left( \frac {\epsilon+ \alpha}{2}-\frac{3}{4} \right)^2+\frac{t^2}{4}}  +i \left \{ \frac{\pi}{2}+\arctan \left[  \frac{3-2(\epsilon+\alpha)} {2t}  \right]   \right \}   %
  \right\}
  +
  \ee

$$
 \frac{\ln(2 \pi)}{2}-\ln\left(\frac{\pi}{q} \right)\left(    \frac {\epsilon+ \alpha_1}{2} +\frac{1}{4}   \right)  +i\frac{t}{2}\ln\left(  \frac{q}{\pi}   \right)
$$

\noindent Let us choose $T_{Asymp}(\alpha)$, in order  for  $ t>T_{Asymp}(\alpha) $ the braces in  (\ref{NoApproxSalvo Bernoulli2})  becomes $ \left\{  \ln \left( \frac{t}{2}\right) +i \frac{\pi}{2}
   \right\}$,so   we have:
 $$
\ln\left[   \left( \frac{\pi}{q}   \right)^{-\frac{s+\alpha_1}{2}} \Gamma  \left( \frac{s+\alpha}{2} \right) \right] \approx \Re_1+i \Im_1=
$$
 $$
 \frac{3}{4} - \frac {\epsilon+ \alpha}{2}-\frac{i t}{2} +\left (\frac {\epsilon+ \alpha}{2}- \frac{1}{4} + \frac{i t}{2}\right)
  \left\{  \ln \left( \frac{t}{2}\right) +i \frac{\pi}{2}
   \right\} +\frac{\ln(2 \pi)}{2}-\ln\left(\frac{\pi}{q} \right)  \left(    \frac {\epsilon+ \alpha_1}{2} +  \frac{1}{4}   \right) +i\frac{t}{2}\ln\left(  \frac{q}{\pi}   \right)
 $$

\noindent   From which we find the  the asymptotic  imaginary parts  $\Im_1 $ :

 \be \label {ImFattDiXiDiSeChi}
\Im_1 =
 \theta_{primitive}(t)
 = 
  -\frac {t}{2}+ \frac {t}{2}\ln\left( \frac{t q}{2 \pi}  \right)- \frac{\pi}{8}+\frac{\pi}{4}(\epsilon+\alpha)
 \ee
 From (\ref {NoApproxSalvo Bernoulli2}) the imaginary part $\Im_2$ that goes to zero as $t\rightarrow \infty$:

 \be \label {IMPartToZero}
 \Im_2 :=  \frac{t}{2}  \ln\sqrt{1+ \left( \frac{2(\epsilon+ \alpha)-3}{2t}  \right)^2} + \left(    \frac {\epsilon+ \alpha}{2} -  \frac{1}{4}   \right) \arctan \left(  \frac{3-2(\epsilon+\alpha)}{2t}      \right)
   \ee

 for $t$ big, posing $y=\epsilon+\alpha$: 
  \be \label{I2Approx} 
   \Im_2  \approx \frac{t}{4}  \left( \frac{2(\epsilon+ \alpha)-3}{2t}  \right)^2 +
   \left(    \frac {\epsilon+ \alpha}{2} -  \frac{1}{4}   \right)  \left(  \frac{3-2(\epsilon+\alpha)}{2t}      \right)
   =
   \frac{t}{4}  \left( \frac{2(y-3)}{2t}  \right)^2 +
   \left(    \frac {y}{2} -  \frac{1}{4}   \right)  \left(  \frac{3-2y}{2t}      \right)
  \ee



 $$ \frac {\partial^2 \Im_2}{ \partial y \partial t} =-\frac{4(2y-3)}{16t^2}-  \frac{1}{2} \frac{3-2y}{2t^2}
+ \left( \frac{y}{2}-\frac{1}{4} \right)\frac{2}{2t^2}
 $$

 For $\epsilon=0$ and $\alpha=2$ we have:
 $$ \left[  \frac {\partial^2 \Im_2}{ \partial \epsilon \partial t}  \right]_{\epsilon=0 \ , \ \alpha=2 , \ \alpha_1=0} =
 - \frac{4}{16t^2}+  \frac{1}{4t^2} + \frac{3}{4} \frac {2}{2 t^2}=\frac{3}{4t^2}$$

For $\epsilon=0$ and  $\alpha=1$ we have:
 $$\left[  \frac {\partial^2 \Im_2}{ \partial \epsilon \partial t}  \right]_{\epsilon=0 \ , \ \alpha_1= \alpha=1} =
 + \frac{4}{16t^2} -  \frac{1}{4t^2} + \frac{1}{4} \frac {2}{2 t^2}=\frac{1}{4t^2}$$
 
 For $\epsilon=0$  and $\alpha=0$ we have:
 $$\left[  \frac {\partial^2 \Im_2}{ \partial \epsilon \partial t}  \right]_{\epsilon=0 \ , \   \alpha_1=\alpha=0} =
 + \frac{12}{16t^2} -  \frac{3}{4t^2} - \frac{1}{4} \frac {2}{2 t^2}= - \frac{1}{4t^2}$$

  

 While $\Im_3=  \Im \left[     \left( \sum_{k=1}^{K-1} \frac{B_{2k}   }{  2 k (2k-1) z^{2k-1}} \right) +R_{2K}(z)  \right]$ for $K=3$.So $k = 1 , and \ 2$

\be \label {Im3TillK2}
\Im_3 -R_{2K}(z)= \frac {-1}{6 t\left[ 1+\left( \frac{2 \epsilon+2 \alpha-3 }{2t}  \right) ^2\right]}
-\frac {1}{45 t^3 \left[1 +\left(\frac{2\epsilon+2\alpha-3}{2t}  \right)^2 \right]^3}
+
\frac {(2\epsilon+2\alpha-3)^2}{60 t^5 \left[1 +\left(\frac{2\epsilon+2\alpha-3}{2t}  \right)^2 \right]^3}     %
\ee
 
  \noindent $\Im_3$  can be neglected because (we take only  the first term)$\Im_3 \approx \frac{-1}{6t} \left( 
 1-\left( \frac{2 \epsilon+2 \alpha-3 }{2t}  \right) ^2 \right) $.
 
 \noindent Besides : $ \frac {\partial \Im_3}{\partial t} \approx \frac{1}{6t^2}  
 +3  \frac{(2 \epsilon+2 \alpha-3)^2 }{24t^4} $. While   $ \frac {\partial^2 \Im_3}{ \partial \epsilon \partial t} \approx  \frac{(2 \epsilon+2 \alpha-3) }{4t^4} $, so, the change with $\epsilon$ is negligible.
\noindent  The  derivative  of (\ref{ImFattDiXiDiSeChi}) with respect to $t$ is:
 
  \be \label {DerAngFattXiChiPrim}
   \frac {\partial \angle\left[ \left( \frac{\pi}{q}   \right)^{-\frac{s+\alpha_1}{2}} \Gamma  \left( \frac{s+\alpha}{2} \right)    \right]}{\partial t} = \frac{\partial \Im_1}{\partial t}+ \frac{\partial \Im_2}{\partial t}+ \frac{\partial \Im_3}{\partial t}
   \ee
   
   where: 

  \be \label {Asympt} 
   \frac{\partial \Im_1}{\partial t}
   = \ln\left( \sqrt{ \frac {t q}{2 \pi} }   \right) \quad ; \quad q=1 \  for \ \zeta(s) \ ; \ q>1 \ for  \  L(s,\chi_{primitive}) \ ; \ t> T_{Asymp}(\alpha)
   \ee
   
  \noindent Notice that  $T_{Asymp}(\alpha)$,  allows us to use the asymptotic relation (\ref{Asympt} ).  Following (\ref {FattorePerXIAChi})  it cannot depend on $q$  or on $\alpha_1$, but only on  $\Gamma  \left( \frac{s+\alpha}{2} \right) $.

\noindent $%
\frac {\partial^2 \Im_3}{ \partial \epsilon \partial t} 
 $ is negligible, and $\frac {\partial}{\partial \epsilon} 
  \frac{\partial \Im_1}{\partial t}%
 =0 $.  So  doing  $\frac{\partial}{\partial \epsilon} $ of (\ref {DerAngFattXiChiPrim}) we have 
the only three cases of $\frac {\partial^2 \Im_2}{ \partial \epsilon \partial t} $ ( from (\ref {DerAngFattXiChiPrim} )  : %

    \be \label {3Cases}
\overbrace{ %
\left[
\frac{\partial}{\partial \epsilon} 
   \frac{\partial\angle\left[ \Gamma\left( \frac{s+\alpha}{2}   \right) \right]}{\partial t}   \right]_{\epsilon=0}
    =  
    \frac {  3 }{4 t^2}  } ^{\zeta(s) \  : \alpha=2}\ \ ; \ \overbrace{
 \left[
\frac{\partial}{\partial \epsilon} 
   \frac{\partial\angle\left[ \Gamma\left( \frac{s+\alpha}{2}   \right) \right]}{\partial t}   \right]_{\epsilon=0}
    =  \frac {  1      }{4 t^2}  }^{ \chi_{primitive}(-1) =-1 \ : \alpha=1} \ \ ; \ \ 
\overbrace{ 
 \left[
\frac{\partial}{\partial \epsilon} 
   \frac{\partial\angle\left[ \Gamma\left( \frac{s+\alpha}{2}   \right) \right]}{\partial t}   \right]_{\epsilon=0} 
   =  - \frac {   1      }{4 t^2} }^{   \chi_{primitive}(-1) =1 \ :  \alpha=0  }
\ee

\begin{figure}[]
\begin{center}
\includegraphics[width=1.0\textwidth]{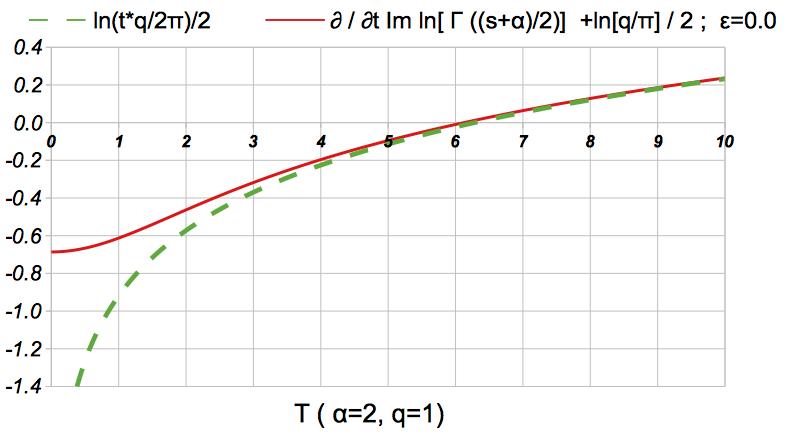} 
\caption{\small {\it   
 Plot of (\ref{DerAngleDeT}) for $\alpha=2$ and $q=1$, i.e. for the $ \frac {\partial \angle[\zeta(s)(s-1)]}{\partial t}$ function. See (\ref{PhaseVarXiAndZeta}).  Notice that the asymptotic phase variation with $t$ (\ref{Asympt} ) is the same surely for $t>10$, i.e.  $T_{Asymp}(\alpha) \approx 10$ here. %
}}
\label {ForZfunc}
\end{center}
\end{figure}

\section {differences between  (\ref {FdipparamChiH})  and (\ref {FdaLgaritmicInt})} \label {ApproxExpr}

The difference between  the infinite sum (\ref {DerLAngleSuTL})  and  the same infinite sum with substitution  of  (\ref  {FdipparamChiH}) with (\ref {FdaLgaritmicInt})  is always a finite number at least for $\epsilon \ge 0$. So the correlation peaks  of  (\ref {DerLAngleSuTL}) stemming from the zeros of $L(s,\chi)$,  see fig.  \ref{Simmetries}, %
are present also in  (\ref {DerLAngleSuTL})  after subsitution.

\noindent PROOF. Using Taylor formula for arctangent the difference between (\ref {DerLAngleSuTL})  using (\ref  {FdipparamChiH}) or (\ref {FdaLgaritmicInt})   can be written as:

$$
-  \left. \frac {\ln(p^*)} {2 \pi}   \sum_{ p<p_{max}   \ : \ gcd(p,q)=1    } \arctan \left( \frac{
 \sin(\ln(p) t  - \angle[\chi(p)]    )}{ p^{1/2+\epsilon}-  
\cos( \ln(p) t  - \angle[\chi(p)] )  }\right)\right|_{t_1}^{t_2}
+
$$
$$
 \frac {\ln(p^*)}{2 \pi} 
 \sum_{  p<p_{max}   \ : \ gcd(p,q)=1   } \frac{2 \cos [\ln(p) t - \angle[\chi(p)] ] \sin\left( \frac{\ln(p) \pi}{\ln(p^*)}   \right) }{p^{1/2+\epsilon}} =
$$

\be \label {SumInPandN}
 -  \left. \frac {\ln(p^*)} {2 \pi}  \ \ \  \sum_{ p<p_{max}   \ : \ gcd(p,q)=1    } \  \ \ \ \  \sum_{n  \ odd >1} \frac{(-1)^{(n-1)/2}}{n}
  \left[   \left( \frac{
 \sin(\ln(p) t    - \angle[\chi(p)]  )}{ p^{1/2+\epsilon  } 
- \cos(\ln(p) t   - \angle[\chi(p)]   )  }\right)^n   \right] \right |_{t_1}^{t_2}
- 
\ee

\be \label {SumPminus2}
\left. \frac {\ln(p^*)} {2 \pi}   \sum_{ p<p_{max}   \ : \ gcd(p,q)=1    }  \left( \frac{
 \sin(\ln(p) t   - \angle[\chi(p)]   )    \cos(\ln(p) t   - \angle[\chi(p)]   )   }{ [p^{1/2+\epsilon}-  
\cos( \ln(p) t - \angle[\chi(p)]  ) ]    p^{1/2+\epsilon}   }\right)\right|_{t_1}^{t_2}
\ee

 \noindent Developing last sum we have, focusing on   main term:

\be \label {SumPminus3}
\left .\frac {\ln(p^*)} {4 \pi}   \sum_{ p<p_{max}   \ : \ gcd(p,q)=1    }  \left( \frac{
 \sin( 2( \ln(p) t   - \angle[\chi(p)]  ) )     }{ [p^{1/2+\epsilon}-  
\cos( \ln(p) t - \angle[\chi(p)]  ) ]    p^{1/2+\epsilon}   }\right)\right|_{t_1}^{t_2} \approx
\ee

\be \label {SumPminus4}
 \frac {\ln(p^*)} {4 \pi}   
\sum_{ p<p_{max}   \ : \ gcd(p,q)=1    } 
 \left( \frac{
 \sin( 2( \ln(p) t_2  - \angle[\chi(p)]  ) )  -
 \sin( 2( \ln(p) t_1   - \angle[\chi(p)]  ) ) 
    }
 { [p^{1/2+\epsilon}-  
\cos( \ln(p) t_2 - \angle[\chi(p)]  ) ]    [p^{1/2+\epsilon}-  
\cos( \ln(p) t_1 - \angle[\chi(p)]  ) ]     }\right)  
\ee

\noindent It is apparent that the difference between computation (\ref {DerLAngleSuTL}                                                                                                                                                                                                                                                                                                                                                                                                                                                                                                                                                                                                             ),  upon substitution of (\ref {FdipparamChiH})  with (\ref {FdaLgaritmicInt})  is absolutely convergent for  sure from $\epsilon>0$  and  conditionally convergent for 
$\epsilon =0 $.

\noindent Proof.         Let us look at (\ref {HalfRotation_GRH}), and, consider only one $h-$class primes till $p_{max}$ ( as in (\ref{DerLAngleSuTL2}   ). 

\noindent Each equal sign interval (see (\ref {HalfRotation_GRH}))  now $\Delta p= \frac {\pi}{2t} p$ alternating sign, in  (\ref{SumPminus4}) is divided by $\approx p$.  But the prime density is $\frac{1}{\ln(p)}$. So we have an alternating sign decreasing module series. The limit must be finite by Leibnitz rule\cite [p.~404] {Apostol:1967}.  So the divergences at $L(s,\chi)$ zeros, i.e. correlation peaks,  
   are the same.

%

 \noindent The double sum, in $p$ and in $n  \ odd >1$, as the neglected term between  (\ref {SumPminus3}) and  (\ref {SumPminus4})   cannot influence divergences all the same, because of the common  factor $ \sim p^{- 3(1/2+\epsilon)}$.%
 
 \noindent If we consider all the $\phi(q)$ terms ( possible if $p_{max} < \infty$), even if we cannot  conclude that the sum of all $\phi(q)$ contributions tend to zero ( as in Lemma 4), we can presume that the poor convergence of a single $h-$class primes  is enanced when all $h-$classes are combined together. This seems confirmed by figg. \ref{ApproxX1Q5}, \ref  {ApprQ3}, and, \ref  {ApprQ4}.

\begin{figure}[]
\begin{center}
\includegraphics[width=1.0\textwidth]{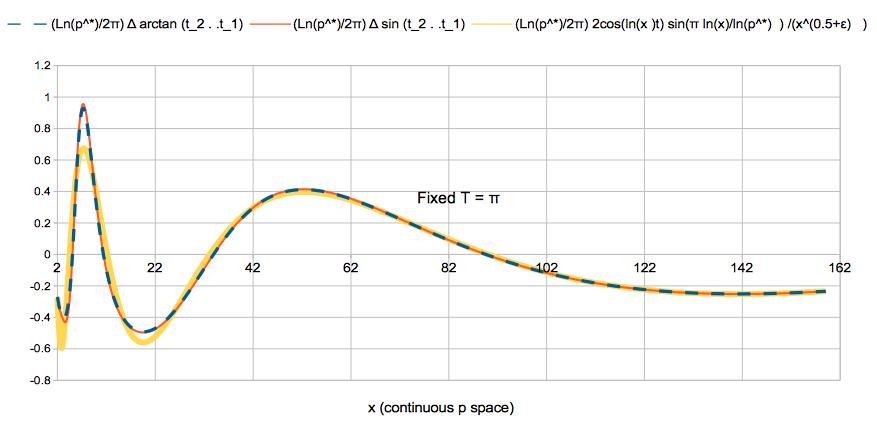} 
\caption{\small {\it   
 Comparative plot of  (\ref {DerLAngleSuTL})  using   (\ref {FdipparamChiH})  or (\ref {FdaLgaritmicInt}), with $t$ fixed, $\angle[\chi]=0$, and continuously varying $p$. It is reported also (\ref {DerLAngleSuTL}) with substitution of  intermediate expression $\left\{  \left( \frac{
 \sin(\ln(p) t    - \angle[\chi(p)]  )}{ p^{1/2+\epsilon}-  
\cos( \ln(p) t - \angle[\chi(p)]  )  }\right) \right \}_{t-\frac{\pi}{\ln(p^*)}   }^{t+\frac{\pi}{\ln(p^*)}  }$%
}}
\label {Approxs}
\end{center}
\end{figure}


\begin{figure}[]
\begin{center}
\includegraphics[width=1.0\textwidth]{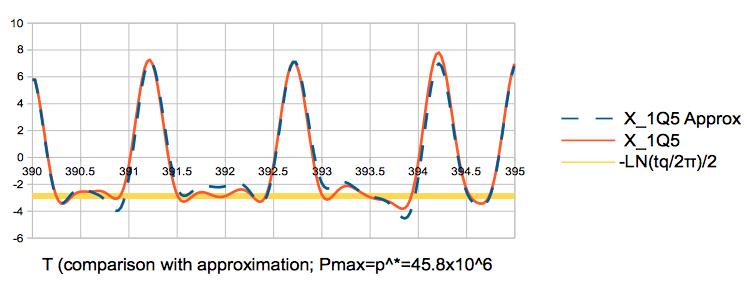} 
\caption{\small {\it   
 Comparative plot of  (\ref {DerLAngleSuTL}) computed with approximation   (\ref {FdaLgaritmicInt}), with $\epsilon=0$, and $390 < t <395 $, for  $\chi_1$ and $q=5$. See tab. \ref {tab1}. We have used the first 2.763.823  primes like in following figures.
}}
\label {ApproxX1Q5}
\end{center}
\end{figure}

\begin{figure}[]
\begin{center}
\includegraphics[width=1.0\textwidth]{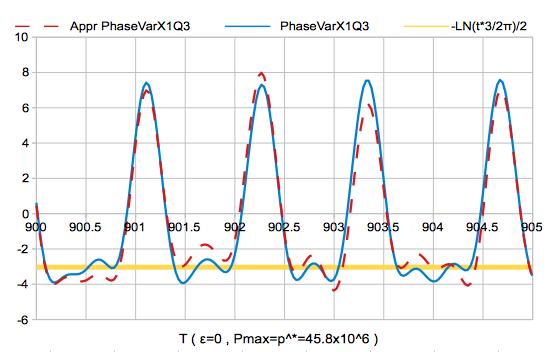} 
\caption{\small {\it   
 Comparative plot of  (\ref {DerLAngleSuTL}) computed with approximation   (\ref {FdaLgaritmicInt}), with $\epsilon=0$, and $900 < t <905 $, for  $q=3$, and,  $\chi_1(n)=(0,1,-1)$,  to comply with format of tab. \ref {tab2} . Or  $\angle[\chi_1(n)]=(undefined, 0 , \pi)$ to comply with format of tab. \ref {tab1}.
}}
\label {ApprQ3}
\end{center}
\end{figure}

\begin{figure}[]
\begin{center}
\includegraphics[width=1.0\textwidth]{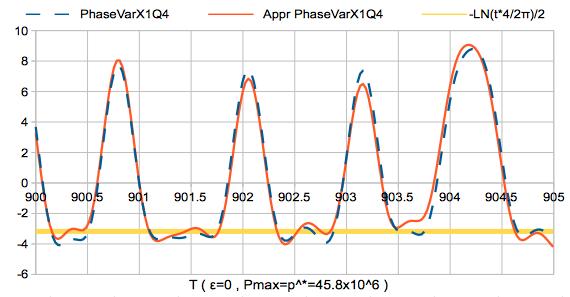} 
\caption{\small {\it   
 Comparative plot of  (\ref {DerLAngleSuTL}) computed with approximation   (\ref {FdaLgaritmicInt}), with $\epsilon=0$, and $900 < t <905 $, for $q=4$,  and $\chi_1(n)=(0,1,0,-1)$,  to comply with format of tab. \ref {tab2}. Or,  $\angle[\chi_1(n)]=(undefined, 0 , undefined, \pi)$ to comply with format of tab. \ref {tab1}.  Notice in interval $904<t<904.5$ there are two zeros not resolved by the implicit  resolution: $\frac{2 \pi}{\ln(45.8 \times10^6) } \approx 0.356 $
}}
\label {ApprQ4}
\end{center}
\end{figure}


\section {
Convergence for primitive characters in critical strip}
\label {EPinCS}


\noindent
Let us see a way to show that  
\be \label {L}
\L(s.\chi_{primitive})=\sum _{n=1}^\infty  \   \chi(n)/n^s
\ee
 converge 
in any compact on critical strip.

 \noindent Abel summation formula , in the version of \cite  {wikiAbelForm}, is reported for easy reading.
 {\it Let $\{ a_n \}_{n=0}^\infty$ be a sequence of real or complex numbers. Define partial sum function A by $A(t)=\sum_{0 \le n \le t}      a_n $  for any real number $t$. Fix real number $x<y$ and let $\phi$ a continuously differentiable function on $[x,y]$. Then $\sum_{x<n\le y}  a_n \phi(n)=A(y) \phi(y)         
 - A(x) \phi(x) - \int_x^y A(u) \phi'(u)du$
 }
 

Posing
\be \label {PartialSum1}
K(X)= \sum_{0<x \le X}\chi(x)
\ee
 we have:
 
\be \label {AbelPerL}
 \sum_{n=1}^{n=N=\lfloor X \rfloor }\chi(n) \times \left(  \frac{1}{n^s} \right)
=
K(X)\times  \left(  \frac{1}{X^s} \right) - K(1) \times \left(  \frac{1}{(1)^s} \right) -(-s) \int_1^X\frac{K(u) du}{u^{s+1}}
\ee

as for $X \rightarrow \infty$:
\begin{itemize}
\item $ K(X)$  is bounded. For sure  it is $<\phi(q)$.
\item $K(X)\times  \left(  \frac{1}{X^s} \right) \rightarrow 0$, for $\Re(s) >0$
\item $K(1)\times  \left(  \frac{1}{1^s} \right) = 1$ 
\item  $\int_1^X\frac{K(u) du}{u^{s+1}}$  converges for $\Re(s) >0$ because of first bullet.
\end{itemize}

\noindent  So  $\L(s.\chi_{primitive})=\sum _{n=1}^\infty  \   \chi(n)/n^s$ converges in $\Re(s) >0$

The  crucial point for series  (\ref {L})  is  the bounded behavior of  (\ref{PartialSum1}).

\section{Some feedback from readers}   \label {FEEDBACK}
I thank all the readers who have send to me comments on previous versions of the paper. Below is given a list of them (unified following subject) with my answers. They are useful especially for me to delve deeper in this hard subject.

\begin{itemize}
\item  {\it  Invalid Use of the Euler Product in the Critical Strip. If  it is correct that  the Dirichlet series $\L(s.\chi_{primitive})=\sum _{n=1}^\infty  \   \chi(n)/n^s$
 converges
for $\Re(s) > 0$ (for non-principal  character) and that the L-function has an Euler product. However,
that these two facts imply the equality $\sum_n =\prod_p$is incorrect for $\Re(s)>0$.
The proof of the
identity, first discovered by Euler, relies on multiplying out the infinite geometric series
for each prime factor and
rearranging the order of an infinite number of terms. This rearrangement is only guaranteed
to preserve the value of the sum if the series is absolutely convergent. So section \ref{StabAndAbsValForL} is baseless.
}. 

\noindent ANSW.I agree on the necessity to fill this gap (not a mistake I think)  in the logical path of previous versions of the paper.  But, if, for $\Re(s)>0$, 
$$\L(s.\chi)- EP(s,\chi)=\left[\sum _{n=1}^{p_{max}}  \   \chi(n)/n^s\right] - \left[ \prod_2^{p_{max}}(1-\chi p^{-s})^{-1} \right]$$
 can be seen as a absolutely convergent series for fixed $p_{max}$. And besides, if  its limit tends to zero as $p_{max}\rightarrow \infty$ as it  is shown in (\ref{DirichletLFuncExpanded}) . Then, it seems to me, that, now,  this gap has been filled.

\noindent In  other words in section \ref  {StabAndAbsValForL}
we have  shown now  that   (\ref{DirichletLFunc})   can be extended in critical strip, like in  (\ref{DirichletLFuncExpanded} ),   at least for primitive $L-$Functions. 

\noindent 

\noindent In conclusion it seems to me that  section \ref {EPInCR} 
is now a valid base for the developments of further  arguments of the paper


\noindent Anyhow I thank the reader for this observation that  prompted me to delve deeper on subject.

\item {\it Incorrect Manipulation of Conditionally Convergent Series. The argument in Lemma 6 is invalid.
The series for the phase is conditionally convergent, not absolutely convergent. In
a conditionally convergent series, the order of summation matters

} 

\noindent ANSW. I agree with the statement, but it is not applicable.The ``manipulation'' of terms  in  (\ref{FinRes} ) (or in others occasions (\ref {EulerTrickFiniteTerms}), (\ref {DerLAngleSuTL2}) and (\ref{WithLogInt}) ) is applied only  at addends with index less than the  running index. So on a set of finite numbers. Partial sums do not change (as it  is stressed thoughout). To explain better. Consider an example  $S_j=-1+1/2-1/3+1/4- . .. . .= \sum_j (-1)^j/(j) \rightarrow - \ln(2) \ ,$ as $j \rightarrow \infty$.
                                                                                                                                                                                                                      
                                                                                                                                                                                                                      \noindent If we take the sum till $J_{max}$    and then we partition the addends with $ j  \le J_{max}$  in positive $\sum^+_{J_{max}} $ and negatives $\sum^-_{J_{max}} $, then we can write $S_{J_{max}}= \sum^+_{J_{max}} -\sum^-_{J_{max}} $, and, we do not alter final result $- \ln(2)$ for  $J_{max} \rightarrow \infty$ transforming the limit in something like $\infty - \infty$. It is exactly the ``manipulation''  in expression : (\ref{FinRes} ). While in (\ref {EulerTrickFiniteTerms}), in (\ref {DerLAngleSuTL2}),  and, in (\ref{WithLogInt}),  we manipulate to highlight some structure of the series, but ,  with { \bf indexes less than the  running index }.  A  completely rearranging of terms  as in  \cite [p.~411]  {Apostol:1967}  without change in final result is possible only where convergence is absolute. But, 
                                                                                                                                                                                                                      it is easy to see  that there are  infinite rearrangement of terms that preserve the sum $- \ln(2)$ of example above. 
                                                                                                                                                                                                                      For example if we change  odd with even terms, i.e.   $S'_j=-1/2+1-1/4+1/3-1/6+ 1/5-1/8+ 1/7- . .. . .=- \ln(2)$.

                                                                                                                                                                                                                      \noindent Partial sums are  the same for even index of the new series  and  different  for  odd index of new series. So partial sums of  $S'_j$ and $S_j $ are different, but  they tend   always to same limit as  $J_{max} \rightarrow \infty$. Notice  however we have done a reordering on infinite terms. So we even can say:  not all , but, some reordering of infinite terms  can be done also on conditional convergent series  and the result  when index  $\rightarrow  \infty$ is preserved. But here we do only ``manipulations'' with { \bf indexes less than the  running index }. So observation does not apply.
                                                                                                                                                                                                                      

\item { \it The 
argument in equations ( \ref{RappSommePosGRH} ) and following  through 
 rests on the idea that because each term in
the series for  $\epsilon > 0$ is smaller in magnitude than the corresponding term for $\epsilon = 0$ (due
to the factor of $p^{- \epsilon}$), the overall sum must also be smaller in magnitude.
This reasoning is incorrect.} 

\noindent
ANSW. ( \ref{RappSommePosGRH} ) is a ratio not of the whole sum, but only of the positive contributions  at different $\epsilon$ values. Each positive contribution of $h-$primes in semi-oscillation $k$: $(h,k)$ , i.e. $ O^+_{k,h }(\epsilon)_{\int dLi}  + O^+_{k,h }(\epsilon)_{\sum} $ decreases with increasing $\epsilon$. 
The same holds for their sum (as are all positive terms). So  ratio ( \ref{RappSommePosGRH} ) also $ \forall p_{max}  <\infty$ is less than 1 and also  (\ref {RatioUnequalEps} ) seems justyfied for big $p_{max}<\infty$. 
The same holds for negative contributions. The rationale is in (\ref{ApproachFromBelow} ), where it is shown an approach from below  to $\rho(\epsilon=0)$, with suitable small $\epsilon \ \ \forall p_{max}$ big enough. Overall sum is not in question at   ( \ref{RappSommePosGRH} ). The factor $p_{max}^{- \epsilon}$ is  borrowed from    expression     (\ref {RatioPrototype}  )            , a reasonable prototype of the ratios with different $\epsilon$. The computation takes in consideration correct contributions. To assure  further that this ratio   approach 1, from below,   identity  (\ref{Id}) has been deeply analyzed. Of course this is not a  practical method to compute the involved quantities \footnote{as the overall sum for example}
, but only,  through corollary to Lemma 6, i.e. 
(\ref{L6ALContinuo}), a way  to get a speedy though  rough estimate for:
  $$
   \left[ \frac {\partial^2 \angle[L(s,\chi_{
    primitive}) )]  } { \partial \epsilon\partial t} \right]_{ \epsilon =0}  \ge 0
   $$

\item {\it Analyzing the sum of positive terms in isolation from the sum of negative (like in  (\ref {FinRes}  )
terms is not a valid operation 
}

\noindent ANSW Generally speaking I can agree but (see  example two bullets above)  $S_{J_{max}}= \sum^+_{J_{max}} -\sum^-_{J_{max}} $ tends exacly to the limit, provided the operation is carried on a finite terms set, with indexes $\le$ running index. So I do not understand the sentence:  it``is not a valid operation''. See also following bullet.
\item {\it  algebraic manipulations  from ( \ref{RappSommePosGRH} ) through (\ref {whereWePointed}) effectively
treat the series as if it were a sum of positive terms, ignoring the crucial role of
cancellation. This is a classic error when dealing with conditional convergence} 

\noindent ANSW.We have chosen to compute the single pieces of equi-sign contribution in absolute value and to put  correct sign  given by (\ref{Campate}) afterward. This is a choice of opportunity because we are interested in a form $\infty-\infty$  that, as explained three bullets   above, can be unusual, but, it is correct because we consider always { \bf indexes less than the  running index }. The criteria to get a correct sum is : `` are the partial sums changed?'' . The answer is negative, so  no mistake seems present.

\item { \it It is true that one can make rearrangements when there are finitely many terms,
but then one cannot take the infinite limit because 
rearrangements are not correct for infinite sums. Therefore, the argument in Lemma 6
and the conclusions drawn from it are invalid}. ANSW. 
  (see example four bullets above)  If we write $S_{J_{max}}= \sum^+_{J_{max}} -\sum^-_{J_{max}} $, we do not alter partial sum and so  final result, as we  use  { \bf indexes less than the  running index}. We transfom the limit in something like $\infty - \infty$ without changing result.
   The reason this form is interesting  is that  (\ref{RappSommePosGRH} ) and (\ref{RappSommeNegGRH} )  tend to be the same value  and with  suitable $\delta $ and $\epsilon$ we can write (\ref {ApproachFromBelow} ) proving the approach to $\rho(\epsilon=0)$ from below  that implies corollary to Lemma 6 i.e. 
(\ref{L6ALContinuo}) . . . .

\item{\it It seems that  in Lemma 4 and the subsequent analysis in Lemma 6,
sums over prime numbers  are replaced with integrals weighted by the logarithmic integral, citing the
Prime Number Theorem in Arithmetic Progressions (PNTAP) in equation (3.7). The
PNTAP provides an asymptotic relationship. It states that $\pi(x)_{h,q} \sim  Li(x)/ \phi(q)$ as $x \rightarrow  \infty$.
This does not mean the two are equal, and the error term is not negligible. Using this as a
direct substitution in an argument that depends on the precise cancellation of oscillating
terms is not rigorous. A valid argument would require the use of an explicit formula for
the distribution of primes, complete with error terms, and demonstrate that these error
terms do not affect the final result.
} 

\noindent ANSW In Lemma 4 it is used the function $Li(x)$, but, with whatsoever well behaved function  the result of (\ref{HDepInt} ) is null  because it depends  from (\ref {XiNonPrincip}). So in Lemma 4 there is no reference to (PNTAP). The rationale of (\ref{WithLogInt}) is that of an algebraic trick to insert, with zero effective effect,  the integral of PNT in (\ref{DerLAngleSuTL}). The difference  in (\ref{WithLogInt})  means nothing  because we are subtracting zero (i.e.(\ref{HDepInt} )) from  (\ref{DerLAngleSuTL}). But it is useful in Lemma 6 to build  { \bf not differences} between primes distribution $\pi(x)$ and $Li(x)$ but {\bf ratios} between functions involving $\pi(x)$ and $Li(x)$ (See Lemma 7). In ratios asymptotic relations can be used  at best. Of course not in differences. So I recommend  a more careful reading and I think that the correct expectations of the reader will be perhaps satisfied.

\noindent I thank  anyway a lot the reader  for this observation that highlight a basic and simple idea. Namely that  it is more convenient  to look for {\bf ratios} between functions involving  $\pi(x)$ and $Li(x)$  rather than to look for { \bf  differences}.

\item{\it The proposed methods are unconventional  } 
 \noindent ANSW  Perhaps you are right, but, my problem is to identify, and, if possible eliminate, eventual flaws.

\item{\it The paper  has not been validated by a rigorous peer review } 
 \noindent ANSW That is true, but before to ask for it ( that by the way it  is no simple matter) I hope to uncover possible problems that could invalidate the paper, and at same time I try to stabilize the treatment of the topic. To this aim readers observations are very important for me. After all, as far as I know,  these ideas have  ``circulated''  only since a few years. Besides I realize  my papers are like ``message in the bottle'' in the  amateur proposals ocean. The probability that a truly expert mathematician cares  them (also to reject them if they deserve it) is close to zero.

\item{\it  Computed results cannot be taken as proofs}
\noindent ANSW Yes I agree. Let us say they are a necessary condition not a sufficient one. 
In other words, for example,  it makes sense to say that fig. \ref{Q13X10Eps} is compatible, i.e verifies, Lemma 6. Of course we cannot say it prove it.       However  I think that the matching of the $t_{peak}$ of the peaks  of  (\ref{DerLAngleSuTL}) with the zeros  of $L(s,\chi)$,
and the matching of the floor with the computed one by 
(\ref {Asympt} )  and (\ref{DerAngleDeT}) (see on figg, 
\ref  {Simmetries} through \ref {Q9Even58M}) is an interesting fact that confirm the exposed theory.

\item{\it The attempt to link mathematical problems to other fields like physics is not inherently frivolous, but, the lack of peer-reviewed validation makes it a speculative effort rather than a serious approach }
\noindent ANSW  The reference to a kinematic quantity as angular momentum is only a short  way to describe  the method I explored. So I used it in the abstract  of \cite {Giovanni Lodone Dec2024} only to be synthetic. The importance of this formal  link (however mathematically  based) with  classical mechanics  in  the logic of the paper is  zero. Relating a serious  peer-review, I'll be glad to face it because I'll have certaily to learn from it. But, I'm afraid,  it is not a likely outcome.
See note two bullets above. 

\item {\it  The term ``spectrum of primes''  is not a rigorously defined concept} 
\noindent ANSW Perhaps you are right. I borrowed  this term from \cite {Mazur:2015}, and stressed it  in (\cite {Giovanni Lodone Dec2024}) . But  anyhow a formula that gather zeros of L-functions and primes ( as (\ref {DerLAngleSuTL}) ), seems interesting.

\end{itemize}

\end{document}